\documentclass[reqno]{amsart}
\usepackage{amsmath, amsthm, amssymb, amscd}

\usepackage[a4paper,width=160mm,top=25mm,bottom=25mm,centering]{geometry}

\usepackage{enumitem}
\usepackage{graphicx} % Required for inserting images
\usepackage[dvipsnames]{xcolor}
\usepackage{comment}
\usepackage{cases}
\usepackage{array}
\usepackage{multirow}
\usepackage{longtable}
\usepackage{microtype}

\makeatletter
\newsavebox{\@brx}
\newcommand{\llangle}[1][]{\savebox{\@brx}{\(\m@th{#1\langle}\)}%
  \mathopen{\copy\@brx\kern-0.5\wd\@brx\usebox{\@brx}}}
\newcommand{\rrangle}[1][]{\savebox{\@brx}{\(\m@th{#1\rangle}\)}%
  \mathclose{\copy\@brx\kern-0.5\wd\@brx\usebox{\@brx}}}
\makeatother

\usepackage{tikz}
\usepackage{pgf}
\usepackage{tikzscale}
\usepackage{tikz-cd}
\usetikzlibrary{decorations.pathreplacing,decorations.markings}
\usetikzlibrary{arrows, knots, external, positioning}

\makeatletter
\tikzset{column sep/.code=\def\pgfmatrixcolumnsep{\pgf@matrix@xscale*(#1)},
  row sep/.code   =\def\pgfmatrixrowsep{\pgf@matrix@yscale*(#1)},
  matrix xscale/.code=%
    \pgfmathsetmacro\pgf@matrix@xscale{\pgf@matrix@xscale*(#1)},
  matrix yscale/.code=%
    \pgfmathsetmacro\pgf@matrix@yscale{\pgf@matrix@yscale*(#1)},
  matrix scale/.style={/tikz/matrix xscale={#1},/tikz/matrix yscale={#1}}}
\def\pgf@matrix@xscale{1}
\def\pgf@matrix@yscale{1}
\makeatother

\usepackage{url}
\usepackage{caption} % to give a padding between a table and its caption

\newtheorem{theorem}{Theorem}[section]

\newtheorem{lemma}[theorem]{Lemma}
\newtheorem{proposition}[theorem]{Proposition}
\newtheorem{corollary}[theorem]{Corollary}

\theoremstyle{definition}

\newtheorem{remark}[theorem]{Remark}

\numberwithin{equation}{section}

\usepackage[bookmarksnumbered=true, colorlinks=true]{hyperref}
\AtBeginDocument{\hypersetup{citecolor=JungleGreen, linkcolor=JungleGreen, urlcolor=JungleGreen}}

\begin{document}

\title[Real Fr{\o}yshov invariants and lattice homology]{Real Fr{\o}yshov invariants and lattice homology}
% temporary

\author{Woohyeok Jo}
\address{HUN-REN R\'enyi Institute of Mathematics, H-1053 Budapest, Re\'altanoda utca 13--15, Hungary}
\email{jo.woohyeok@renyi.hu}

%\thanks{}
\date\today

\begin{abstract}
For a class of almost rational plumbing $4$-manifolds equipped with smooth involutions, we give a description of the Real Fr{\o}yshov invariants of their boundaries in terms of lattice homology. Our approach uses the $S^1$-stable equivalence between Seiberg--Witten Floer spectra and lattice spectra constructed by Dai--Sasahira--Stoffregen, and the $\operatorname{O}(2)$-equivariant versions of these spectra developed by Kang--Park--Taniguchi. As applications, we derive a formula for these invariants for every $\operatorname{spin}^c$ structure on the double branched covers of even torus knots $T_{2p,q}$, and determine the smooth nonorientable $4$-genus of an infinite family of torus knots $T_{4,q}$.
\end{abstract}

\maketitle

\section{Introduction}

Real Seiberg--Witten theory provides invariants for $3$- and $4$-manifolds equipped with smooth orientation-preserving involutions, and has a wide range of variants and applications \cite{Tian-Wang-2009, Nakamura-2013, Nakamura-2015, Kato-2022, Li-2022, Miyazawa-2023, Baraglia-Hekmati-2024, Konno-Miyazawa-Taniguchi-2024, Konno-Miyazawa-Taniguchi-2025, Kang-Park-Taniguchi-2025, Baraglia-2026, Jo-Tagliente-2026, Kang-Park-Taniguchi-2026, Xiao-2026}. Among such involutions, \emph{odd involutions}, whose fixed point sets are nonempty and have codimension two, are of particular interest, as they arise naturally as the deck transformations of double branched covers of knots in $S^3$, or double branched covers along smoothly embedded surfaces in $4$-manifolds. In this framework, Miyazawa proved the existence of an infinite family of exotic $\mathbb{RP}^2$'s in $S^4$ \cite{Miyazawa-2023}, and Kang--Park--Taniguchi proved that $(2n,1)$-cables of the figure-eight knot are not smoothly slice \cite{Kang-Park-Taniguchi-2025, Kang-Park-Taniguchi-2026}.

In the latter result, the invariants \[
\delta_R (Y,\mathfrak{s},\tau)\in \mathbb{Q},
\]
called \emph{Real Fr{\o}yshov invariants}, defined for a rational homology $3$-sphere $Y$ with an odd involution $\tau$ and a Real $\operatorname{spin}^c$ structure $\mathfrak{s}$,\footnote{For a rational homology sphere $Y$ with an odd involution $\tau$, a $\operatorname{spin}^c$ structure $\mathfrak{s}$ is said to be \emph{Real} if $\tau^*\mathfrak{s}\cong \overline{\mathfrak{s}}$; see Section \ref{subsec:O(2)_SWF} for details.} played an important role. These are Fr{\o}yshov-type invariants introduced by Konno--Miyazawa--Taniguchi \cite{Konno-Miyazawa-Taniguchi-2024}, and are extracted from the Real version of Manolescu's Seiberg--Witten Floer spectrum \cite{Manolescu-2003}; see Section \ref{subsec:O(2)_SWF} for details. Although these invariants are generally difficult to compute, Kang--Park--Taniguchi \cite{Kang-Park-Taniguchi-2026} developed a method for computing them in the Real spin case under suitable hypotheses, using the combinatorial description of Seiberg--Witten Floer spectra in terms of lattice homology due to Dai--Sasahira--Stoffregen \cite{Dai-Sasahira-Stoffregen-2026}. The first goal of this paper is to extend these computations to general Real $\operatorname{spin}^c$ structures, and to establish a formula for Real Fr{\o}yshov invariants in terms of lattice homology.

For a negative definite plumbing graph $\Gamma$, let $W_\Gamma$ denote the associated $4$-manifold and let $Y_\Gamma:=\partial W_\Gamma$. As in \cite{Dai-Sasahira-Stoffregen-2026, Kang-Park-Taniguchi-2026}, we focus on the case where $\Gamma$ is \emph{almost rational}: a negative definite plumbing tree is almost rational if it can be made rational by decreasing the weight of at most one vertex \cite{Nemethi-2005}. For example, every negative definite star-shaped graph is almost rational. The definition of a rational graph is recalled in Section \ref{subsec:carrying_paths}.

\begin{theorem}\label{thm:main}
Let $\Gamma$ be an almost rational graph, let $\tau:W_\Gamma \to W_\Gamma$ be a smooth orientation-preserving involution whose restriction $\tau|_{Y_\Gamma}:Y_\Gamma\to Y_\Gamma$ is odd, and let $\mathfrak{s}$ be a Real spin$^c$ structure for $(Y_\Gamma,\tau|_{Y_\Gamma})$. Suppose there is an almost $I$-invariant path \[
\gamma= (\mathfrak{t}_{-n},\dots,\mathfrak{t}_{-1}, \mathfrak{t}_1,\dots,\mathfrak{t}_n)
\]
which carries the lattice homology of $(\Gamma,\mathfrak{s})$.\footnote{Almost $I$-invariant paths and the carrying condition are defined in Section \ref{subsec:almost_I_invariant_paths}. These notions are introduced in \cite{Kang-Park-Taniguchi-2026} and \cite{Dai-Sasahira-Stoffregen-2026}, respectively.} Then \[ 
\delta_R(Y_\Gamma,\mathfrak{s},\tau|_{Y_{\Gamma}}) = \frac{c_1(\mathfrak{t}_{\pm 1})^2-\sigma(W_\Gamma)}{16}. 
\]
\end{theorem} 

Theorem \ref{thm:main} reduces the computation of $\delta_R$ to finding an almost $I$-invariant path carrying the relevant lattice homology. Our second main result gives a sufficient condition for the existence of such a path for every Real $\operatorname{spin}^c$ structure. To formulate this, let $S_v\subset W_\Gamma$ denote the vertex sphere corresponding to a vertex $v$ of $\Gamma$. We say that an involution $\tau:W_\Gamma\to W_\Gamma$ \emph{permutes the vertices of $\Gamma$} if there is a permutation $\sigma$ of the vertices of $\Gamma$ such that, for every vertex $v$, $\tau(S_v)=S_{\sigma(v)}$ and the restriction $\tau|_{S_v}:S_v\to S_{\sigma(v)}$ is orientation-preserving. 

\begin{theorem}\label{thm:main_existence}
Let $\Gamma$ be an almost rational graph, and let $\tau:W_\Gamma\to W_\Gamma$ be a smooth orientation-preserving involution whose restriction $\tau_Y:=\tau|_{Y_\Gamma}$ is odd. Suppose that $\tau$ permutes the vertices of $\Gamma$ and that there is a $\tau$-fixed vertex $v_0$ such that decreasing the weight of $v_0$ produces a rational graph. Then, for every Real spin$^c$ structure $\mathfrak{s}$ for $(Y_\Gamma,\tau_Y)$, there exists an almost $I$-invariant path carrying the lattice homology of $(\Gamma,\mathfrak{s})$.
\end{theorem}

Combining the proofs of Theorems \ref{thm:main} and \ref{thm:main_existence} also yields the following formula for $\delta_R$ in terms of characteristic covectors. It is a Real analog of the formulas for Heegaard Floer $d$-invariants in \cite[Corollary 1.5]{Ozsvath-Szabo-2003} and \cite[Theorem 8.3]{Nemethi-2005}.

\begin{corollary}\label{cor:Froyshov_formula}
Under the hypotheses of Theorem \textup{\ref{thm:main_existence}}, let $\mathfrak{s}$ be a Real spin$^c$ structure for $(Y_\Gamma,\tau_Y)$, and let $[k]\subset \operatorname{Char}(\Gamma)$ be the corresponding orbit.\footnote{There is a natural bijection between $\operatorname{spin}^c$ structures on $Y_\Gamma$ and characteristic orbits in $\operatorname{Char}(\Gamma)$; see Section \ref{subsec:lattice_homology}.} Then there exists a unique set $B_{\mathfrak{s}}$ of pairwise nonadjacent $\tau$-fixed vertices of $\Gamma$ for which the set \[
\mathcal{K}_{[k]}:=\left\{\ell\in[k]:\ell+\tau^*\ell=2\sum_{v\in B_{\mathfrak{s}}}\operatorname{PD}[S_v] \right\}
\]
is nonempty. Moreover, we have the following properties:  \begin{enumerate}[label=\textup{(\alph*)}, leftmargin=2em]
    \item The Real Fr{\o}yshov invariant is given by  \[
\delta_R(Y_\Gamma,\mathfrak{s},\tau_Y)= \max_{\ell\in\mathcal{K}_{[k]}} \frac{\ell^2+|\Gamma|}{16},
\]
where $|\Gamma|$ denotes the number of vertices of $\Gamma$.
\item Let $\Gamma_{\mathrm{fix}}$ denote the full subgraph of $\Gamma$ spanned by the $\tau$-fixed vertices. Then \[
b_{\mathfrak{s}}:=\sum_{v\in B_{\mathfrak{s}}}v
\]
is a Wu class of $\Gamma_{\mathrm{fix}}$. Consequently, if $\det Q_{\Gamma_{\mathrm{fix}}}$ is odd, then $B_{\mathfrak{s}}$ is the support of the unique Wu class of $\Gamma_{\mathrm{fix}}$ and is therefore independent of $\mathfrak{s}$. Furthermore, $\det Q_{\Gamma_{\mathrm{fix}}}$ is odd whenever $Y_\Gamma$ is a $\mathbb{Z}_2$-homology sphere.
\end{enumerate} \end{corollary}

Here, $Q_{\Gamma_{\mathrm{fix}}}$ denotes the intersection form associated to $\Gamma_{\mathrm{fix}}$. Wu classes can be easily computed by several methods \cite{Neumann-1980, Saveliev-2002, Stipsicz-2008}. Thus, when $\det Q_{\Gamma_{\mathrm{fix}}}$ is odd, Corollary \ref{cor:Froyshov_formula} provides a practical method for computing the Real Fr{\o}yshov invariants: one first determines the Wu class of $\Gamma_{\mathrm{fix}}$ and then maximizes $\ell^2$ over $\mathcal{K}_{[k]}$. 

\begin{remark}
    Let $(Y,\mathfrak{s})$ be a $\operatorname{spin}^c$ rational homology sphere, and suppose $(W,\mathfrak{t})$ is a negative definite $\operatorname{spin}^c$ filling of $(Y,\mathfrak{s})$. Then \[
    \frac{c_1(\mathfrak{t})^2-\sigma(W)}{4}\leq d(Y,\mathfrak{s}),
    \] 
    and the two sides are congruent modulo $2$ \cite{Ozsvath-Szabo-2003-02}. Therefore, under the hypotheses of Theorem \ref{thm:main}, we have \begin{equation}\label{eq:d_and_delta_R}
        \frac{d(Y_\Gamma,\mathfrak{s})}{4} -\delta_R(Y_\Gamma,\mathfrak{s},\tau|_{Y_\Gamma}) \in \frac{1}{2}\mathbb{Z}_{\geq 0}.
    \end{equation}
    Under the hypotheses of Theorem \ref{thm:main_existence}, \eqref{eq:d_and_delta_R} also follows from the formula of Corollary \ref{cor:Froyshov_formula}(a), together with \cite[Theorem 8.3]{Nemethi-2005}.
\end{remark}

We next apply Corollary \ref{cor:Froyshov_formula} to the nonorientable $4$-genera of knots. For a knot $K\subset S^3$, its smooth \emph{nonorientable $4$-genus} is defined by \[
\gamma_4(K):=\min_F b_1(F),
\]
where the minimum is taken over all smoothly embedded compact nonorientable surfaces $F\subset B^4$ with
$\partial F=K$. This invariant was first introduced in \cite{Murakami-Yasuhara-2000}, and has since been extensively studied \cite{Gilmer-Livingston-2011, Batson-2014, Ozsvath-Stipsicz-Szabo-2017, Golla-Marengon-2018, Jabuka-Kelly-2018, Jabuka-VanCott-2021, Ghanbarian-2022, Feller-Golla-2023, Hoste-Shanahan-VanCott-2023, BKST-2025, Fairchild-2025, Sinha-2025, FGMP-2026, Lee-2026, Lee-Sabloff-2026}. 

However, even for torus knots, the exact values of $\gamma_4$ remain unknown in general. For example, while $\gamma_4(T_{p,q})=1$ for $p\leq 3$, the case $p=4$ is only partially understood. The following results are known; see \cite{BKST-2025} and the references therein. \begin{enumerate}[label=\textup{(\arabic*)}]
\item $\gamma_4(T_{4,q})\leq2$.
\item If $q\equiv 5,7\pmod{8}$, then $\gamma_4(T_{4,q})=2$.
\item If $q\in\{1,3,9,11\}$, then $\gamma_4(T_{4,q})=1$.
\end{enumerate}
Also, when $q\equiv 1,3\pmod{8}$, the linking form obstruction gives the following necessary conditions for $T_{4,q}$ to bound a M\"obius band: \begin{itemize}
    \item If $q\equiv 1\pmod{8}$ and $\gamma_4(T_{4,q})=1$, then $q=x^2\pm 8y^2$ for some $x,y\in \mathbb{Z}$.
    \item If $q\equiv 3\pmod{8}$ and $\gamma_4(T_{4,q})=1$, then $q=x^2+2y^2$ for some $x,y\in \mathbb{Z}$.
\end{itemize}
These follow from the linking form obstruction in \cite{BKST-2025}; see Corollary \ref{cor:torus_linking_form} and Remark \ref{rmk:BKST} for details. Observe that $\gamma_4(T_{4,q})$ remains undetermined for infinitely many values of $q$, beginning with $q=17,19,25,27,33,41,\dots$. Furthermore, as pointed out in \cite{BKST-2025}, none of the other previously known lower bounds for $\gamma_4$ provide any additional information for these cases.

We use Real Fr{\o}yshov invariants to obtain further obstructions. We first compute these invariants for the double branched covers \[
Y_{p,q}:=\Sigma_2(S^3,T_{2p,q})=\Sigma(2,2p,q) \quad (p\geq 1, \, q>1 \text{ odd}, \, \gcd(p,q)=1).
\]
Every $\operatorname{spin}^c$ structure on $Y_{p,q}$ is Real with respect to the deck transformation $\tau_{p,q}$. In Section \ref{sec:even_torus_knots}, we obtain a formula for all the invariants \[
\delta_R(Y_{p,q},\mathfrak s,\tau_{p,q})
\]
which can be computed by a finite algorithm; see Proposition \ref{prop:torus_Froyshov_formula}. For the unique spin structure $\mathfrak{s}_0$ of $Y_{p,q}$, this recovers the formula \[
\delta_R(Y_{p,q},\mathfrak{s}_0,\tau_{p,q}) =-\frac{1}{2}\bar\mu(Y_{p,q})
\]
proved in \cite[Theorem 4.7]{Kang-Park-Taniguchi-2026}, where $\bar{\mu}$ denotes the Neumann--Siebenmann invariant \cite{Neumann-1980, Siebenmann-1980, Saveliev-2002}; see Proposition \ref{prop:spin_case}. When $p=2$, we derive the following closed formula; see Proposition \ref{prop:T4q_Froyshov_formula}.
\begin{theorem}\label{thm:intro_T4q_formula}
Write $q=4k+l\geq3$, where $l\in\{1,3\}$. For the labeling of $\operatorname{Spin}^c(Y_{2,q})$ defined in \eqref{eq:torus_labels}, we have \[
\delta_R(Y_{2,q},\mathfrak{s}_i,\tau_{2,q}) =\frac{(q-2i)^2}{8q}-\frac{3-l}{16} -\max\left\{0,\left\lceil\frac{k-i}{2}\right\rceil\right\} \quad \left(0\leq i\leq\frac{q-1}{2}\right),
\]
and $\delta_R(Y_{2,q},\mathfrak{s}_{-i},\tau_{2,q}) =\delta_R(Y_{2,q},\mathfrak{s}_i,\tau_{2,q})$.
\end{theorem}

This leads to the following result.

\begin{theorem}\label{thm:T4q_nonorientable_genus}
Let $q\geq 1$ be an integer with $q\equiv1\pmod{8}$. If $\gamma_4(T_{4,q})=1$, then $q=x^2+32y^2$ for some $x,y\in \mathbb{Z}$.
\end{theorem}

This determines $\gamma_4(T_{4,q})=2$ for $q=17,73,89,97,161,193,\dots$, which was previously unknown. In fact, Theorem \ref{thm:T4q_nonorientable_genus} determines $\gamma_4(T_{4,q})=2$ for an infinite family of $q$ where $\gamma_4(T_{4,q})$ was previously unknown. For example, it yields $\gamma_4(T_{4,q})=2$ for $q=17\cdot 25^r$ ($r\geq 0$); see Corollary \ref{cor:T4q_infinite_family}.

\subsection*{Outline of the proofs}
We first explain the proof of Theorem \ref{thm:main}. We regard $\operatorname{O}(2)$ as $S^1\rtimes \langle I\rangle$, where $I:=\operatorname{diag}(1,-1)\in \operatorname{O}(2)$. Following the constructions of Kang--Park--Taniguchi \cite{Kang-Park-Taniguchi-2026}, we consider: \begin{itemize}[leftmargin=2em]
    \item The $\operatorname{O}(2)$-equivariant path lattice spectrum $\mathcal{H}_{\operatorname{O}(2)}(\gamma,\mathfrak{s})$ associated to an almost $I$-invariant path $\gamma$.
    \item The $\operatorname{O}(2)$-equivariant Seiberg--Witten Floer spectrum $\operatorname{SWF}_{\operatorname{O}(2)}(Y_\Gamma,\mathfrak{s},\tau|_{Y_{\Gamma}})$.
\end{itemize}
Restricting the $\operatorname{O}(2)$-action on $\mathcal{H}_{\operatorname{O}(2)}(\gamma,\mathfrak{s})$ to the subgroup $S^1$ yields the lattice spectrum $\mathcal{H}(\Gamma,\mathfrak{s})$ constructed by Dai--Sasahira--Stoffregen \cite{Dai-Sasahira-Stoffregen-2026}. Similarly, restricting $\operatorname{SWF}_{\operatorname{O}(2)}(Y_\Gamma,\mathfrak{s},\tau|_{Y_{\Gamma}})$ to $S^1$ recovers Manolescu's Seiberg--Witten Floer spectrum $\operatorname{SWF}_{S^1}(Y_\Gamma,\mathfrak{s})$. Dai--Sasahira--Stoffregen constructed an $S^1$-stable equivalence \begin{equation}\label{eq:DSS_equivalence}
    \mathcal{T}:\mathcal{H}(\Gamma,\mathfrak{s}) \longrightarrow \operatorname{SWF}_{S^1}(Y_\Gamma,\mathfrak{s})
\end{equation}
by gluing relative Bauer--Furuta maps using homotopies supplied by adjunction relations; see Section \ref{subsec:spectra_comparison} for details. On the other hand, taking the $I$-fixed-point set of $\operatorname{SWF}_{\operatorname{O}(2)}(Y_\Gamma,\mathfrak{s},\tau|_{Y_\Gamma})$ gives the Real Seiberg--Witten Floer spectrum $\operatorname{SWF}_R(Y_\Gamma,\mathfrak{s},\tau|_{Y_\Gamma})$ of Konno--Miyazawa--Taniguchi \cite{Konno-Miyazawa-Taniguchi-2024}. The Real Fr{\o}yshov invariant $\delta_R(Y_\Gamma,\mathfrak{s},\tau|_{Y_\Gamma})$ can be obtained from the $\mathbb{Z}_2$-equivariant (reduced) Borel cohomology of $\operatorname{SWF}_R(Y_\Gamma,\mathfrak{s},\tau|_{Y_\Gamma})$. Therefore, to prove the theorem, we have to determine the $\mathbb{Z}_2$-equivariant Borel cohomology of the $I$-fixed-point set of $\operatorname{SWF}_{\operatorname{O}(2)}(Y_\Gamma,\mathfrak{s},\tau|_{Y_\Gamma})$.

Let us compare our approach with the argument of Kang--Park--Taniguchi. The path lattice spectrum $\mathcal{H}_{\operatorname{O}(2)}(\gamma,\mathfrak{s})$ is represented by a finite $\operatorname{O}(2)$-CW complex, which is built from two ``arms'' joined by a central cylinder $S^{m\mathbb{C}}\wedge[-1,1]_+$ for some integer $m\geq 0$; see Figure \ref{fig:schematic_description_2}. The action of $I\in\operatorname{O}(2)$ exchanges the two arms, while on the central cylinder it is given by $[z,t]\mapsto[\bar z,-t]$. In the case where $Y_\Gamma$ is the double branched cover of a knot $K\subset S^3$ and $\mathfrak{s}$ is the unique spin structure, Kang--Park--Taniguchi constructed an $\operatorname{O}(2)$-equivariant map \[
\mathcal{T}_{\operatorname{O}(2)}: \mathcal{H}_{\operatorname{O}(2)}(\gamma,\mathfrak{s}) \longrightarrow \operatorname{SWF}_{\operatorname{O}(2)}(Y_\Gamma,\mathfrak{s},\tau|_{Y_\Gamma}),
\]
such that restricting the $\operatorname{O}(2)$-action to $S^1$ recovers \eqref{eq:DSS_equivalence}. To do this, they first construct an $S^1$-stable equivalence $\mathcal{H}_{\operatorname{O}(2)}(\gamma,\mathfrak{s}) \to \operatorname{SWF}_{\operatorname{O}(2)}(Y_\Gamma,\mathfrak{s},\tau|_{Y_\Gamma})$ that is $\operatorname{O}(2)$-equivariant on the union of the two arms. Then they choose an oriented surface $S\subset B^4$ with $\partial S=K$, take the double branched cover of $B^4$ along $S$, and consider the associated relative Bauer--Furuta map for the unique spin structure. Using this Bauer--Furuta map together with the argument in \cite{Dai-Sasahira-Stoffregen-2026}, they obtain the desired $\operatorname{O}(2)$-equivariant map $\mathcal{T}_{\operatorname{O}(2)}$. Finally, they show that the restriction of $\mathcal{T}_{\operatorname{O}(2)}$ to the $I$-fixed-point sets induces an isomorphism on mod $2$ singular homology, and use this to compute $\delta_R(Y_\Gamma,\mathfrak{s},\tau|_{Y_\Gamma})$. 

For our proof of Theorem \ref{thm:main}, let $X$ denote the $\operatorname{O}(2)$-CW complex representing $\mathcal{H}_{\operatorname{O}(2)}(\gamma,\mathfrak{s})$, as described above, and let $A\subset X$ be the union of the two arms. For an $\operatorname{O}(2)$-space $F$ representing $\operatorname{SWF}_{\operatorname{O}(2)}(Y_\Gamma,\mathfrak{s},\tau|_{Y_\Gamma})$, we first follow Kang--Park--Taniguchi's argument to construct an $S^1$-stable equivalence \[
f:X\longrightarrow F
\]
such that $f|_A:A\to F$ is $\operatorname{O}(2)$-equivariant. Instead of modifying the map $f$ to make it $\operatorname{O}(2)$-equivariant on all of $X$, we focus on the mapping cone \[
\Theta:=\operatorname{Cone}(f|_A:A\longrightarrow F)
\]
with its natural $\operatorname{O}(2)$-action. The mapping cone $\Theta$ is $S^1$-stably equivalent to $X/A$, which is an $S^1$-representation sphere, and the $I$-fixed-point sets of $\Theta$ and $F$ agree, i.e., $\Theta^I=F^I$. We apply Smith theory and Borel's dimension formula (Section \ref{subsec:equivariant_algebraic_topology}) to the subgroup $\langle-1,I\rangle\cong\mathbb{Z}_2\times\mathbb{Z}_2$ of $\operatorname{O}(2)$ to determine the mod $2$ singular cohomology of $\Theta^I=F^I$. This leads to the formula for $\delta_R$ in Theorem \ref{thm:main}.

We now briefly explain the proofs of the other main results. The proof of Theorem \ref{thm:main_existence} is combinatorial and adapts N\'emethi's arguments \cite{Nemethi-2005} to our equivariant setting. To construct an almost $I$-invariant path $\gamma= (\mathfrak{t}_{-n},\dots,\mathfrak{t}_{-1}, \mathfrak{t}_1,\dots,\mathfrak{t}_n)$ which carries the lattice homology, the key additional step is to choose the central pair $(\mathfrak{t}_{-1},\mathfrak{t}_1)$ of $\operatorname{spin}^c$ structures that maximize the weight function in a certain sense (Section \ref{subsec:construction_of_central_pair}). Corollary \ref{cor:Froyshov_formula} is a direct consequence of the construction in the proof of Theorem \ref{thm:main_existence}, together with the formula of Theorem \ref{thm:main}. 

In Section \ref{sec:even_torus_knots}, we study even torus knots  $T_{2p,q}$. The double branched cover of $T_{2p,q}$ is the Brieskorn manifold $\Sigma(2,2p,q)$, which is the boundary of a negative definite star-shaped plumbing. We examine the $\operatorname{spin}^c$ structures and use Corollary \ref{cor:Froyshov_formula} to derive a formula for $\delta_R$ that can be evaluated by a finite algorithm; see Proposition \ref{prop:torus_Froyshov_formula}. We use this formula to obtain a closed formula for $\delta_R$ in the case $p=2$, proving Theorem \ref{thm:intro_T4q_formula}; see Proposition \ref{prop:T4q_Froyshov_formula}.

Finally, to prove Theorem \ref{thm:T4q_nonorientable_genus}, we consider the double cover of $B^4$ branched along a M\"obius band bounded by $T_{4,q}$. This cover has $b_2=1$. We first examine the linking form of its boundary, which gives constraints on $q$. (This was also considered in \cite{BKST-2025}.) We then apply the Real Fr{\o}yshov inequality (Corollary \ref{cor:real_froyshov_filling}) to obtain further constraints, thereby proving Theorem \ref{thm:T4q_nonorientable_genus}.

\subsection*{Acknowledgments} 
The author would like to thank Alessio di Prisa, Marco Golla, Marco Marengon, Simone Tagliente, and Masaki Taniguchi for helpful conversations and comments. The author was supported by ERC Advanced Grant KnotSurf4d. No AI tools were used in this project.

\section{Seiberg--Witten Floer spectra}\label{sec:SWF}

We recall the constructions of Seiberg--Witten Floer spectra, Real Fr{\o}yshov invariants, and relative Bauer--Furuta maps. General references are \cite{Manolescu-2003, Manolescu-2007, Khandhawit-2015, Konno-Miyazawa-Taniguchi-2024, Kang-Park-Taniguchi-2026, Dai-Sasahira-Stoffregen-2026}.

\subsection{Representations}\label{subsec:representations}
We fix notation for the representations of $\mathbb{Z}_2$, $S^1$, and $\operatorname{O}(2)$ used below. Put $I:=\operatorname{diag}(1,-1)$ so that $I^2=1$, $IzI^{-1}=\bar{z}$ for $z\in S^1$, and $\operatorname{O}(2)=S^1\rtimes \langle I\rangle$. 

For $\mathbb{Z}_2$, we use the following representations: \begin{itemize}[leftmargin=2em]
    \item $\mathbb{R}$ denotes the trivial real $1$-dimensional representation.
    \item $\widetilde{\mathbb{R}}$ denotes the real $1$-dimensional representation where $-1$ acts by $x\mapsto -x$.
\end{itemize}
For $S^1$, we use: \begin{itemize}[leftmargin=2em]
    \item $\mathbb{R}$ denotes the trivial real $1$-dimensional representation.
    \item $\mathbb{C}$ denotes the real $2$-dimensional representation where $S^1$ acts by scalar multiplication.
\end{itemize}
For $\operatorname{O}(2)$, we use: \begin{itemize}[leftmargin=2em]
    \item $\mathbb{R}$ denotes the trivial real $1$-dimensional representation.
    \item $\widetilde{\mathbb{R}}$ denotes the real $1$-dimensional representation where $S^1$ acts trivially and $I$ acts by $x\mapsto -x$. 
    \item $\mathbb{C}$ denotes the real $2$-dimensional representation where $S^1$ acts by scalar multiplication and $I$ acts by conjugation $z\mapsto \bar{z}$. 
\end{itemize}

\subsection{Equivariant stable homotopy categories}
We introduce the equivariant stable homotopy categories $\mathfrak{C}_G$ used throughout this section. We will consider $G\in \{\mathbb{Z}_2,S^1,\operatorname{O}(2)\}$. The category $\mathfrak{C}_{S^1}$ was first introduced in \cite{Manolescu-2003} as an equivariant analog of the Spanier--Whitehead category. The categories $\mathfrak{C}_{\mathbb{Z}_2}$ and $\mathfrak{C}_{\operatorname{O}(2)}$ were introduced in \cite{Konno-Miyazawa-Taniguchi-2024} and \cite{Kang-Park-Taniguchi-2026}, respectively. 

We work with based $G$-spaces and based $G$-maps, and homotopies of based $G$-maps are assumed to be based and $G$-equivariant. We assume that the $G$-action fixes the basepoint. For an unbased $G$-space $T$, we put $T_+:=T\amalg \{*\}$ with basepoint $*$. For based $G$-spaces $X$ and $Y$, we denote by $[X,Y]_*^G$ the set of homotopy classes of based $G$-maps $X\to Y$. For a finite-dimensional real $G$-representation $V$, we write $S^V:=V\cup \{\infty\}$ for the one-point compactification with basepoint $\infty$, and define $\Sigma^VX:=S^V\wedge X$ for a based $G$-space $X$. 

For $G=\mathbb{Z}_2$, an object of the category $\mathfrak{C}_{\mathbb{Z}_2}$ is a triple $(X,m,n)$, where $X$ is a based $\mathbb{Z}_2$-space having the homotopy type of a based $\mathbb{Z}_2$-CW complex, $m\in \mathbb{Z}$, and $n\in \mathbb{Q}$. We interpret this triple as the formal desuspension $\Sigma^{-m\mathbb{R}-n\widetilde{\mathbb{R}}} X$. For objects $(X,m,n)$ and $(X',m',n')$, if $n-n'\in \mathbb{Z}$, the set of morphisms is defined by \[
\operatorname{Mor}_{\mathfrak{C}_{\mathbb{Z}_2}}((X,m,n),(X',m',n')):=\underset{p,q\to \infty}{\operatorname{colim}} \left[ \Sigma^{p\mathbb R\oplus q\widetilde{\mathbb{R}}}X,\, \Sigma^{(p+m-m')\mathbb{R}\oplus (q+n-n')\widetilde{\mathbb{R}}}X' \right]^{\mathbb{Z}_2}_*.
\]
If $n-n'\notin \mathbb{Z}$, the morphism set is empty. 

For $G=S^1$, the category $\mathfrak{C}_{S^1}$ is defined in the same way, using the representation $\mathbb{C}$ instead of $\widetilde{\mathbb{R}}$. Thus, an object $(X,m,n)$ in $\mathfrak{C}_{S^1}$ is interpreted as the formal desuspension $\Sigma^{-m\mathbb{R}-n\mathbb{C}}X$. 

For $G=\operatorname{O}(2)$, an object of $\mathfrak{C}_{\operatorname{O}(2)}$ is a quadruple $(X,a,b,c)$, where $X$ has the homotopy type of a based $\operatorname{O}(2)$-CW complex, $a,b\in \mathbb{Z}$, and $c\in \mathbb{Q}$. This object is interpreted as $\Sigma^{-a\mathbb{R}-b\widetilde{\mathbb{R}}-c\mathbb{C}}X$. When $c-c'\in\mathbb{Z}$, the morphism set is defined by \[
\operatorname{Mor}_{\mathfrak{C}_{\operatorname{O}(2)}} ((X,a,b,c),(X',a',b',c')):=\underset{p,q,r\to \infty}{\operatorname{colim}} \left[\Sigma^{p\mathbb{R}\oplus q\widetilde{\mathbb{R}}\oplus r\mathbb{C}}X,\, \Sigma^{(p+a-a')\mathbb{R}\oplus (q+b-b')\widetilde{\mathbb{R}}\oplus (r+c-c')\mathbb{C}}X' \right]^{\operatorname{O}(2)}_*,
\]
and the morphism set is empty if $c-c'\notin\mathbb{Z}$.

For any $G\in \{\mathbb{Z}_2,S^1,\operatorname{O}(2)\}$, an isomorphism in $\mathfrak{C}_G$ is called a \emph{$G$-stable equivalence}, and we write $\simeq_{G,\mathrm{st}}$ for this relation. For a based $G$-CW complex $X$, we identify it with the object $(X,0,0)$ (or $(X,0,0,0)$ if $G=\operatorname{O}(2)$) in $\mathfrak{C}_G$. 

For a $\mathbb{Z}_2$-representation $V$ isomorphic to $\mathbb{R}^p\oplus \widetilde{\mathbb{R}}^q$, we define \[
\Sigma^V (X,m,n):=(\Sigma^VX,m,n), \quad \Sigma^{-V}(X,m,n):=(\Sigma^VX,m+2p, n+2q)
\]
in $\mathfrak{C}_{\mathbb{Z}_2}$. Similarly, for an $S^1$-representation $V$ isomorphic to $\mathbb{R}^p\oplus \mathbb{C}^q$, define \begin{equation}\label{eq:suspension_of_object}
    \Sigma^V (X,m,n):=(\Sigma^VX,m,n), \quad \Sigma^{-V}(X,m,n):=(\Sigma^{V^{S^1}}X,m+2p, n+q)
\end{equation}
in $\mathfrak{C}_{S^1}$, where $V^{S^1}$ is the linear subspace of $V$ fixed by $S^1$. For an $\operatorname{O}(2)$-representation $V$ isomorphic to $\mathbb{R}^p\oplus \widetilde{\mathbb{R}}^q\oplus \mathbb{C}^r$, define \[
\Sigma^V (X,a,b,c):=(\Sigma^VX,a,b,c), \quad \Sigma^{-V}(X,a,b,c):=(\Sigma^V X,a+2p,b+2q,c+2r)
\]
in $\mathfrak{C}_{\operatorname{O}(2)}$. An explanation for these conventions is given in \cite[Section 3.1]{Dai-Sasahira-Stoffregen-2026}.

We have a restriction functor $\operatorname{Res}^{\operatorname{O}(2)}_{S^1}:\mathfrak{C}_{\operatorname{O}(2)}\to \mathfrak{C}_{S^1}$ given by \begin{equation}\label{eq:restriction_functor}
    (X,a,b,c)\longmapsto (X|_{S^1},a+b,c),
\end{equation}
where $X|_{S^1}$ denotes $X$ with the restricted $S^1$-action. 

Observe that $-1\in S^1\subset \operatorname{O}(2)$ commutes with $I\in \operatorname{O}(2)$. Therefore, if $X$ is a based $\operatorname{O}(2)$-space, then $\mathbb{Z}_2=\{\pm1\}\subset S^1$ acts on the $I$-fixed-point set $X^I$. For the $\operatorname{O}(2)$-representations $\mathbb{R}$, $\widetilde{\mathbb{R}}$, and $\mathbb{C}$, we have \[
\mathbb{R}^I=\mathbb{R}, \quad \widetilde{\mathbb{R}}^I=0, \quad \mathbb{C}^I=\widetilde{\mathbb{R}}
\]
as $\mathbb{Z}_2$-representations. Following \cite{Kang-Park-Taniguchi-2026}, we define a functor $\Phi^I:\mathfrak{C}_{\operatorname{O}(2)}\to \mathfrak{C}_{\mathbb{Z}_2}$ by \begin{equation}\label{eq:mod_2_restriction}
    (X,a,b,c)\longmapsto (X^I,a,c).
\end{equation}

For $q\in \mathbb{Q}$, we use the notation \begin{equation}\label{eq:rational_sphere}
    S^{q\mathbb{C}}:= (S^0,0,-q)
\end{equation}
in the category $\mathfrak{C}_{S^1}$. More generally, for any object $(X,m,n)$, we write \[
\Sigma^{q\mathbb{C}}(X,m,n):=(X,m,n-q).
\]
If $q$ is an integer, this convention agrees with \eqref{eq:suspension_of_object} up to $S^1$-stable equivalence. 

\subsection{Spaces of type SWF}\label{subsec:type_SWF}
Let $s\in \mathbb{Z}_{\geq 0}$ be a nonnegative integer, and let $X$ be a based $\mathbb{Z}_2$-space having the homotopy type of a finite based $\mathbb{Z}_2$-CW complex. We say that $X$ is a space of \emph{type $\mathbb{Z}_2$-$\operatorname{SWF}$ at level} $s$ if the fixed-point set $X^{\mathbb{Z}_2}$ is $\mathbb{Z}_2$-homotopy equivalent to the sphere $S^{s\mathbb{R}}$. 

Let $X$ be a space of type $\mathbb{Z}_2$-SWF at level $s$, and consider its $\mathbb{Z}_2$-equivariant (reduced) Borel cohomology $\widetilde{H}^*_{\mathbb{Z}_2}(X;\mathbb{F}_2)$, where $\mathbb{F}_2$ is the field of two elements. Since $\widetilde{H}^*_{\mathbb{Z}_2}(S^0;\mathbb{F}_2)\cong \mathbb{F}_2[W]$ with $|W|=1$, $\widetilde{H}^*_{\mathbb{Z}_2}(X;\mathbb{F}_2)$ has a graded $\mathbb{F}_2[W]$-module structure. Define \[
d(X):=\min\left\{m\in\mathbb{Z}:\text{there exists $x\in \widetilde{H}^m_{\mathbb{Z}_2}(X;\mathbb{F}_2)$ such that $W^kx\neq 0$ for all $k\geq 0$} \right\}.
\]
The equivariant localization theorem \cite[Theorem III.3.8]{tom_Dieck-1987} implies that $d(X)$ is a well-defined finite integer. 

For any finite-dimensional $\mathbb{Z}_2$-representation $V$, the suspension $\Sigma^V X$ is also of type $\mathbb{Z}_2$-SWF, and equivariant Thom isomorphism gives \[
d(\Sigma^V X)=d(X)+\dim_{\mathbb{R}}V.
\]
For an object $(X,m,n)$ in $\mathfrak{C}_{\mathbb{Z}_2}$, where $X$ is of type  $\mathbb{Z}_2$-$\operatorname{SWF}$, we define \begin{equation}\label{eq:spectrum_cohomology}
    \widetilde{H}^*_{\mathbb{Z}_2}((X,m,n);\mathbb{F}_2):= \widetilde{H}^{*+m+n}_{\mathbb{Z}_2}(X;\mathbb{F}_2), \quad d(X,m,n):=d(X)-m-n.
\end{equation}
These depend only on the $\mathbb{Z}_2$-stable equivalence class of $(X,m,n)$.

\subsection{Seiberg--Witten Floer spectra}\label{subsec:SWF}

We briefly recall the construction of $S^1$-equivariant Seiberg--Witten Floer spectra for rational homology spheres; see \cite{Manolescu-2003} for details. 

Let $Y$ be an oriented rational homology $3$-sphere with a Riemannian metric $g$ and a $\operatorname{spin}^c$ structure $\mathfrak{s}$. Let $\mathbb{S}$ denote the spinor bundle, and let $\rho:T^*Y\to \operatorname{End}_{\mathbb{C}}(\mathbb{S})$ denote Clifford multiplication. Its extension $\rho:T^*Y\otimes \mathbb{C}\to \operatorname{End}_{\mathbb{C}}(\mathbb{S})$ identifies $T^*Y\otimes i\mathbb{R}$ with $i\cdot \mathfrak{su}(\mathbb{S})$, the bundle of trace-free Hermitian endomorphisms of $\mathbb{S}$. Fix a reference $\operatorname{spin}^c$ connection $A_0$ whose induced connection on the determinant line bundle is flat, and let $D_{A_0}:\Gamma(\mathbb{S})\to \Gamma(\mathbb{S})$ denote the corresponding Dirac operator. 

Consider the configuration space \[
\mathcal{C}:=\Omega^1(Y,i\mathbb{R})\oplus \Gamma(\mathbb{S}).
\]
The gauge group $\mathcal{G}:=C^\infty(Y,S^1)$ acts on $\mathcal{C}$ by $u\cdot (a,\phi)=(a-u^{-1}du,u\phi)$. The normalized gauge group $\mathcal{G}_0$ is the subgroup of $\mathcal{G}$ consisting of elements of the form $e^{\xi}$ with $\xi\in C^\infty(Y,i\mathbb{R})$ and $\int_Y \xi\operatorname{dvol}=0$. The global Coulomb slice is given by \[
\mathcal{V}:= \ker d^*\oplus \Gamma(\mathbb{S})\subset \mathcal{C}.
\]
For any $(a,\phi)\in \mathcal{C}$, there exists a unique $u\in \mathcal{G}_0$ such that $u\cdot (a,\phi)\in \mathcal{V}$. The assignment $(a,\phi)\mapsto u\cdot (a,\phi)$ defines a (nonlinear) projection $\Pi:\mathcal{C}\to \mathcal{V}$. There is a restricted action on $\mathcal{V}$ of the group $S^1$ of constant gauge transformations, given by $z\cdot (a,\phi)=(a,z\phi)$.

The Chern--Simons--Dirac functional $\operatorname{CSD}:\mathcal{C}\to \mathbb{R}$ is defined by \[
\operatorname{CSD}(a,\phi)=-\frac{1}{2}\int_Y a\wedge da +\frac{1}{2}\int_Y \langle \phi, D_{A_0} \phi+\rho(a)\phi\rangle \operatorname{dvol}.
\]
It is gauge invariant. With respect to the $L^2$-metric, the formal gradient of $\operatorname{CSD}$ is given by \[
\nabla_{L^2}\operatorname{CSD}(a,\phi)=(* da+q(\phi),\, D_{A_0}\phi+\rho(a)\phi),
\]
where $q(\phi):=\rho^{-1}((\phi\otimes \phi^*)_0)\in \Omega^1(Y,i\mathbb{R})$.

For $(a,\phi)\in \mathcal{V}$, we have \begin{equation}\label{eq:decomposition_l+c}
    d\Pi_{(a,\phi)}\nabla_{L^2}\operatorname{CSD}(a,\phi)=l(a,\phi)+c(a,\phi),
\end{equation}
where $l:=*d\oplus D_{A_0}:\mathcal{V}\to \mathcal{V}$ is a linear, $L^2$-self-adjoint operator, and $c:\mathcal{V}\to \mathcal{V}$ is the nonlinear part given by \[
c(a,\phi)=(q(\phi)-d\xi(q(\phi)), \,\rho(a)\phi+\xi(q(\phi))\phi).
\] 
Here $\xi:\Omega^1(Y,i\mathbb{R})\to \Omega^0(Y,i\mathbb{R})$ is the linear operator characterized by \[
d^*d\xi(a)=d^*a, \quad \int_Y \xi(a)\operatorname{dvol}=0 \quad (a\in \Omega^1(Y,i\mathbb{R})). 
\]
Both $l$ and $c$ are $S^1$-equivariant. 

For $\lambda,\mu\in \mathbb{R}$ with $\lambda<\mu$, let $\tilde{p}^\mu_\lambda$ denote the $L^2$-orthogonal projection from $\mathcal{V}$ onto the finite-dimensional subspace $\mathcal{V}^\mu_\lambda$ spanned by the eigenvectors of $l$ with eigenvalues in $(\lambda,\mu]$. Choose a smooth function $\beta:\mathbb{R}\to \mathbb{R}$ with $\beta \geq 0$, $\beta^{-1}((0,\infty))=(0,1)$, and $\int_{\mathbb{R}} \beta=1$. Then define $p^\mu_\lambda:\mathcal{V}\to \mathcal{V}$ by \[
p^\mu_\lambda =\int_0^1 \beta(\tau)\tilde{p}^{\mu-\tau}_{\lambda+\tau} \,d\tau \quad (\mu,-\lambda>1).
\]

Fix an integer $k\geq 4$ and let $\mathcal{V}_{k+1}$ denote the Sobolev completion with respect to the $L^2_{k+1}$-norm. There exists $R>0$ such that, for sufficiently large $\mu, -\lambda$, if $x:\mathbb{R}\to \mathcal{V}^\mu_\lambda$ satisfies \[
\frac{dx}{dt}=-(l+p^\mu_\lambda c)(x(t)), 
\]
and if the image of $x$ is contained in the closed ball $\overline{B}(2R)$ in $\mathcal{V}_{k+1}$ of radius $2R$, then it is in fact contained in the open ball $B(R)$ in $\mathcal{V}_{k+1}$ of radius $R$. 

Let $u:\mathcal{V}^\mu_\lambda \to [0,1]$ be a smooth $S^1$-invariant function that vanishes outside $B(3R)\cap \mathcal{V}^\mu_\lambda$ and is identically $1$ on $\overline{B}(2R)\cap \mathcal{V}^\mu_\lambda$. The vector field $x\mapsto -u(x)(l+p^\mu_\lambda c)(x)$ is $S^1$-equivariant and compactly supported, so it generates a complete smooth $S^1$-equivariant flow $\varphi:\mathbb{R}\times \mathcal{V}^\mu_\lambda\to \mathcal{V}^\mu_\lambda$. The ball $\overline{B}(2R)\cap \mathcal{V}^\mu_\lambda$ is an $S^1$-invariant isolating neighborhood of the flow $\varphi$, so there exists an $S^1$-equivariant index pair $(N^\mu_\lambda,L^\mu_\lambda)$ for its invariant set. Such a pair can be chosen so that the $S^1$-space $N^\mu_\lambda/L^\mu_\lambda$ with basepoint $L^\mu_\lambda/L^\mu_\lambda$ has the homotopy type of a finite $S^1$-CW complex. 

The \emph{$S^1$-equivariant Seiberg--Witten Floer spectrum} of $(Y,\mathfrak{s})$ is defined as the object \[
\operatorname{SWF}_{S^1}(Y,\mathfrak{s}):= \Sigma^{-\mathcal{V}_\lambda^0}\left(N^\mu_\lambda/L^\mu_\lambda,\,0,\,n(Y,\mathfrak{s},g)\right)
\]
in $\mathfrak C_{S^1}$, for sufficiently large $\mu$ and $-\lambda$. The number $n(Y,\mathfrak{s},g)$ is defined in \cite{Manolescu-2003} by \begin{equation}\label{eq:correction_term}
    n(Y,\mathfrak{s},g):=\operatorname{ind}_{\mathbb{C}}D^+_{X,\mathrm{APS}} - \frac{c_1(\mathfrak{t})^2-\sigma(X)}8 \in\mathbb{Q},
\end{equation}
where: \begin{itemize}[leftmargin=2em]
    \item $(X,\mathfrak{t})$ is a compact $\operatorname{spin}^c$ $4$-manifold with boundary $(Y,\mathfrak{s})$, with a Riemannian metric that is isometric to $[0,1]\times Y$ on a boundary collar. 
    \item $\operatorname{ind}_{\mathbb{C}}D_{X,\operatorname{APS}}^+$ is the complex index of the Dirac operator with APS boundary conditions, for a $\operatorname{spin}^c$ connection $\widehat{A}$ on $X$ that is pulled back from the connection $A_0$ on the collar. 
\end{itemize}
The number $n(Y,\mathfrak{s},g)$ is independent of the choice of $X$, and the isomorphism class of $\operatorname{SWF}_{S^1}(Y,\mathfrak{s})$ in $\mathfrak{C}_{S^1}$ is an invariant of $(Y,\mathfrak{s})$. 

Write \begin{equation}\label{eq:decomposition_eigenspace}
    \mathcal{V}_\lambda^0 = V_\lambda^0\oplus W_\lambda^0 \subset \ker d^*\oplus \Gamma(\mathbb{S}),
\end{equation}
so that $V_\lambda^0$ is a trivial real $S^1$-representation, and $W_\lambda^0$ is a complex $S^1$-representation with the scalar $S^1$-action. Then we have an $S^1$-stable equivalence  \[
\operatorname{SWF}_{S^1}(Y,\mathfrak{s}) \simeq_{S^1,\mathrm{st}} \left(N^\mu_\lambda/L^\mu_\lambda,\, \dim_{\mathbb{R}}V_\lambda^0,\, \dim_{\mathbb{C}}W_\lambda^0+n(Y,\mathfrak{s},g) \right).
\]

\subsection{O(2)-equivariant and Real Seiberg--Witten Floer spectra}\label{subsec:O(2)_SWF}
We now recall the construction of the $\operatorname{O}(2)$-equivariant and Real Seiberg--Witten Floer spectra; see \cite{Konno-Miyazawa-Taniguchi-2024, Kang-Park-Taniguchi-2026} for details.

Let $Y$ be a rational homology $3$-sphere with an \emph{odd involution} $\tau:Y\to Y$; this means that $\tau$ is an orientation-preserving smooth involution with nonempty fixed-point set of codimension $2$. Suppose $Y$ is equipped with a $\tau$-invariant metric $g$ and a $\operatorname{spin}^c$ structure $\mathfrak{s}$ with $\tau^*\mathfrak{s}\cong \overline{\mathfrak{s}}$. Let $\mathbb{S}$ and $\rho$ denote the spinor bundle and Clifford multiplication, respectively, as in Section \ref{subsec:SWF}. The involution $\tau$ lifts to an anti-linear fiberwise isometry $\mathcal{I}:\mathbb{S}\to\mathbb{S}$ with $\mathcal{I}^2=\operatorname{id}$, compatible with Clifford multiplication in the sense that $\mathcal{I}(\rho(a)\phi)=\rho(\tau^*a)\mathcal{I}(\phi)$ for $a\in \Omega^1(Y,\mathbb{R})$ and $\phi\in \Gamma(\mathbb{S})$. The pair $(\mathfrak{s},\mathcal{I})$ is called a \emph{Real spin$^c$ structure} for $(Y,\tau)$; we usually suppress the lift $\mathcal{I}$ from the notation and say that $\mathfrak{s}$ is a Real $\operatorname{spin}^c$ structure for $(Y,\tau)$.\footnote{On a rational homology sphere, such lifts are unique up to gauge equivalence. We fix a lift throughout the constructions.} 

Choose an $\mathcal{I}$-compatible $\operatorname{spin}^c$ connection $A_0$ whose induced connection on the determinant line bundle is flat. The element $I\in\operatorname{O}(2)$ (see Section \ref{subsec:representations}) defines an involution on $\mathcal{C}$ by \begin{equation}\label{eq:involution_on_Coulomb_slice}
    I\cdot(a,\phi)=(-\tau^*a,\mathcal{I}\phi).
\end{equation}
Together with the scalar $S^1$ action $z\cdot (a,\phi)=(a,z\phi)$, this gives an action of $\operatorname{O}(2)=S^1\rtimes\langle I\rangle$ on $\mathcal{C}$, and this action preserves the global Coulomb slice $\mathcal{V}$. The functional $\operatorname{CSD}:\mathcal{C}\to\mathbb{R}$ is $\operatorname{O}(2)$-invariant, and the maps $l,c:\mathcal{V}\to\mathcal{V}$ are $\operatorname{O}(2)$-equivariant. Taking the cutoff function $u$ to be $\operatorname{O}(2)$-invariant, the construction of Section \ref{subsec:SWF} gives an $\operatorname{O}(2)$-equivariant flow. Choose an $\operatorname{O}(2)$-equivariant index pair $(N^\mu_\lambda,L^\mu_\lambda)$ for its isolated invariant set such that $N^\mu_\lambda/L^\mu_\lambda$ has the homotopy type of a finite based $\operatorname{O}(2)$-CW complex. The \emph{$\operatorname{O}(2)$-equivariant Seiberg--Witten Floer spectrum} is defined as the object \[
\operatorname{SWF}_{\operatorname{O}(2)}(Y,\mathfrak{s},\tau):= \Sigma^{-\mathcal{V}_\lambda^0} \left(N^\mu_\lambda/L^\mu_\lambda,\,0,\,0,\,n(Y,\mathfrak{s},g)\right)
\]
in $\mathfrak{C}_{\operatorname{O}(2)}$, for sufficiently large $\mu$ and $-\lambda$. Its isomorphism class in $\mathfrak{C}_{\operatorname{O}(2)}$ is an invariant of $(Y,\mathfrak{s},\tau)$.  

If we decompose $\mathcal{V}^0_\lambda$ as in \eqref{eq:decomposition_eigenspace}, then there exist $a_+,a_-,b\in \mathbb{Z}_{\geq 0}$ such that \[
V_\lambda^0 \cong \mathbb{R}^{a_+}\oplus\widetilde{\mathbb{R}}^{a_-}, \quad W_\lambda^0\cong\mathbb{C}^b
\]
as $\operatorname{O}(2)$-representations, and we have an $\operatorname{O}(2)$-stable equivalence \begin{equation}\label{eq:O(2)-stable_equivalence}
   \operatorname{SWF}_{\operatorname{O}(2)}(Y,\mathfrak{s},\tau)\simeq_{\operatorname{O}(2),\mathrm{st}} \left(N^\mu_\lambda/L^\mu_\lambda,\,a_+,\,a_-,\,b+n(Y,\mathfrak{s},g)\right). 
\end{equation}
For the restriction \eqref{eq:restriction_functor}, we clearly have \[
\operatorname{Res}^{\operatorname{O}(2)}_{S^1} \operatorname{SWF}_{\operatorname{O}(2)} (Y,\mathfrak{s},\tau) \simeq_{S^1,\mathrm{st}} \operatorname{SWF}_{S^1}(Y,\mathfrak{s}).
\]

Now recall the functor $\Phi^I:\mathfrak{C}_{\operatorname{O}(2)}\to \mathfrak{C}_{\mathbb{Z}_2}$ given in \eqref{eq:mod_2_restriction}, which takes $I$-fixed-point sets. The \emph{Real Seiberg--Witten Floer spectrum} is defined as the object \[
\operatorname{SWF}_R(Y,\mathfrak{s},\tau):= \Phi^I\operatorname{SWF}_{\operatorname{O}(2)} (Y,\mathfrak{s},\tau)
\]
in $\mathfrak{C}_{\mathbb{Z}_2}$, which is $\mathbb{Z}_2$-stably equivalent to the object \[
\left( (N^\mu_\lambda/L^\mu_\lambda)^I,\,a_+,\,b+n(Y,\mathfrak{s},g)\right).
\]
The following lemma can be proved by the same argument as in the proof of  \cite[Lemma 3.13]{Dai-Sasahira-Stoffregen-2026}.
\begin{lemma}\label{lem:O2_SWF_fixed_points}
The $S^1$-action on $N^\mu_\lambda/L^\mu_\lambda$ is free outside $(N^\mu_\lambda/L^\mu_\lambda)^{S^1}$, and there is a based $\langle I\rangle$-homotopy equivalence \[
(N^\mu_\lambda/L^\mu_\lambda)^{S^1} \simeq_{\langle I\rangle} S^{V_\lambda^0} \cong_{\langle I\rangle} S^{a_+\mathbb{R}\oplus a_-\widetilde{\mathbb{R}}}.
\]
Consequently, the space $(N^\mu_\lambda/L^\mu_\lambda)^I$ is of type $\mathbb{Z}_2$-$\operatorname{SWF}$ at level $a_+$.
\end{lemma}
Using \eqref{eq:spectrum_cohomology}, the \emph{Real Seiberg--Witten Floer cohomology} is defined as the $\mathbb{Q}$-graded $\mathbb{F}_2[W]$-module \[
\operatorname{SWFH}_R^*(Y,\mathfrak{s},\tau) := \widetilde{H}_{\mathbb{Z}_2}^*(\operatorname{SWF}_R(Y,\mathfrak{s},\tau);\mathbb{F}_2)\cong \widetilde{H}^{*+a_++b+n(Y,\mathfrak{s},g)}_{\mathbb{Z}_2}((N^\mu_\lambda/L^\mu_\lambda)^I;\mathbb{F}_2).
\]
The \emph{Real Fr{\o}yshov invariant} \cite[Definition 3.19]{Konno-Miyazawa-Taniguchi-2024} is defined by  \[
\delta_R(Y,\mathfrak{s},\tau) := \frac{1}{2}\,d(\operatorname{SWF}_R(Y,\mathfrak{s},\tau))=\frac{1}{2}\left( d((N^\mu_\lambda/L^\mu_\lambda)^I)-a_+-b-n(Y,\mathfrak s,g) \right).
\]
The isomorphism class of $\operatorname{SWF}_R(Y,\mathfrak{s},\tau)$, the $\mathbb{Q}$-graded $\mathbb{F}_2[W]$-module $\operatorname{SWFH}_R^*(Y,\mathfrak{s},\tau)$, and  $\delta_R(Y,\mathfrak{s},\tau)$ are invariants of $(Y,\mathfrak{s},\tau)$.

For later use, we collect some basic properties of Real Fr{\o}yshov invariants; see \cite[Section 3]{Konno-Miyazawa-Taniguchi-2024} for details. First, consider $S^3$ as the unit sphere in $\mathbb{C}^2$. Let $\tau_{\mathrm{std}}:S^3\to S^3$ be the involution $\tau_{\mathrm{std}}(z_1,z_2)=(\bar{z}_1,\bar{z}_2)$, and let $\mathfrak{s}_0$ denote the unique spin structure on $S^3$. Then \[
\operatorname{SWF}_R(S^3,\mathfrak{s}_0,\tau_{\mathrm{std}})\simeq_{\mathbb{Z}_2,\mathrm{st}} S^0, \quad \delta_R(S^3,\mathfrak{s}_0,\tau_{\mathrm{std}})=0.
\]
Next, if $-Y$ denotes the orientation reversal of $Y$, then \[
\delta_R(-Y,\mathfrak{s},\tau) = -\delta_R(Y,\mathfrak{s},\tau).
\]

Now we consider cobordisms. For $i=0,1$, let $Y_i$ be a rational homology $3$-sphere with an odd involution $\tau_i$ and a Real $\operatorname{spin}^c$ structure $\mathfrak{s}_i$ for $(Y_i,\tau_i)$. Let $W$ be a smooth compact oriented cobordism from $Y_0$ to $Y_1$, with an odd involution $\tau$ extending $\tau_0$ and $\tau_1$. Suppose $\mathfrak{t}$ is a Real $\operatorname{spin}^c$ structure for $(W,\tau)$ with $\mathfrak{t}|_{Y_i}=\mathfrak{s}_i$ for $i=0,1$.\footnote{Real $\operatorname{spin}^c$ structures on a $4$-manifold are defined similarly to the $3$-manifold case, with the additional assumption that the lift $\mathcal{I}$ preserves the $\pm$-spinor bundles.} In this case, we say that $(W,\mathfrak{t},\tau)$ is a \emph{Real spin$^c$ cobordism} from $(Y_0,\mathfrak{s}_0,\tau_0)$ to $(Y_1,\mathfrak{s}_1,\tau_1)$. 

\begin{theorem}[{\cite[Theorem 3.23]{Konno-Miyazawa-Taniguchi-2024}}]\label{thm:Real_cobordism}
    Let $(W,\mathfrak{t},\tau)$ be a Real spin$^c$ cobordism from $(Y_0,\mathfrak{s}_0,\tau_0)$ to $(Y_1,\mathfrak{s}_1,\tau_1)$, where $Y_0,Y_1$ are rational homology $3$-spheres. If $b_1(W)=0$ and $b_+^{-\tau}(W)=0$, then \[
\delta_R(Y_0,\mathfrak{s}_0,\tau_0) + \frac{c_1(\mathfrak{t})^2-\sigma(W)}{16} \leq \delta_R(Y_1,\mathfrak{s}_1,\tau_1).
\]
\end{theorem}
Here, $b_+^{-\tau}(W)$ denotes the dimension of the maximal positive definite subspace of the $(-1)$-eigenspace of $\tau^*:H^2(W;\mathbb{R})\to H^2(W;\mathbb{R})$. The following corollary is a direct consequence, where a \emph{Real spin$^c$ filling} is defined in the obvious way. 

\begin{corollary}\label{cor:real_froyshov_filling}
Let $(W,\mathfrak{t},\tau)$ be a Real spin$^c$ filling of $(Y,\mathfrak{s},\tau_Y)$. If $b_1(W)=0$ and $b_+^{-\tau}(W)=0$, then \[
\frac{c_1(\mathfrak{t})^2-\sigma(W)}{16} \leq \delta_R(Y,\mathfrak{s},\tau_Y).
\] \end{corollary}
\begin{proof} Removing a small invariant ball centered at a $\tau$-fixed point in the interior of $W$ yields a Real $\operatorname{spin}^c$ cobordism from $(S^3,\mathfrak{s}_0,\tau_{\mathrm{std}})$ to $(Y,\mathfrak{s},\tau_Y)$. Applying Theorem \ref{thm:Real_cobordism} gives the result. 
\end{proof}

\subsection{Relative Bauer--Furuta maps}\label{subsec:relative_BF}

We recall the relative Bauer--Furuta maps and some properties; see \cite{Manolescu-2003, Manolescu-2007, Khandhawit-2015, Dai-Sasahira-Stoffregen-2026} for details. 

Let $(Y,\mathfrak{s})$ be a $\operatorname{spin}^c$ rational homology sphere. We use the notation in Section \ref{subsec:SWF} regarding $(Y,\mathfrak{s})$. Let $(W,\mathfrak{t})$ be a compact, connected, oriented $\operatorname{spin}^c$ $4$-manifold with $\partial(W,\mathfrak{t})=(Y,\mathfrak{s})$ and $b_1(W)=0$. Choose a Riemannian metric that is a product near the boundary, and a $\operatorname{spin}^c$ connection $\widehat{A}_0$ on $W$ that is pulled back from $A_0$ on a boundary collar. We write $\mathbb{S}_W=\mathbb{S}_W^+\oplus \mathbb{S}^-_W$ for the spinor bundle over $W$.

For sufficiently large $k>0$ and $\mu>0$, the Seiberg--Witten map with boundary terms is given by: \begin{align}\label{eq:SW_with_boundary_terms} \begin{split}
        \operatorname{SW}^\mu_W: L^2_k(\Omega^1_{\operatorname{CC}}(W,i\mathbb{R})\oplus \Gamma(\mathbb{S}^+_W))&\longrightarrow L^2_{k-1}(\Omega^+(W,i\mathbb{R})\oplus \Gamma(\mathbb{S}^-_W))\oplus \mathcal{V}^\mu_{k-1/2}, \\
    (\alpha,\Phi) & \longmapsto \left(d^+\alpha+\frac{1}{2}F_{\widehat{A}_0}^+-q_W(\Phi), \, D^+_{\widehat{A}_0}\Phi+\rho_W(\alpha)\Phi, \, \tilde{p}^\mu(i^*\alpha,\Phi|_Y) \right),
\end{split} \end{align}
where: \begin{itemize}[leftmargin=2em]
    \item $\Omega^1_{\mathrm{CC}}(W,i\mathbb{R})$ denotes the space of imaginary-valued $1$-forms on $W$ satisfying the double Coulomb condition \cite{Khandhawit-2015}. 
    \item $F_{\widehat{A}_0}\in \Omega^2(W,i\mathbb{R})$ is the curvature of the connection on the determinant line bundle $\det \mathbb{S}^+_W$ induced by $\widehat{A}_0$. $F_{\widehat{A}_0}^+\in \Omega^+(W,i\mathbb{R})$ denotes its self-dual part.
    \item $q_W(\Phi):=\rho_W^{-1}((\Phi\otimes \Phi^*)_0)\in \Omega^+(W,i\mathbb{R})$, where $\rho_W:\Lambda^* T^*W\otimes \mathbb{C}\to \operatorname{End}(\mathbb{S}_W)$ is Clifford multiplication for $W$, which identifies $\Lambda^+T^*W\otimes i\mathbb{R}$ with $i\cdot \mathfrak{su}(\mathbb{S}_W^+)$. 
    \item $i:Y\hookrightarrow W$ denotes the inclusion.
    \item $\mathcal{V}^\mu$ denotes the spectral subspace of the global Coulomb slice $\mathcal{V}$ for $Y$ with respect to the operator $l=*d\oplus D_{A_0}$ in \eqref{eq:decomposition_l+c}, corresponding to eigenvalues in $(-\infty,\mu]$. $\tilde{p}^\mu$ denotes its spectral projection, and $\mathcal{V}^\mu_{k-1/2}$ denotes its Sobolev completion with respect to the $L^2_{k-1/2}$-norm.
\end{itemize}
The map $\operatorname{SW}^\mu_W$ is the sum of the Fredholm operator \[
\ell_W^\mu(\alpha,\Phi) :=\left(d^+\alpha, \, D_{\widehat{A}_0}^+\Phi,\, \tilde{p}^\mu (i^*\alpha,\Phi|_Y)\right)
\]
and a compact map. Choose a finite-dimensional $S^1$-invariant subspace $U'\subset L^2_{k-1}(\Omega^+(W,i\mathbb{R})\oplus \Gamma(\mathbb{S}^-_W))$ and $\lambda<0$ such that \[
\operatorname{im}\ell_W^\mu + \left(U' \oplus\mathcal{V}^\mu_\lambda\right)=L^2_{k-1}(\Omega^+(W,i\mathbb{R})\oplus \Gamma(\mathbb{S}^-_W)) \oplus\mathcal{V}^\mu_{k-1/2};
\]
recall from Section \ref{subsec:SWF} that $\mathcal{V}^\mu_\lambda$ is the finite-dimensional subspace of $\mathcal{V}$ spanned by the eigenvectors of $l$ with eigenvalues in $(\lambda,\mu]$. Then \[
U:=(\ell_W^\mu)^{-1}\left(U'\oplus\mathcal{V}^\mu_\lambda\right)
\]
is finite-dimensional. Let $\pi_{U'}:L^2_{k-1}(\Omega^+(W,i\mathbb{R})\oplus \Gamma(\mathbb{S}^-_W))\to U'$ denote the $L^2$-orthogonal projection, and define \begin{equation}\label{eq:SW_finite_dim_approximation}
    \operatorname{SW}^{\lambda,\mu}_{U'}:=(\pi_{U'}\oplus\operatorname{id}) \circ\operatorname{SW}_W^\mu|_U :U \longrightarrow U'\oplus\mathcal{V}^\mu_\lambda.
\end{equation}
By taking sufficiently large $U'$ and $\mu,-\lambda$, we can choose an $S^1$-equivariant index pair $(N^\mu_\lambda,L^\mu_\lambda)$ as in Section \ref{subsec:SWF} so that $\operatorname{SW}^{\lambda,\mu}_{U'}$ induces a based $S^1$-map \begin{equation}\label{eq:relative-BF-representative}
f_{U'}^{\lambda,\mu}: S^{U}\longrightarrow S^{U'} \wedge \left(N^\mu_\lambda/L^\mu_\lambda\right).
\end{equation}
This map represents the \emph{relative Bauer--Furuta map}, which is a morphism \[
\operatorname{BF}_{W,\mathfrak{t}}: S^{\frac{c_1(\mathfrak{t})^2-\sigma(W)}{8}\mathbb{C}} \longrightarrow \Sigma^{b_+(W)\mathbb{R}} \operatorname{SWF}_{S^1}(Y,\mathfrak{s})
\]
in the category $\mathfrak{C}_{S^1}$, where we are using the convention \eqref{eq:rational_sphere}. 

For a $\operatorname{spin}^c$ cobordism $(W,\mathfrak{t})$ between $\operatorname{spin}^c$ rational homology spheres $(Y_0,\mathfrak{s}_0)$ and $(Y_1,\mathfrak{s}_1)$, with $b_1(W)=0$, the same construction gives a morphism \[
S^{\frac{c_1(\mathfrak{t})^2-\sigma(W)}{8}\mathbb{C}} \longrightarrow \Sigma^{b_+(W)\mathbb{R}} \operatorname{SWF}_{S^1}(-Y_0,\mathfrak{s}_0) \wedge \operatorname{SWF}_{S^1}(Y_1,\mathfrak{s}_1).
\]
Composing with the duality morphism  \[
\operatorname{SWF}_{S^1}(Y_0,\mathfrak{s}_0)\wedge \operatorname{SWF}_{S^1}(-Y_0,\mathfrak{s}_0)\longrightarrow S^0
\]
in $\mathfrak{C}_{S^1}$ \cite{Manolescu-2007}, we obtain the \emph{relative Bauer--Furuta map} for the cobordism $(W,\mathfrak{t})$: \begin{equation}\label{eq:relative_BF}
    \operatorname{BF}_{W,\mathfrak{t}}: \Sigma^{\frac{c_1(\mathfrak{t})^2-\sigma(W)}{8}\mathbb{C}} \operatorname{SWF}_{S^1}(Y_0,\mathfrak{s}_0) \longrightarrow \Sigma^{b_+(W)\mathbb{R}} \operatorname{SWF}_{S^1}(Y_1,\mathfrak{s}_1).
\end{equation}

The relative Bauer--Furuta map satisfies the following adjunction relation. For any object $Z$ in the category $\mathfrak{C}_{S^1}$, let $U:Z\to  \Sigma^{\mathbb{C}}Z$ denote the morphism induced by the standard inclusion $S^0\hookrightarrow S^{\mathbb{C}}$. 
\begin{proposition}[{\cite[Proposition 3.15]{Dai-Sasahira-Stoffregen-2026}}]\label{prop:BF_adjunction}
Let $(W,\mathfrak{t})$ be a spin$^c$ filling or cobordism as above, and let $S\subset\operatorname{int}W$ be a smoothly embedded sphere with $S^2<0$. Set \[
\mathfrak{t}'=\mathfrak{t}+\operatorname{PD}[S], \quad n:=\frac{\langle c_1(\mathfrak{t}),[S]\rangle+S^2}{2} \in\mathbb{Z}. \]
Then \[
\begin{cases}
    \operatorname{BF}_{W,\mathfrak{t}'} \circ\, U^n= \operatorname{BF}_{W,\mathfrak{t}}  & \text{if $n\geq 0$}, \\
    \operatorname{BF}_{W,\mathfrak{t}'}= \operatorname{BF}_{W,\mathfrak{t}} \circ\, U^{-n} & \text{if $n<0$}.
\end{cases} \]
In particular, $\operatorname{BF}_{W,\mathfrak{t}}=\operatorname{BF}_{W,\mathfrak{t}'}$ if $n=0$.
\end{proposition}

Note that $\mathfrak{t}$ and $\mathfrak{t}'$ have the same boundary restriction. Also, since $c_1(\mathfrak{t}')=c_1(\mathfrak{t})+2\operatorname{PD}[S]$, we have \[
\frac{c_1(\mathfrak{t}')^2-\sigma(W)}{8}-\frac{c_1(\mathfrak{t})^2-\sigma(W)}{8}=n.
\]

\subsection{Compatibility}\label{subsec:compatibility}
For later use, we describe a certain type of compatibility of relative Bauer--Furuta maps; see \eqref{eq:BF_naturality} below. The argument is essentially based on \cite[Section 3.4]{Dai-Sasahira-Stoffregen-2026}.

Let $Y$ be a rational homology sphere, and let $W$ be a negative definite filling of $Y$ with $b_1(W)=0$. Suppose $\tau:W\to W$ is a smooth orientation-preserving involution such that $\tau|_Y:Y\to Y$ is an odd involution. We assume that the Riemannian metric on $W$ is $\tau$-invariant and is a product near the boundary. Let $(\mathfrak{s},\mathcal{I}_Y)$ be a Real $\operatorname{spin}^c$ structure for $(Y,\tau|_Y)$; see Section \ref{subsec:O(2)_SWF}. Let $\mathfrak{t}_+$ be a $\operatorname{spin}^c$ structure on $W$ such that $\mathfrak{t}_+|_Y=\mathfrak{s}$, and put $\mathfrak{t}_-:=\tau^* \overline{\mathfrak{t}_+}$. In particular, we have anti-linear fiberwise isometries $J^\pm:\mathbb{S}_{W,\mathfrak{t}_+}^\pm \to \mathbb{S}_{W,\mathfrak{t}_-}^\pm$ of spinor bundles covering $\tau$, compatible with Clifford multiplications. We may identify the restriction $J^+|_Y$ with $\mathcal{I}_Y:\mathbb{S}\to \mathbb{S}$. Fix an $\mathcal{I}_Y$-compatible reference $\operatorname{spin}^c$ connection $A_0$ for $\mathfrak{s}$ as in Section \ref{subsec:O(2)_SWF}. Then fix reference $\operatorname{spin}^c$ connections $\widehat{A}_\pm$ for $\mathfrak{t}_\pm$, such that $\widehat{A}_-=\tau^*\overline{\widehat{A}_+}$, and both are pulled back from $A_0$ on a boundary collar.

Now recall the construction in Section \ref{subsec:relative_BF}. Writing $\operatorname{SW}_{W,\pm}^{\mu}$ for the Seiberg--Witten maps \eqref{eq:SW_with_boundary_terms} associated to $\mathfrak{t}_\pm$, we have a commutative diagram \begin{equation}\label{eq:SW_diagram}
    \begin{tikzcd}
        L^2_k(\Omega^1_{\mathrm{CC}}(W,i\mathbb{R})\oplus \Gamma(\mathbb{S}^+_{W,\mathfrak{t}_+})) \ar[d, "\mathcal{J}"'] \ar[r, "\operatorname{SW}^\mu_{W,+}"] & L^2_{k-1}(\Omega^+(W,i\mathbb{R})\oplus \Gamma(\mathbb{S}^-_{W,\mathfrak{t}_+}))\oplus \mathcal{V}^\mu_{k-1/2} \ar[d, "\mathcal{J}'\oplus I"] \\
        L^2_k(\Omega^1_{\mathrm{CC}}(W,i\mathbb{R})\oplus \Gamma(\mathbb{S}^+_{W,\mathfrak{t}_-})) \ar[r, "\operatorname{SW}^\mu_{W,-}"] & L^2_{k-1}(\Omega^+(W,i\mathbb{R})\oplus \Gamma(\mathbb{S}^-_{W,\mathfrak{t}_-}))\oplus \mathcal{V}^\mu_{k-1/2}, 
    \end{tikzcd}
\end{equation}
where $I$ is the involution defined by \eqref{eq:involution_on_Coulomb_slice}, and $\mathcal{J}, \mathcal{J}'$ are defined by $\mathcal{J}(\alpha,\Phi)=(-\tau^*\alpha,J^+\Phi)$ and $\mathcal{J}'(\omega,\Psi)=(-\tau^*\omega,J^-\Psi)$. The maps $\mathcal{J},\mathcal{J}'$ are $\mathbb{R}$-linear isometries satisfying $\mathcal{J}(z\cdot (\alpha,\Phi))=\bar{z}\cdot \mathcal{J}(\alpha,\Phi)$ and $\mathcal{J}'(z\cdot (\omega,\Psi))=\bar{z}\cdot \mathcal{J}'(\omega,\Psi)$ for $z\in S^1$. The diagram \eqref{eq:SW_diagram} still commutes if we replace $\operatorname{SW}^\mu_{W,\pm}$ by their linear parts $\ell^\mu_{W,\pm}$.

Choose $U'_+\subset L^2_{k-1}(\Omega^+(W,i\mathbb{R})\oplus \Gamma(\mathbb{S}^-_{W,\mathfrak{t}_+}))$ and $\lambda<0<\mu$ as in Section \ref{subsec:relative_BF}, and set $U'_-:=\mathcal{J}'(U'_+)$. The corresponding inverse images $U_+$ and $U_-$ satisfy $U_-=\mathcal{J}(U_+)$. Moreover, the diagram \[
\begin{tikzcd}
    U_+ \ar[d, "\mathcal{J}|_{U_+}"'] \ar[rr, "\operatorname{SW}^{\lambda,\mu}_{U'_+}"] && U'_+\oplus \mathcal{V}^\mu_\lambda \ar[d, "\mathcal{J}'|_{U'_+}\oplus I"] \\
    U_- \ar[rr, "\operatorname{SW}^{\lambda,\mu}_{U'_-}"] && U'_-\oplus \mathcal{V}^\mu_\lambda
\end{tikzcd}
\]
commutes, where the horizontal maps are given by \eqref{eq:SW_finite_dim_approximation}. 

Now choose an $\operatorname{O}(2)$-equivariant index pair $(N^\mu_\lambda,L^\mu_\lambda)$ as in Section \ref{subsec:O(2)_SWF}. The involution $I$ induces an involution $I_Y$ on $N^\mu_\lambda/L^\mu_\lambda$ satisfying $I_Y(z\cdot x)=\bar{z}\cdot I_Y(x)$ for $x\in N^\mu_\lambda/L^\mu_\lambda$ and $z\in S^1$. Then the representatives (see \eqref{eq:relative-BF-representative}) \[
f_\pm:=f_{U'_\pm}^{\lambda,\mu}: S^{U_\pm}\longrightarrow S^{U'_\pm}\wedge N^\mu_\lambda/L^\mu_\lambda
\]
of the Bauer--Furuta maps $\operatorname{BF}_{W,\mathfrak{t}_\pm}$ may be chosen so that the following diagram commutes: \[
\begin{tikzcd}
    S^{U_+} \ar[rr, "f_+"] \ar[d, "\mathcal{J}"'] && S^{U'_+}\wedge N^\mu_\lambda/L^\mu_\lambda \ar[d, "\mathcal{J}' \wedge I_Y"] \\
    S^{U_-} \ar[rr, "f_-"] && S^{U'_-}\wedge N^\mu_\lambda/L^\mu_\lambda.
\end{tikzcd}
\]
Then the following diagram also commutes: \begin{equation}\label{eq:BF_representative_commutativity}
    \begin{tikzcd}
   S^{U'_-}\wedge S^{U_+} \ar[rr, "\operatorname{id}\wedge f_+"] \ar[d, "(\mathcal{J}')^{-1}\wedge \mathcal{J}"'] && S^{U'_-}\wedge S^{U'_+}\wedge N^\mu_\lambda/L^\mu_\lambda \ar[d, "(\mathcal{J}')^{-1}\wedge \mathcal{J}' \wedge I_Y"] \\
    S^{U'_+}\wedge S^{U_-} \ar[rr, "\operatorname{id}\wedge f_-"] &&  S^{U'_+}\wedge S^{U'_-}\wedge N^\mu_\lambda/L^\mu_\lambda.
\end{tikzcd}
\end{equation}
Since $c_1(\mathfrak{t}_-)=-\tau^*c_1(\mathfrak{t}_+)$, we have $c_1(\mathfrak{t}_+)^2=c_1(\mathfrak{t}_-)^2$. The commutative diagram \eqref{eq:BF_representative_commutativity} implies that \begin{equation}\label{eq:BF_naturality}
\operatorname{BF}_{W,\mathfrak{t}_-} = I_Y\circ \operatorname{BF}_{W,\mathfrak{t}_+} \circ \,I:S^{\frac{c_1(\mathfrak{t}_\pm)^2-\sigma(W)}{8}\mathbb{C}}\longrightarrow \operatorname{SWF}_{S^1}(Y,\mathfrak{s}),
\end{equation}
in $\mathfrak{C}_{S^1}$, where $I_Y$ denotes the involution coming from $\operatorname{SWF}_{\operatorname{O}(2)}(Y,\mathfrak{s},\tau|_Y)$, and $I$ denotes complex conjugation. (Since $I_Y$ and $I$ are conjugate-equivariant, the composition $I_Y\circ \operatorname{BF}_{W,\mathfrak{t}_+}\,\circ I$ is a well-defined morphism in $\mathfrak C_{S^1}$, although $I_Y$ and $I$ are not themselves morphisms in $\mathfrak{C}_{S^1}$. Also, we are using the fact that if $\mathbb{C}^m$ is an $\operatorname{O}(2)$-representation where $S^1$ acts by scalar multiplication, then the action of the element $I\in \operatorname{O}(2)$ is complex conjugation up to a change of basis by \cite[Lemma A.1]{Kang-Park-Taniguchi-2026}.) We note that the relation \eqref{eq:BF_naturality} is also used in \cite[Section 4.5]{Kang-Park-Taniguchi-2026}. 

\section{Lattice spectra}\label{sec:lattice_spectra}

We recall the construction of lattice spectra, following the conventions of \cite{Dai-Sasahira-Stoffregen-2026}.

\subsection{Lattice homology}\label{subsec:lattice_homology}

We begin by recalling the construction of lattice homology, which was developed in \cite{Nemethi-2005, Nemethi-2008}. We use the homological conventions of \cite{Dai-Sasahira-Stoffregen-2026}, which differ from the cohomological conventions in \cite{Nemethi-2005, Nemethi-2008}.

Let $\Gamma$ be a negative definite plumbing tree of spheres. Write $W_\Gamma$ for the associated $4$-manifold and give $Y_\Gamma:=\partial W_\Gamma$ its boundary orientation. Note that $Y_\Gamma$ is a rational homology sphere. Put \[
L:=H_2(W_\Gamma;\mathbb{Z}), \quad L^*:=H^2(W_\Gamma;\mathbb{Z}) \cong\operatorname{Hom}(L,\mathbb{Z}).
\]
For each vertex $v$ of $\Gamma$, we also write $v$ for the class in $L$ of the corresponding sphere. Write $Q_\Gamma$ for the intersection pairing and define $q:L\to L^*$ by $q(x)=Q_\Gamma(x,-)$. We write $v^*:=q(v)$, which denotes the Poincar\'e dual of $v\in L$. The $\mathbb{Q}$-linear extension $q_{\mathbb{Q}}:L\otimes \mathbb{Q}\to L^*\otimes \mathbb{Q}$ of $q$ is an isomorphism. We use this to define $k^2$ for $k\in L^*$. 

The set of \emph{characteristic covectors} is
\[
\operatorname{Char}(\Gamma):=\left\{k\in L^*:k(x)\equiv Q_\Gamma(x,x)\pmod{2} \text{ for all }x\in L\right\}.
\]
The first Chern class map $c_1:\operatorname{Spin}^c(W_\Gamma)\to \operatorname{Char}(\Gamma)$ is a bijection. Restriction to the boundary induces an identification $\operatorname{Spin}^c(Y_\Gamma)\cong \operatorname{Char}(\Gamma)/2q(L)$. Thus, for any $k\in \operatorname{Char}(\Gamma)$, its orbit $[k]:=k+2q(L)\subset \operatorname{Char}(\Gamma)$ represents a $\operatorname{spin}^c$ structure on $Y_\Gamma$.

Fix an orbit $[k]$. We regard it as an affine lattice in $L^*\otimes \mathbb{R}$. This gives a cubical decomposition of $L^*\otimes \mathbb{R}$. We call any $\ell\in [k]$ a \emph{lattice point}. For $\ell\in[k]$ and a subset $A$ of the vertices of $\Gamma$, put \[
\square(\ell,A):= \left\{ \ell+2\sum_{v\in A}t_vv^*: t_v\in [0,1]\right\}.
\]
This is a cube of dimension $|A|$, called a \emph{lattice cube}. A lattice point is a lattice cube of dimension $0$.

Let $|\Gamma|$ denote the number of vertices of $\Gamma$. The \emph{weights} of a lattice point $\ell$ and a lattice cube $\sigma$ are defined by \[
w(\ell):=\frac{\ell^2+|\Gamma|}{4}, \quad w(\sigma):= \min\{w(\ell):\ell\text{ is a vertex of }\sigma\}.
\]
All the weights lie in the single coset $\mathcal{A}_{[k]}:=w(k)+2\mathbb{Z}\subset \mathbb{Q}$. 

Choose orientations of the lattice cubes. For a codimension-one face $\tau$ of a cube $\sigma$, define $\epsilon(\sigma,\tau)$ to be $\pm 1$ according to whether the orientation of $\tau$ agrees with the boundary orientation induced by $\sigma$. Let $C_d(\Gamma,[k])$ be the free $\mathbb{Z}[U]$-module generated by the $d$-dimensional lattice cubes, and put $C_*(\Gamma,[k]):= \bigoplus_{d\geq0}C_d(\Gamma,[k])$. This complex has two gradings: the \emph{dimensional grading} and the \emph{Maslov grading}, given by \[
\dim (U^n\sigma):=d, \quad \operatorname{gr}(U^n\sigma):=w(\sigma)+d-2n \quad(n\geq 0),
\]
for a $d$-dimensional lattice cube $\sigma$. Define a $\mathbb{Z}[U]$-linear differential by \begin{equation}\label{eq:differential_lattice_homology}
    \partial\sigma := \sum_{\substack{\tau\subset\sigma\\ \dim\tau=\dim\sigma-1}} \epsilon(\sigma,\tau)\, U^{(w(\tau)-w(\sigma))/2}\tau.
\end{equation}
The differential satisfies $\partial^2=0$ and decreases both gradings by one. Different choices of cube orientations give isomorphic chain complexes. The \emph{lattice homology} is defined by \[
\mathbb{H}_d(\Gamma,[k]) := H_d(C_*(\Gamma,[k]),\partial), \quad \mathbb{H}_*(\Gamma,[k]) := \bigoplus_{d\geq0}\mathbb{H}_d(\Gamma,[k]).
\]

Here is another description of lattice homology. For $a\in\mathcal A_{[k]}$, consider the superlevel set: \begin{equation}\label{eq:superlevel_sets}
    S_a:= \bigcup_{\text{lattice cubes $\sigma$ with $w(\sigma)\geq a$}}\sigma \subset L^*\otimes \mathbb{R}.
\end{equation}
The negative definiteness of $Q_\Gamma$ implies that $S_a$ is a finite cubical complex. The summand of $\mathbb{H}_d(\Gamma,[k])$ in Maslov grading $a+d$ is naturally isomorphic to $H_d(S_a;\mathbb{Z})$. Under these identifications, multiplication by $U$ is induced by the inclusion $S_a\hookrightarrow S_{a-2}$; see \cite[Theorem 3.1.12]{Nemethi-2008}.

\subsection{Almost rational graphs and carrying paths}\label{subsec:carrying_paths}

We first recall the definition of rational and almost rational graphs \cite{Nemethi-2005}. Let $\Gamma$ be a negative definite plumbing tree, and let $K\in\operatorname{Char}(\Gamma)$ be the \emph{canonical characteristic covector}, determined by \[
K(v)=-Q_\Gamma(v,v)-2
\]
for every vertex $v$ of $\Gamma$. Put \[
\chi_K(x):=-\frac{Q_\Gamma(x,x)+K(x)}2 \quad(x\in L).
\]
An element $x=\sum_v n_vv\in L$ is called \emph{effective} if $n_v\geq0$ for every $v$. The graph $\Gamma$ is called \emph{rational} if $\chi_K(x)\geq 1$ for every nonzero effective $x\in L$. We say that $\Gamma$ is \emph{almost rational} (abbreviated AR) if decreasing the weight of at most one vertex produces a rational graph. In particular, rational graphs are AR, and every negative definite star-shaped plumbing tree is AR \cite[Example 8.2]{Nemethi-2005}. For an AR graph $\Gamma$, lattice homology is concentrated in dimensional grading zero \cite[Theorem 4.3.3]{Nemethi-2008}, i.e., $\mathbb{H}_d(\Gamma,[k])=0$ ($d>0$) for every orbit $[k]$. From now on, we assume that $\Gamma$ is AR. 

A \emph{lattice path} in $[k]$ is a finite sequence $\gamma=(k_0,\dots,k_N)$ ($N\geq 0$) of distinct lattice points such that \[ 
k_{i+1}=k_i+2v_i^* \quad (0\leq i<N)
\]
for some vertices $v_i$ of $\Gamma$. We also write $\gamma$ for the cubical subcomplex consisting of these lattice points and the edges $e_i:=[k_i,k_{i+1}]$. When $N=0$, this is the one-point path $\{k_0\}$.

For a lattice path $\gamma$, let $C_*(\gamma,[k])\subset C_*(\Gamma,[k])$ be the subcomplex generated over $\mathbb{Z}[U]$ by the vertices and edges of $\gamma$, and write \[
\mathbb{H}_d(\gamma,[k]):=H_d(C_*(\gamma,[k]),\partial).
\]
We have $\mathbb{H}_d(\gamma,[k])=0$ for $d>0$ \cite[Lemma 4.9]{Dai-Sasahira-Stoffregen-2026}. We say that the path $\gamma$ \emph{carries} the lattice homology of $(\Gamma,[k])$ if the inclusion $C_*(\gamma,[k])\hookrightarrow C_*(\Gamma,[k])$ induces an isomorphism on homology. Since $\Gamma$ is AR, this is equivalent to requiring that the induced map on $\mathbb{H}_0$ is an isomorphism. Recall from the previous subsection the superlevel set description of lattice homology; the following criterion will be used in Section \ref{sec:proof_main_existence}.

\begin{lemma}\label{lem:carrying_criterion}
Let $\Gamma$ be an AR graph, and let $[k]\subset \operatorname{Char}(\Gamma)$ be an orbit. A lattice path $\gamma$ in $[k]$ carries the lattice homology of $(\Gamma,[k])$ if and only if, for every $a\in\mathcal{A}_{[k]}$, the inclusion $\gamma\cap S_a\hookrightarrow S_a$ induces a bijection \[
\pi_0(\gamma\cap S_a)\longrightarrow\pi_0(S_a).
\]
\end{lemma}
\begin{proof} The superlevel set description gives a graded $\mathbb{Z}[U]$-module isomorphism \[
\mathbb{H}_0(\Gamma,[k])\cong \bigoplus_{a\in \mathcal{A}_{[k]}} H_0(S_a;\mathbb{Z}), \]
where: \begin{itemize}[leftmargin=2em]
    \item $\mathbb{H}_0(\Gamma,[k])$ is graded by the Maslov grading.
    \item $H_0(S_a;\mathbb{Z})$ has grading $a$, and the $U$-action is induced by the inclusion $S_a\hookrightarrow S_{a-2}$. 
\end{itemize}
This also applies to the subcomplex $\gamma$ and is natural with respect to inclusion \cite[3.1.14]{Nemethi-2008}, so we have the following commutative diagram of graded $\mathbb{Z}[U]$-modules: \[
\begin{tikzcd}
    \mathbb{H}_0(\gamma,[k]) \ar[d] \ar[r, "\cong"] & \bigoplus_{a\in \mathcal{A}_{[k]}}H_0(\gamma\cap S_a;\mathbb{Z}) \ar[d] \\
    \mathbb{H}_0(\Gamma,[k]) \ar[r, "\cong"] & \bigoplus_{a\in \mathcal{A}_{[k]}} H_0(S_a;\mathbb{Z}).
\end{tikzcd}
\]
Now, the lattice path $\gamma$ carries the lattice homology if and only if the left vertical map is an isomorphism. By commutativity, this is equivalent to the right vertical map being an isomorphism. It is clear that the right vertical map is an isomorphism if and only if $\gamma\cap S_a\hookrightarrow S_a$ induces a bijection on $\pi_0$ for each $a$. 
\end{proof}

Finally, we note that every orbit $[k]\subset \operatorname{Char}(\Gamma)$ admits a lattice path carrying the lattice homology of $(\Gamma,[k])$ \cite[Theorem 4.10]{Dai-Sasahira-Stoffregen-2026}.

\subsection{Lattice spectra}\label{subsec:lattice_spectra}

We next describe the construction of lattice spectra, following \cite{Dai-Sasahira-Stoffregen-2026}; see also \cite[Section 4.2]{Kang-Park-Taniguchi-2026}.

Let $\Gamma$ be a negative definite AR graph, and fix an orbit $[k]\subset \operatorname{Char}(\Gamma)$. Let $\gamma=(k_0,\dots,k_N)$ be a finite lattice path in $[k]$, which consists of vertices $k_i$ and edges $e_i:=[k_i,k_{i+1}]$. Choose a cutoff $h\in \mathcal{A}_{[k]}=w(k)+2\mathbb{Z}$ so that $h\leq w(k_i)$ for each $i=0,\dots,N$. Then $(w(k_i)-h)/2$ and $(w(e_i)-h)/2$ are nonnegative integers. We define a based $S^1$-space $X_{\gamma,h}$ by (see Figure \ref{fig:lattice_spectra}): \begin{equation}\label{eq:path_lattice_space}
    X_{\gamma,h} := \left(\coprod_{i=0}^N S^{\frac{w(k_i)-h}{2}\mathbb{C}} \amalg \coprod_{i=0}^{N-1} \left( S^{\frac{w(e_i)-h}{2}\mathbb{C}}\wedge [0,1]_+\right) \right)/{\sim}, 
\end{equation}
where: \begin{itemize}[leftmargin=2em]
    \item $[0,1]$ has the trivial $S^1$-action. 
    \item For each $i=0,\dots,N-1$, the two ends $S^{\frac{w(e_i)-h}{2}\mathbb{C}}\times \{0\}$ and $S^{\frac{w(e_i)-h}{2}\mathbb{C}}\times \{1\}$ of the (reduced) cylinder $S^{\frac{w(e_i)-h}{2}\mathbb{C}}\wedge [0,1]_+$ are attached to the spheres $S^{\frac{w(k_i)-h}{2}\mathbb{C}}$ and $S^{\frac{w(k_{i+1})-h}{2}\mathbb{C}}$, respectively, via the standard inclusions. 
\end{itemize}
The space $X_{\gamma,h}$ has a finite $S^1$-CW complex structure. For any integer $r\geq 0$, this construction gives an $S^1$-equivariant homeomorphism $X_{\gamma,h-2r} \cong \Sigma^{r\mathbb{C}}X_{\gamma,h}$. Consequently, the $S^1$-stable equivalence class of
\begin{equation}\label{eq:path_lattice_spectrum}
\mathcal{H}(\gamma,[k]) := \left(X_{\gamma,h},0,-h/2\right) = \Sigma^{\frac{h}{2}\mathbb C}X_{\gamma,h}
\end{equation}
in the category $\mathfrak{C}_{S^1}$ is independent of the choice of the cutoff $h$. We call the object $\mathcal{H}(\gamma,[k])$ the \emph{$S^1$-equivariant path lattice spectrum} of $\gamma$. If $\gamma$ carries the lattice homology of $(\Gamma,[k])$, then the $S^1$-stable equivalence class of $\mathcal{H}(\gamma,[k])$ depends only on $(\Gamma,[k])$; see \cite[Section 5.3]{Dai-Sasahira-Stoffregen-2026}. We therefore define the \emph{$S^1$-equivariant lattice spectrum} of $(\Gamma,[k])$ by \[
\mathcal{H}(\Gamma,[k]):= \mathcal{H}(\gamma,[k])
\]
for any carrying path $\gamma$.

\subsection{Lattice spectra and Seiberg--Witten Floer spectra}\label{subsec:spectra_comparison}
We now describe the equivalence between lattice spectra and Seiberg--Witten Floer spectra established in \cite[Theorem 1.1]{Dai-Sasahira-Stoffregen-2026}. 

\begin{figure}[t]
    \centering
    \begin{tikzpicture}[scale=0.8]
 \draw (-7,5) node[above]{$S^{\frac{w(k_0)-h}{2}\mathbb{C}}$};
  \draw (-3,5) node[above]{$S^{\frac{w(k_1)-h}{2}\mathbb{C}}$};
   \draw (3,5) node[above]{$S^{\frac{w(k_{N-1})-h}{2}\mathbb{C}}$};
    \draw (7,5) node[above]{$S^{\frac{w(k_N)-h}{2}\mathbb{C}}$};
    \draw[thick] (-3,4) circle (1);
    \draw[thick] (-7,4) circle (1);
    \draw[thick] (3,4) circle (1);
    \draw[thick] (7,4) circle (1);
    \draw (0,5) node{\textcolor{Emerald}{$\cdots$}};
    \draw (0,3) node{$\cdots$};
    \draw (0,0) node{$\operatorname{SWF}_{S^1}(Y_\Gamma,\mathfrak{s})$};
    \draw[thick]  (-7,3)--(-1,3)  (7,3)--(1,3)  ;
    \draw[thick, Emerald] (-7,5)--(-1,5) (7,5)--(1,5);
    \draw (-7,5) node[circle, fill, inner sep=1.2pt, black]{};
    \draw (-3,5) node[circle, fill, inner sep=1.2pt, black]{};
    \draw (3,5) node[circle, fill, inner sep=1.2pt, black]{};
    \draw (7,5) node[circle, fill, inner sep=1.2pt, black]{};
    \draw (-7,3) node[circle, fill, inner sep=1.2pt, black]{};
    \draw (-3,3) node[circle, fill, inner sep=1.2pt, black]{};
    \draw (3,3) node[circle, fill, inner sep=1.2pt, black]{};
    \draw (7,3) node[circle, fill, inner sep=1.2pt, black]{};
   \draw[->] (-7,2.8)  -- node[pos=0.45,anchor=north east,inner sep=2pt, font=\small] {$\operatorname{BF}_{W_\Gamma,\mathfrak{t}_0}$} (-1.5,0);
    \draw[->] (7,2.8) -- node[pos=0.45,anchor=north west,inner sep=2pt, font=\small] {$\operatorname{BF}_{W_\Gamma,\mathfrak{t}_N}$} (1.5,0);
    \draw[->] (-3,2.8) -- node[pos=0.30,anchor=north east,inner sep=2pt, font=\small]{$\operatorname{BF}_{W_\Gamma,\mathfrak{t}_1}$}  (-0.5,0.5);
    \draw[->] (3,2.8) -- node[pos=0.30,anchor=north west,inner sep=2pt, font=\small]   {$\operatorname{BF}_{W_\Gamma,\mathfrak{t}_{N-1}}$}  (0.5,0.5);
    \draw[thin,densely dashed] (-7,5) arc[start angle=90,end angle=270, x radius=0.35,y radius=1];
    \draw[thin] (-7,5) arc[start angle=90,end angle=-90, x radius=0.35,y radius=1];
    \draw[thin,densely dashed] (-3,5) arc[start angle=90,end angle=270, x radius=0.35,y radius=1];
    \draw[thin] (-3,5) arc[start angle=90,end angle=-90, x radius=0.35,y radius=1];
    \draw[thin,densely dashed] (3,5) arc[start angle=90,end angle=270, x radius=0.35,y radius=1];
    \draw[thin] (3,5) arc[start angle=90,end angle=-90, x radius=0.35,y radius=1];
    \draw[thin,densely dashed] (7,5) arc[start angle=90,end angle=270, x radius=0.35,y radius=1];
    \draw[thin] (7,5) arc[start angle=90,end angle=-90, x radius=0.35,y radius=1];
\end{tikzpicture}
    \caption{Schematic description of $X_{\gamma,h}$ and the morphism $\mathcal{T}_\gamma$. The colored edge is collapsed to a single basepoint.}
    \label{fig:lattice_spectra}
\end{figure}
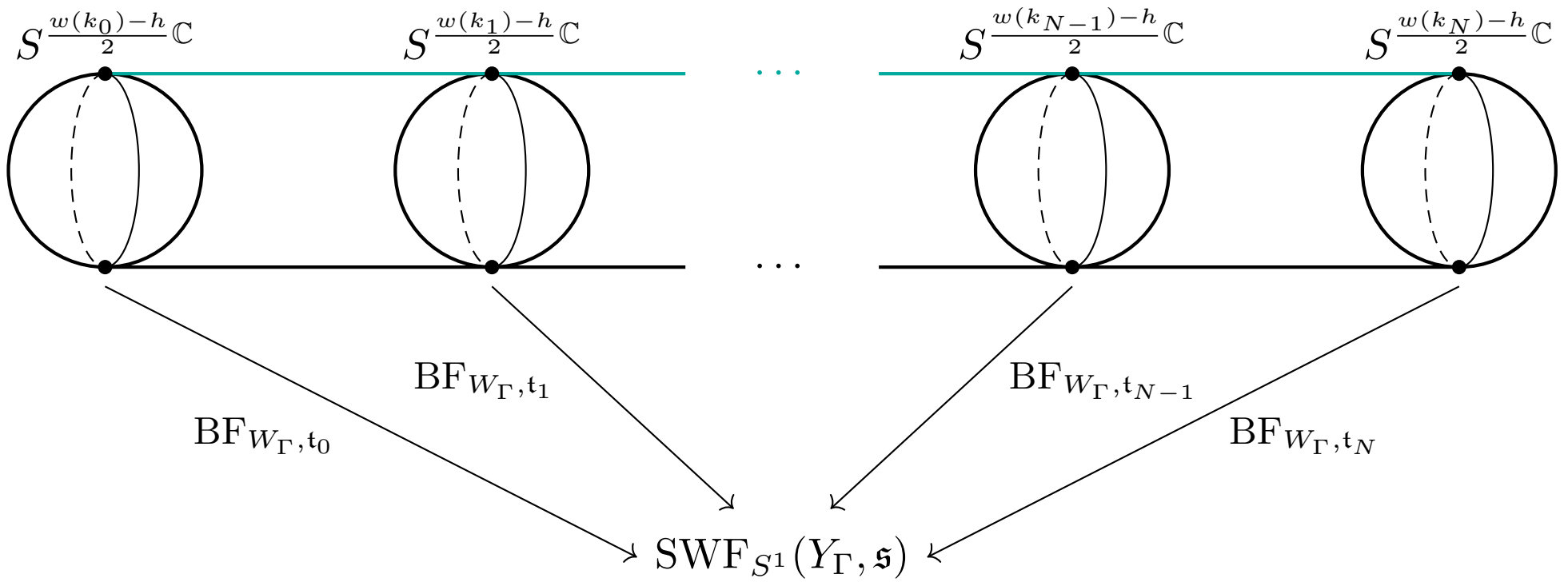

Again, let $\Gamma$ be a negative definite AR graph, and fix an orbit $[k]\subset \operatorname{Char}(\Gamma)$. Let $\gamma=(k_0,\dots,k_N)$ be a lattice path carrying the lattice homology of $(\Gamma,[k])$. Let $\mathfrak{s}$ be the $\operatorname{spin}^c$ structure on $Y_\Gamma$ corresponding to the orbit $[k]$, and let $\mathfrak{t}_i$ be the $\operatorname{spin}^c$ structure on $W_\Gamma$ with $c_1(\mathfrak t_i)=k_i$. Then we have \[
\frac{c_1(\mathfrak{t}_i)^2-\sigma(W_\Gamma)}{8}=\frac{k_i^2+|\Gamma|}{8}=\frac{w(k_i)}{2}.
\]
Thus, the relative Bauer--Furuta map has the form \[
\operatorname{BF}_{W_\Gamma,\mathfrak{t}_i}: S^{\frac{w(k_i)}{2}\mathbb{C}} \longrightarrow \operatorname{SWF}_{S^1}(Y_\Gamma,\mathfrak{s}).
\]
On the other hand, in the construction \eqref{eq:path_lattice_space}, the inclusion $S^{\frac{w(k_i)-h}{2}\mathbb{C}}\hookrightarrow X_{\gamma,h}$ induces a morphism \[
j_i:S^{\frac{w(k_i)}{2}\mathbb{C}} \longrightarrow\mathcal H(\gamma,[k]) \quad (i=0,\dots,N)
\]
in $\mathfrak{C}_{S^1}$. For successive vertices $k_{i+1}=k_i+2v_i^*$ of the path $\gamma$, the corresponding $\operatorname{spin}^c$ structures satisfy $\mathfrak{t}_{i+1}=\mathfrak{t}_i+v_i^*=\mathfrak{t}_i+\operatorname{PD}(v_i)$, and we have \[
\frac{w(k_{i+1})-w(k_i)}{2} = \frac{\langle c_1(\mathfrak{t}_i),v_i \rangle +Q_\Gamma(v_i,v_i)}{2}.
\]
Applying Proposition \ref{prop:BF_adjunction} gives \begin{equation}\label{eq:adjunction}
  \operatorname{BF}_{W_\Gamma,\mathfrak{t}_i} \circ \, U^{\frac{w(k_i)-w(e_i)}{2}} = \operatorname{BF}_{W_\Gamma,\mathfrak{t}_{i+1}} \circ \, U^{\frac{w(k_{i+1})-w(e_i)}{2}}: S^{\frac{w(e_i)}{2}\mathbb{C}}\longrightarrow \operatorname{SWF}_{S^1}(Y_\Gamma,\mathfrak{s}) \quad (i=0,\dots,N-1).  
\end{equation}
After a common stabilization, choose representatives of the Bauer--Furuta maps and based $S^1$-homotopies representing \eqref{eq:adjunction}. Now, gluing the representatives of $\operatorname{BF}_{W_\Gamma,\mathfrak{t}_i}$ on the vertex spheres $S^{\frac{w(k_i)-h}{2}\mathbb{C}}$ and these homotopies on the edge cylinders $S^{\frac{w(e_i)-h}{2}\mathbb{C}}\wedge [0,1]_+$ defines a morphism \begin{equation}\label{eq:lattice_comparison_map}
\mathcal{T}_\gamma: \mathcal{H}(\gamma,[k]) \longrightarrow \operatorname{SWF}_{S^1}(Y_\Gamma,\mathfrak{s})
\end{equation}
in $\mathfrak{C}_{S^1}$ satisfying \[
\mathcal{T}_\gamma\circ j_i = \operatorname{BF}_{W_\Gamma,\mathfrak{t}_i} : S^{\frac{w(k_i)}{2}\mathbb{C}} \longrightarrow \operatorname{SWF}_{S^1}(Y_\Gamma,\mathfrak{s}) \quad (i=0,\dots, N).
\]
This is the construction of \cite[Lemma 5.10]{Dai-Sasahira-Stoffregen-2026}; see Figure \ref{fig:lattice_spectra}. 

\begin{theorem}[{\cite{Dai-Sasahira-Stoffregen-2026}}]\label{thm:lattice_SWF_equivalence}
Suppose that $\Gamma$ is almost rational and that $\gamma$ carries the lattice homology of $(\Gamma,[k])$.
Then the morphism $\mathcal{T}_\gamma$ is an $S^1$-stable equivalence. In particular, we have \[
\mathcal{H}(\Gamma,[k]) \simeq_{S^1,\mathrm{st}} \operatorname{SWF}_{S^1}(Y_\Gamma,\mathfrak{s}).
\]
\end{theorem}

\section{Proof of Theorem \texorpdfstring{\ref{thm:main}}{1.1}}\label{sec:proof_main}
\subsection{Almost \texorpdfstring{$I$}{I}-invariant paths}\label{subsec:almost_I_invariant_paths}

Let $\Gamma$ be a negative definite AR graph. Write $W:=W_\Gamma$ and $Y:=\partial W_\Gamma$. Let $\tau:W\to W$ be a smooth orientation-preserving involution whose restriction $\tau_Y:=\tau|_Y$ is odd. Fix a Real $\operatorname{spin}^c$ structure $\mathfrak{s}$ for $(Y,\tau_Y)$, and let $[k]\subset \operatorname{Char}(\Gamma)$ be the corresponding orbit of characteristic covectors.

Consider the involution $\mathfrak{t}\mapsto \tau^*\overline{\mathfrak{t}}$ on $\operatorname{Spin}^c(W)$. If $\mathfrak{t}|_Y=\mathfrak{s}$, then $(\tau^*\overline{\mathfrak{t}})|_Y\cong \tau^*\overline{\mathfrak{s}}\cong \mathfrak{s}$. On $\operatorname{Char}(\Gamma)$, the corresponding involution is given by $I(\ell):=-\tau^*\ell$. It preserves the orbit $[k]$ and the weight function $w$. 

Following \cite[Section 4.3]{Kang-Park-Taniguchi-2026}, we define an \emph{almost $I$-invariant path} to be a finite sequence \[
\gamma=(\mathfrak{t}_{-n},\dots,\mathfrak{t}_{-1}, \mathfrak{t}_1,\dots,\mathfrak{t}_n) \quad (n\geq 1)
\]
of $\operatorname{spin}^c$ structures on $W$ satisfying the following conditions, where $k_i:=c_1(\mathfrak{t}_i)$ ($i\in \{\pm1,\dots,\pm n\}$).  
\begin{enumerate}[label=(\arabic*), leftmargin=2em]
\item $\mathfrak{t}_i|_Y=\mathfrak{s}$ for every $i\in \{\pm 1,\dots, \pm n\}$. 
\item\label{almost_path_(2)} $\mathfrak{t}_{-i}=\tau^*\overline{\mathfrak{t}}_i$ for every $i\in \{1,\dots,n\}$. 
\item The sequences $\gamma_-:=(k_{-n},\dots,k_{-1})$ and $\gamma_+:=(k_1,\dots,k_n)$ are lattice paths in the sense of Section \ref{subsec:carrying_paths}. This is equivalent to requiring that, for every $i\in \{-n,\dots,-2\}\cup \{1,\dots,n-1\}$, there exists a vertex $v$ in $\Gamma$ such that $\mathfrak{t}_{i+1}=\mathfrak{t}_i+v^*$. 
\item\label{almost_path_(4)} There exists a (possibly empty) set $B$ consisting of pairwise nonadjacent vertices of $\Gamma$ such that, for every $v\in B$, the corresponding sphere $S_v\subset W$ is preserved setwise and orientation-preservingly by $\tau$, and $\mathfrak{t}_1-\mathfrak{t}_{-1}=\sum_{v\in B} v^*$. 
\end{enumerate}

We also use $\gamma$ to denote the corresponding sequence $(k_{-n},\dots,k_{-1},k_1,\dots,k_n)$ of characteristic covectors. Note that this sequence is not necessarily a lattice path in the sense of Section \ref{subsec:carrying_paths}. Condition \ref{almost_path_(2)} implies $w(k_{-i})=w(k_i)$ for $i=1,\dots,n$.

\begin{lemma}\label{lem:central_adjunction}
For an almost $I$-invariant path $\gamma$, the relative Bauer--Furuta maps satisfy \[
\operatorname{BF}_{W,\mathfrak{t}_{-1}} = \operatorname{BF}_{W,\mathfrak{t}_1} : S^{\frac{w(k_{\pm 1})}{2}\mathbb{C}}\longrightarrow \operatorname{SWF}_{S^1}(Y,\mathfrak{s})
\]
in the category $\mathfrak C_{S^1}$.
\end{lemma}
\begin{proof} If the set $B$ in \ref{almost_path_(4)} is empty, then $\mathfrak{t}_{-1}=\mathfrak{t}_1$, so the result is obvious. Assume that $B$ is nonempty. Then $k_1-k_{-1}=2\sum_{j=1}^r\operatorname{PD}[S_j]$, where $S_1,\dots,S_r$ are pairwise disjoint vertex spheres in $W$. It follows that $\langle k_1- k_{-1},[S_j]\rangle = 2S_j^2$. On the other hand, since $k_{-1}=-\tau^*k_1$ by \ref{almost_path_(2)} and $\tau_*[S_j]=[S_j]$ by \ref{almost_path_(4)}, $\langle k_{-1},[S_j]\rangle=-\langle k_1,[S_j]\rangle$. Consequently, we have \[
\langle k_{-1},[S_j]\rangle=-S_j^2, \quad \langle k_1,[S_j]\rangle=S_j^2.
\]
Now, for $j=0,\dots,r$, define
\[
\mathfrak{t}^{(j)}:=\mathfrak{t}_{-1} + \sum_{\ell=1}^j\operatorname{PD}[S_\ell].
\]
Note that $\mathfrak{t}^{(0)}=\mathfrak{t}_{-1}$ and $\mathfrak{t}^{(r)}=\mathfrak{t}_1$. Since the spheres $S_1,\dots,S_r$ are pairwise disjoint, we have \[
\langle c_1(\mathfrak{t}^{(j-1)}),[S_j]\rangle = \langle k_{-1},[S_j]\rangle = -S_j^2 \quad (j=1,\dots,r).
\]
Proposition \ref{prop:BF_adjunction} therefore gives \[
\operatorname{BF}_{W,\mathfrak{t}^{(j-1)}} = \operatorname{BF}_{W,\mathfrak{t}^{(j)}} \quad (j=1,\dots,r),
\]
completing the proof. 
\end{proof}

\begin{remark}\label{rmk:central_cube}
    The set $B$ in condition \ref{almost_path_(4)} is uniquely determined by $\gamma$, since the covectors $v^*$ are linearly independent. Put \[
w_c:=w(k_{-1})=w(k_1), \quad \square_\gamma:=\square(k_{-1},B).
\]
For each subset $A\subset B$, write $k_A:=k_{-1}+2\sum_{v\in A}v^*$. In particular, $k_\varnothing=k_{-1}$ and $k_B=k_1$. The proof of Lemma \ref{lem:central_adjunction} shows that $\langle k_{-1},v\rangle=-Q_\Gamma(v,v)$ for every $v\in B$. Since the vertices in $B$ are pairwise nonadjacent, it follows that
\[
k_A^2 = k_{-1}^2 + 4\sum_{v\in A} \left(\langle k_{-1},v\rangle+Q_\Gamma(v,v)\right) = k_{-1}^2.
\]
Consequently, $w(k_A)=w_c$ for every subset $A\subset B$. Thus every vertex, and hence every face, of the lattice cube $\square_\gamma$ has weight $w_c$. Moreover, condition \ref{almost_path_(4)} implies $\tau^*v^*=v^*$ for every $v\in B$. Therefore, $I(k_A)=k_1-2\sum_{v\in A}v^*=k_{B\setminus A}$. In particular, $I$ preserves $\square_\gamma$ and exchanges its opposite corners $k_{-1}$ and $k_1$. When $B=\varnothing$, the cube $\square_\gamma$ consists of the single $I$-fixed lattice point $k_{-1}=k_1$.
\end{remark}

Although an almost $I$-invariant path $\gamma$ is not necessarily a lattice path in the sense of Section \ref{subsec:carrying_paths}, we can complete it to a lattice path as follows. Choose an ordering $B=\{v_1,\dots,v_r\}$, where $r=|B|$, and put \[
\kappa_j:=k_{-1}+2\sum_{\ell=1}^j v_\ell^* \quad (j=0,\dots, r).
\]
Thus $\kappa_0=k_{-1}$ and $\kappa_r=k_1$. Replace $(k_{-1},k_1)$ by $(\kappa_0,\dots,\kappa_r)$: \begin{equation}\label{eq:completion}
    \widehat{\gamma} := (k_{-n},\dots,k_{-2},  \kappa_0,\dots,\kappa_r,  k_2,\dots,k_n).
\end{equation}
When $B=\varnothing$, this replaces the repeated central vertex $k_{-1}=k_1$ by a single vertex. Every successive difference in $\widehat{\gamma}$ is of the form $2v^*$ for a vertex $v$ of $\Gamma$. It follows that the vertices of $\widehat{\gamma}$ are distinct; hence $\widehat{\gamma}$ is a lattice path in the sense of Section \ref{subsec:carrying_paths}. We call it a \emph{completion} of $\gamma$.

We say that an almost $I$-invariant path $\gamma$ \emph{carries} the lattice homology of $(\Gamma,[k])$ if a completion $\widehat{\gamma}$ carries the lattice homology of $(\Gamma,[k])$. Using Remark \ref{rmk:central_cube}, it is easy to see that this condition is independent of the ordering of $B$.

If $\gamma$ is a carrying almost $I$-invariant path, then Theorem \ref{thm:lattice_SWF_equivalence} applies to every completion $\widehat{\gamma}$ and gives an $S^1$-stable equivalence \[
\mathcal{T}_{\widehat{\gamma}}: \mathcal{H}(\widehat{\gamma},[k]) \longrightarrow \operatorname{SWF}_{S^1}(Y,\mathfrak{s}).
\]

\begin{remark}\label{rmk:Wu_class}
    For a negative definite plumbing tree $\Gamma$, recall that a \emph{Wu class} is an element of the form $\omega=\sum_v \epsilon_v v$ with $\epsilon_v\in \{0,1\}$, such that \[
Q_\Gamma(\omega,v)\equiv Q_\Gamma(v,v)\pmod{2}
\]
for every vertex $v$. Wu classes always exist, and there is a unique Wu class if and only if $\det Q_\Gamma$ is odd. Moreover, the vertices $v$ with $\epsilon_v=1$ are pairwise nonadjacent; see \cite{Neumann-1980} for details. 

Suppose that $\tau$ permutes the vertices of $\Gamma$, and let $\Gamma_{\mathrm{fix}}$ denote the full subgraph spanned by the $\tau$-fixed vertices. For an almost $I$-invariant path $\gamma$, put \[
b:=\sum_{v\in B}v.
\]
Conditions \textup{\ref{almost_path_(2)}} and \textup{\ref{almost_path_(4)}} imply $k_1+\tau^*k_1=2q(b)$. Thus, for every $\tau$-fixed vertex $v$, \[
Q_{\Gamma_{\mathrm{fix}}}(b,v) = k_1(v) \equiv Q_{\Gamma_{\mathrm{fix}}}(v,v) \pmod{2}.
\]
Consequently, $b$ is a Wu class of $\Gamma_{\mathrm{fix}}$.
\end{remark}

\subsection{The O(2)-equivariant path lattice spectrum}\label{subsec:O2_path_lattice_spectrum}

Let $\gamma=(k_{-n},\dots,k_{-1},k_1,\dots,k_n)$ be an almost $I$-invariant path as in Section \ref{subsec:almost_I_invariant_paths}. We construct an $\operatorname{O}(2)$-equivariant path lattice spectrum of $\gamma$, following \cite[Sections 4.4 and 4.5]{Kang-Park-Taniguchi-2026}.

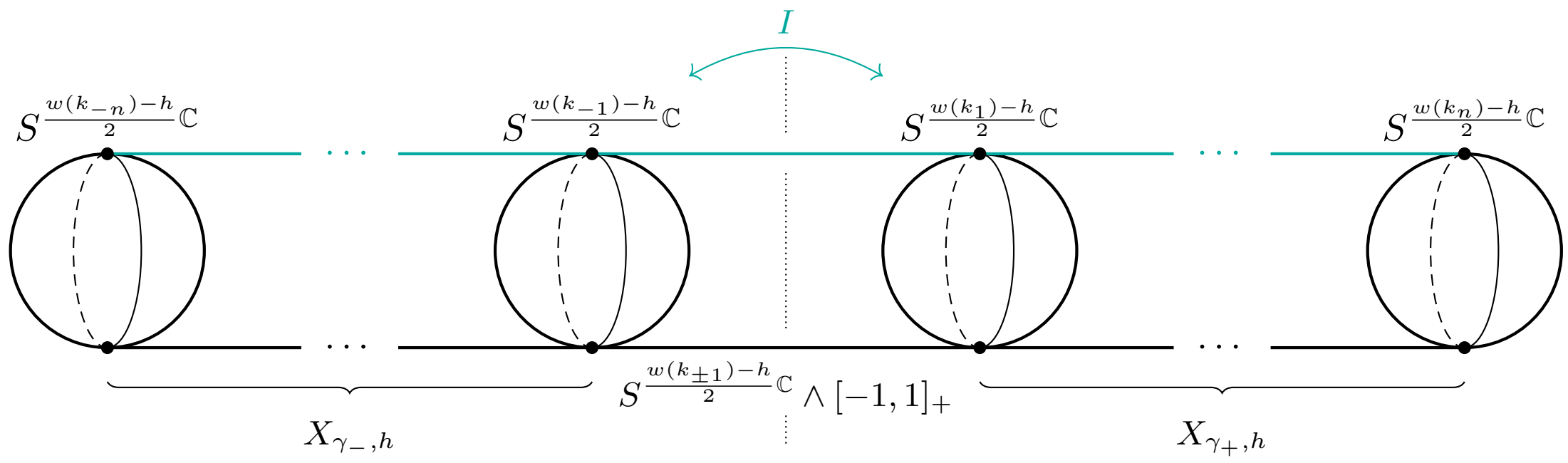
\begin{figure}[t]
    \centering
    \begin{tikzpicture}[scale=0.9]
 \draw (-7,5) node[above]{$S^{\frac{w(k_{-n})-h}{2}\mathbb{C}}$};
  \draw (-2,5) node[above]{$S^{\frac{w(k_{-1})-h}{2}\mathbb{C}}$};
   \draw (2,5) node[above]{$S^{\frac{w(k_{1})-h}{2}\mathbb{C}}$};
    \draw (7,5) node[above]{$S^{\frac{w(k_n)-h}{2}\mathbb{C}}$};
    \draw[thick] (-7,4) circle (1);
    \draw[thick] (-2,4) circle (1);
    \draw[thick] (2,4) circle (1);
    \draw[thick] (7,4) circle (1);
    \draw (-4.5,5) node{\textcolor{Emerald}{$\cdots$}};
    \draw (4.5,5) node{\textcolor{Emerald}{$\cdots$}};
    \draw (-4.5,3) node{$\cdots$};
    \draw (4.5,3) node{$\cdots$};
    \draw[thick]  (-7,3)--(-5,3) (-4,3)--(4,3)  (7,3)--(5,3)  ;
    \draw[thick, Emerald] (-7,5)--(-5,5) (-4,5)--(4,5)  (7,5)--(5,5);
    \draw (-7,5) node[circle, fill, inner sep=1.2pt, black]{};
    \draw (-2,5) node[circle, fill, inner sep=1.2pt, black]{};
    \draw (2,5) node[circle, fill, inner sep=1.2pt, black]{};
    \draw (7,5) node[circle, fill, inner sep=1.2pt, black]{};
    \draw (-7,3) node[circle, fill, inner sep=1.2pt, black]{};
    \draw (-2,3) node[circle, fill, inner sep=1.2pt, black]{};
    \draw (2,3) node[circle, fill, inner sep=1.2pt, black]{};
    \draw (7,3) node[circle, fill, inner sep=1.2pt, black]{};
    \draw[thin,densely dashed] (-7,5) arc[start angle=90,end angle=270, x radius=0.35,y radius=1];
    \draw[thin] (-7,5) arc[start angle=90,end angle=-90, x radius=0.35,y radius=1];
    \draw[thin,densely dashed] (-2,5) arc[start angle=90,end angle=270, x radius=0.35,y radius=1];
    \draw[thin] (-2,5) arc[start angle=90,end angle=-90, x radius=0.35,y radius=1];
    \draw[thin,densely dashed] (2,5) arc[start angle=90,end angle=270, x radius=0.35,y radius=1];
    \draw[thin] (2,5) arc[start angle=90,end angle=-90, x radius=0.35,y radius=1];
    \draw[thin,densely dashed] (7,5) arc[start angle=90,end angle=270, x radius=0.35,y radius=1];
    \draw[thin] (7,5) arc[start angle=90,end angle=-90, x radius=0.35,y radius=1];
    \draw (0,3) node[below]{\small{$S^{\frac{w(k_{\pm1})-h}{2}\mathbb{C}}\wedge[-1,1]_+$}};
    \draw [decorate,decoration={mirror,brace,amplitude=3pt},xshift=0pt,yshift=-10pt]
	(-7,3) -- (-2,3) node [black,midway,xshift=0pt,yshift=-15pt]{\small{$X_{\gamma_-,h}$}};
    \draw [decorate,decoration={mirror,brace,amplitude=3pt},xshift=0pt,yshift=-10pt]
	(2,3) -- (7,3) node [black,midway,xshift=0pt,yshift=-15pt]{\small{$X_{\gamma_+,h}$}};
    \draw[<->, Emerald] (-1,5.8) to[bend left] (1,5.8);
    \draw[densely dotted] (0,6)--(0,5.2) (0,4.8)--(0,3.2) (0,2.3)--(0,2);
   \draw (0,6.1) node[above] {\small{\textcolor{Emerald}{$I$}}};
\end{tikzpicture}
    \caption{Schematic description of $X_{\gamma,h}$ and the involution $I$. The colored edge is collapsed to a single basepoint.}
    \label{fig:schematic_description_2}
\end{figure}

Choose a cutoff $h\in\mathcal{A}_{[k]}$ such that $h\leq w(k_i)$ for every $i\in \{\pm1,\dots,\pm n\}$. Let $X_{\gamma_-,h}$ and $X_{\gamma_+,h}$ denote the based $S^1$-spaces constructed as in \eqref{eq:path_lattice_space}, associated to the paths $\gamma_-=(k_{-n},\dots,k_{-1})$ and $\gamma_+=(k_1,\dots,k_n)$, respectively. Define a based $S^1$-space by (see Figure \ref{fig:schematic_description_2}):
\begin{equation}\label{eq:almost_I_path_space}
X_{\gamma,h} := \left( X_{\gamma_-,h} \amalg \left(S^{\frac{w(k_{\pm1})-h}{2}\mathbb{C}}\wedge[-1,1]_+\right)\amalg X_{\gamma_+,h} \right)/{\sim},
\end{equation}
where: \begin{itemize}[leftmargin=2em]
    \item $[-1,1]$ has the trivial $S^1$-action.
    \item The two ends $S^{\frac{w(k_{\pm1})-h}{2}\mathbb{C}} \times \{-1\}$ and $S^{\frac{w(k_{\pm1})-h}{2}\mathbb{C}}\times \{1\}$ of the cylinder $S^{\frac{w(k_{\pm1})-h}{2}\mathbb{C}}\wedge[-1,1]_+$ are attached to the spheres $S^{\frac{w(k_{-1})-h}{2}\mathbb{C}}\subset X_{\gamma_-,h}$ and $S^{\frac{w(k_1)-h}{2}\mathbb{C}}\subset X_{\gamma_+,h}$, respectively, via the identity maps.
\end{itemize}
We define an involution $I$ on $X_{\gamma,h}$ as follows: \begin{itemize}[leftmargin=2em]
    \item  For each $i=1,\dots,n$, $I$ exchanges the spheres $S^{\frac{w(k_{-i})-h}{2}\mathbb{C}}\subset X_{\gamma_-,h}$ and $S^{\frac{w(k_i)-h}{2}\mathbb{C}}\subset X_{\gamma_+,h}$ by complex conjugation $z\mapsto \bar{z}$. 
    \item For each $i=1,\dots,n-1$, $I$ exchanges the cylinders $S^{\frac{w(e_{-i})-h}{2}\mathbb{C}}\wedge [0,1]_+\subset X_{\gamma_-,h}$ and $S^{\frac{w(e_{i})-h}{2}\mathbb{C}}\wedge [0,1]_+\subset X_{\gamma_+,h}$ via $[z,t]\leftrightarrow [\bar{z},1-t]$, where $e_{-i}$ and $e_i$ are the edges $[k_{-(i+1)},k_{-i}]\subset \gamma_-$ and $[k_i,k_{i+1}]$, respectively.  
    \item On the central cylinder $S^{\frac{w(k_{\pm1})-h}{2}\mathbb{C}}\wedge[-1,1]_+$, $I$ is given by $[z,t]\mapsto [\bar{z},-t]$. 
\end{itemize}
The involution $I$ satisfies $I(\lambda\cdot x)=\bar{\lambda}\cdot I(x)$ for $\lambda \in S^1$ and $x\in X_{\gamma,h}$. It follows that $X_{\gamma,h}$ is a based $\operatorname{O}(2)=S^1\rtimes\langle I\rangle$-space, and $A_{\gamma,h}:=X_{\gamma_-,h}\cup X_{\gamma_+,h}$ is an $\operatorname{O}(2)$-invariant based subspace. For every integer $m\geq 0$, there are $\operatorname{O}(2)$-equivariant homeomorphisms $X_{\gamma,h-2m}\cong \Sigma^{m\mathbb{C}}X_{\gamma,h}$ and $A_{\gamma,h-2m}\cong\Sigma^{m\mathbb{C}}A_{\gamma,h}$. We therefore define \begin{equation}\label{eq:O2_path_lattice_spectrum}
\mathcal{H}_{\operatorname{O}(2)}(\gamma,[k]) :=(X_{\gamma,h},0,0,-h/2), \quad \mathcal{A}_\gamma:=(A_{\gamma,h},0,0,-h/2)
\end{equation}
as objects in the category $\mathfrak{C}_{\operatorname{O}(2)}$. Their $\operatorname{O}(2)$-stable equivalence classes are independent of the cutoff $h$. We call $\mathcal{H}_{\operatorname{O}(2)}(\gamma,[k])$ the \emph{$\operatorname{O}(2)$-equivariant path lattice spectrum} of $\gamma$. 

\begin{lemma}\label{lem:central_quotient}
    The $I$-fixed-point set of $A_{\gamma,h}$ is the basepoint, i.e., $A_{\gamma,h}^I=\{*\}$. Moreover, there is an $\operatorname{O}(2)$-equivariant based homeomorphism \[
    X_{\gamma,h}/A_{\gamma,h} \cong_{\operatorname{O}(2)} S^{\frac{w(k_{\pm 1})-h}{2}\mathbb{C}\oplus\widetilde{\mathbb{R}}}.
    \]
\end{lemma}
\begin{proof}
    The first statement is obvious since $I$ exchanges the two summands of $A_{\gamma,h}=X_{\gamma_-,h}\cup X_{\gamma_+,h}$. For the second statement, note that collapsing $A_{\gamma,h}$ collapses both ends of the central cylinder $S^{\frac{w(k_{\pm1})-h}{2}\mathbb{C}}\wedge[-1,1]_+$, which gives: \[
X_{\gamma,h}/A_{\gamma,h} \cong_{\operatorname{O}(2)} S^{\frac{w(k_{\pm1})-h}{2}\mathbb{C}}\wedge \left([-1,1]/\{-1,1\}\right)\cong_{\operatorname{O}(2)} S^{\frac{w(k_{\pm1})-h}{2}\mathbb{C}}\wedge  S^{\widetilde{\mathbb{R}}}.
\]
\end{proof}

The following lemma shows that the $\operatorname{O}(2)$-equivariant path lattice spectrum of $\gamma$ is $S^1$-stably equivalent to the $S^1$-equivariant path lattice spectrum of any completion of $\gamma$.
\begin{lemma}\label{lem:completion}
    For an almost $I$-invariant path $\gamma$, there is an $S^1$-stable equivalence \[
    \operatorname{Res}^{\operatorname{O}(2)}_{S^1} \mathcal{H}_{\operatorname{O}(2)}(\gamma,[k]) \simeq_{S^1,\mathrm{st}} \mathcal{H}(\widehat{\gamma},[k])
    \]
    for any completion $\widehat{\gamma}$ of $\gamma$.
\end{lemma}
\begin{proof}
    Let $r=|B|$ and write \[
    \widehat{\gamma}=(\underbrace{k_{-n},\dots,k_{-2},\kappa_0}_{\gamma_-}, \dots, \overbrace{\kappa_r,k_2,\dots,k_n}^{\gamma_+})
    \]
    as in \eqref{eq:completion}. Consider first the case $r>0$. The construction \eqref{eq:path_lattice_space} shows that \[
    X_{\widehat{\gamma},h}= \left(X_{\gamma_-,h} \amalg B_{\widehat{\gamma},h}  \amalg X_{\gamma_+,h}   \right)/{\sim}= A_{\gamma,h}\cup B_{\widehat{\gamma},h},
    \]
    where \[
    B_{\widehat{\gamma},h}:= \left( S^{\frac{w(\kappa_0)-h}{2}\mathbb{C}} \amalg \left(S^{\frac{w([\kappa_0,\kappa_1])-h}{2}\mathbb{C}}\wedge [0,1]_+\right) \amalg \cdots \amalg  \left(S^{\frac{w([\kappa_{r-1},\kappa_r])-h}{2}\mathbb{C}}\wedge [0,1]_+\right) \amalg S^{\frac{w(\kappa_r)-h}{2}\mathbb{C}}\right)/{\sim}.
    \]
   On the other hand, Remark \ref{rmk:central_cube} shows that the $\kappa_0,\dots,\kappa_r$ have the same weight $w_c=w(k_{\pm 1})$. Therefore, we have obvious $S^1$-equivariant based homeomorphisms \begin{equation}\label{eq:cylinder_homeomorphism}
          B_{\widehat{\gamma},h}\cong_{S^1} S^{\frac{w_c-h}{2}\mathbb{C}} \wedge [0,r]_+\cong_{S^1}  S^{\frac{w(k_{\pm 1})-h}{2}\mathbb{C}}\wedge [-1,1]_+,
   \end{equation}
   which gives an $S^1$-equivariant based homeomorphism \begin{equation}\label{eq:path_space_homeomorphism}
        X_{\widehat\gamma,h}\cong_{S^1}X_{\gamma,h}
   \end{equation}
  which is the identity on $A_{\gamma,h}$. 

If $r=0$, then $X_{\widehat{\gamma},h}$ is obtained from $A_{\gamma,h}$ by identifying $S^{\frac{w(k_{-1})-h}{2}\mathbb{C}}\subset X_{\gamma_-,h}$ with $S^{\frac{w(k_1)-h}{2}\mathbb{C}}\subset X_{\gamma_+,h}$ via the identity map. Therefore, the map \begin{equation}\label{eq:collapsing_map}
    X_{\gamma,h}=A_{\gamma,h}\cup \left(S^{\frac{w(k_{\pm 1})-h}{2}\mathbb{C}}\wedge [-1,1]_+ \right)\longrightarrow A_{\gamma,h}/{\sim}=X_{\widehat{\gamma},h}
\end{equation}
collapsing the $[-1,1]_+$ direction of the cylinder $S^{\frac{w(k_{\pm 1})-h}{2}\mathbb{C}}\wedge [-1,1]_+ $ is clearly a based $S^1$-homotopy equivalence. \end{proof}

Let \[
\iota_\gamma: \mathcal{A}_\gamma \longrightarrow \mathcal{H}_{\operatorname{O}(2)}(\gamma,[k])
\]
denote the morphism in $\mathfrak{C}_{\operatorname{O}(2)}$ induced by the inclusion $A_{\gamma,h}\hookrightarrow X_{\gamma,h}$. Also, for each $i\in\{\pm1,\dots,\pm n\}$, let \[
j_i: S^{\frac{w(k_i)}{2}\mathbb{C}} \longrightarrow \operatorname{Res}^{\operatorname{O}(2)}_{S^1} \mathcal{H}_{\operatorname{O}(2)}(\gamma,[k])
\]
denote the morphism in $\mathfrak{C}_{S^1}$ induced by the inclusion $S^{\frac{w(k_i)-h}{2}\mathbb{C}}\hookrightarrow X_{\gamma,h}$. 

\begin{proposition}\label{prop:almost_I_comparison}
There exist morphisms \[
\mathcal{T}_0: \mathcal{A}_\gamma \longrightarrow \operatorname{SWF}_{\operatorname{O}(2)}(Y,\mathfrak{s},\tau_Y), \quad \mathcal{T}_\gamma: \operatorname{Res}^{\operatorname{O}(2)}_{S^1} \mathcal{H}_{\operatorname{O}(2)}(\gamma,[k]) \longrightarrow \operatorname{SWF}_{S^1}(Y,\mathfrak{s})
\]
in $\mathfrak{C}_{\operatorname{O}(2)}$ and $\mathfrak{C}_{S^1}$, respectively, such that \begin{equation}\label{eq:almost_I_comparison}
    \mathcal{T}_\gamma\circ \operatorname{Res}^{\operatorname{O}(2)}_{S^1}\iota_\gamma = \operatorname{Res}^{\operatorname{O}(2)}_{S^1}\mathcal{T}_0, \quad \mathcal{T}_\gamma\circ j_i = \operatorname{BF}_{W,\mathfrak{t}_i} \quad (i\in\{\pm1,\dots,\pm n\})
\end{equation}
in $\mathfrak{C}_{S^1}$. If $\gamma$ carries the lattice homology of $(\Gamma,[k])$, then $\mathcal{T}_\gamma$ is an $S^1$-stable equivalence.
\end{proposition}
\begin{proof} 
We will obtain the required morphisms from ordinary maps on representatives. In the following, we take suspensions implicitly when needed, suppressing the suspension factors. Also, we denote by $\operatorname{SWF}(Y)$ a based $\operatorname{O}(2)$-space representing $\operatorname{SWF}_{\operatorname{O}(2)}(Y,\mathfrak{s},\tau_Y)$. 

Recall from \eqref{eq:almost_I_path_space} that \[
X_{\gamma,h}=A_{\gamma,h}\cup \left(S^{\frac{w(k_{\pm1})-h}{2}\mathbb{C}}\wedge[-1,1]_+\right)= X_{\gamma_-,h}\cup \left(S^{\frac{w(k_{\pm1})-h}{2}\mathbb{C}}\wedge[-1,1]_+\right) \cup X_{\gamma_+,h}.
\]
We first define $f_+:X_{\gamma_+,h}\to \operatorname{SWF}(Y)$ as in Section \ref{subsec:spectra_comparison}. More precisely, recall that \[
X_{\gamma_+,h}= \bigcup_{i=1}^n S^{\frac{w(k_i)-h}{2}\mathbb{C}} \cup \bigcup_{i=1}^{n-1} \left( S^{\frac{w(e_i)-h}{2}\mathbb{C}}\wedge [0,1]_+\right),
\]
where $e_i$ is the edge $[k_i,k_{i+1}]$. We define $f_+$ (after suitable stabilizations) so that the following hold: \begin{itemize}[leftmargin=2em]
    \item On the vertex sphere $S^{\frac{w(k_i)-h}{2}\mathbb{C}}$, $f_+$ represents $\operatorname{BF}_{W,\mathfrak{t}_i}$.
    \item On the edge cylinder $S^{\frac{w(e_i)-h}{2}\mathbb{C}}\wedge [0,1]_+$, $f_+$ is a based $S^1$-homotopy representing the adjunction relation $\operatorname{BF}_{W,\mathfrak{t}_i} \circ \, U^{\frac{w(k_i)-w(e_i)}{2}} = \operatorname{BF}_{W,\mathfrak{t}_{i+1}} \circ \, U^{\frac{w(k_{i+1})-w(e_i)}{2}}$.
\end{itemize}
Next, define $f_-:X_{\gamma_-,h}\to \operatorname{SWF}(Y)$ by \[
f_-:=I_Y\circ f_+\circ I,
\]
where $I$ is the involution on $X_{\gamma,h}$ defined above, and $I_Y$ is the involution on $\operatorname{SWF}(Y)$ coming from the $\operatorname{O}(2)$-action. Recall the compatibility relation \eqref{eq:BF_naturality} from Section \ref{subsec:compatibility}: \[
\operatorname{BF}_{W,\mathfrak{t}_{-i}} = I_Y\circ\operatorname{BF}_{W,\mathfrak{t}_i}\circ \,I.
\]
Writing \[
X_{\gamma_-,h}= \bigcup_{i=1}^n S^{\frac{w(k_{-i})-h}{2}\mathbb{C}} \cup \bigcup_{i=1}^{n-1} \left( S^{\frac{w(e_{-i})-h}{2}\mathbb{C}}\wedge [0,1]_+\right),
\]
where $e_{-i}$ is the edge $[k_{-(i+1)},k_{-i}]$, it follows that the following also hold for $f_-$: \begin{itemize}[leftmargin=2em]
    \item On the vertex sphere $S^{\frac{w(k_{-i})-h}{2}\mathbb{C}}$, $f_-$ represents $\operatorname{BF}_{W,\mathfrak{t}_{-i}}$.
    \item On the edge cylinder $S^{\frac{w(e_{-i})-h}{2}\mathbb{C}}\wedge [0,1]_+$, $f_-$ is a based $S^1$-homotopy representing the adjunction relation $\operatorname{BF}_{W,\mathfrak{t}_{-i}} \circ \, U^{\frac{w(k_{-i})-w(e_{-i})}{2}} = \operatorname{BF}_{W,\mathfrak{t}_{-(i+1)}} \circ \, U^{\frac{w(k_{-(i+1)})-w(e_{-i})}{2}}$.
\end{itemize}
Gluing $f_-$ and $f_+$ gives a based $\operatorname{O}(2)$-map $f:A_{\gamma,h}\to \operatorname{SWF}(Y)$. We denote the induced morphism $\mathcal{A}_\gamma\to \operatorname{SWF}_{\operatorname{O}(2)}(Y,\mathfrak{s},\tau_Y)$ in $\mathfrak{C}_{\operatorname{O}(2)}$ by $\mathcal{T}_0$. 

We extend $f$ to $\tilde{f}:X_{\gamma,h}\to \operatorname{SWF}(Y)$ by defining $\tilde{f}$ on the cylinder $S^{\frac{w(k_{\pm1})-h}{2}\mathbb{C}}\wedge[-1,1]_+$ to be a based $S^1$-homotopy representing the relation $\operatorname{BF}_{W,\mathfrak{t}_{-1}}=\operatorname{BF}_{W,\mathfrak{t}_1}$ of Lemma \ref{lem:central_adjunction}. Thus, $\tilde{f}$ is a based $S^1$-map. We denote the induced morphism $\operatorname{Res}^{\operatorname{O}(2)}_{S^1}\mathcal{H}_{\operatorname{O}(2)}(\gamma,[k])\to \operatorname{SWF}_{S^1}(Y,\mathfrak{s})$ in $\mathfrak{C}_{S^1}$ by $\mathcal{T}_{\gamma}$. The relations \eqref{eq:almost_I_comparison} automatically hold by construction. 

Suppose now that $\gamma$ carries lattice homology, and choose a completion $\widehat{\gamma}$. We use the notation in the proof of Lemma \ref{lem:completion}. First, assume that $B\neq\varnothing$. Recall the $S^1$-homeomorphism $X_{\hat{\gamma},h}\cong_{S^1} X_{\gamma,h}$ given in \eqref{eq:path_space_homeomorphism}. Composing with $\tilde{f}$ gives a based $S^1$-map $\hat{f}:X_{\widehat{\gamma},h}\to \operatorname{SWF}(Y)$. From our construction of $\tilde{f}$, we see that $\hat{f}$ is precisely the map constructed in Section \ref{subsec:spectra_comparison}, i.e., the following hold: \begin{itemize}[leftmargin=2em]
    \item On each vertex sphere in $X_{\widehat{\gamma},h}$, $\hat{f}$ represents the relative Bauer--Furuta map.
    \item On each edge cylinder in $X_{\widehat{\gamma},h}$, $\hat{f}$ represents the adjunction relation.
\end{itemize}
Hence $\hat{f}$ induces an $S^1$-stable equivalence by Theorem \ref{thm:lattice_SWF_equivalence}. Consequently, the morphism $\mathcal{T}_\gamma$ induced by $\tilde{f}$ is also an $S^1$-stable equivalence. 

Finally, assume that $B=\varnothing$. Let $q:X_{\gamma,h}\to X_{\widehat\gamma,h}$ be the $S^1$-homotopy equivalence in \eqref{eq:collapsing_map}. Write the restriction of $\tilde{f}$ to the central cylinder $S^{\frac{w(k_{\pm1})-h}{2}\mathbb{C}}\wedge[-1,1]_+$ as $H(z,t)$. Consider the homotopy \[
H_s(z,t):=H(z,-1+(1-s)(1+t)) \quad (s\in [0,1]),
\]
which satisfies $H_0=H$, $H_s(z,-1)=H(z,-1)$, and $H_1(z,t)=H(z,-1)$. Since the inclusion $S^{\frac{w(k_{\pm 1})-h}{2}\mathbb{C}}\times \{1\}\hookrightarrow X_{\gamma_+,h}$ is an $S^1$-cofibration, the homotopy $H_s$ restricted to $S^{\frac{w(k_{\pm 1})-h}{2}\mathbb{C}}\times \{1\}$ extends to a homotopy $\tilde{H}_s$ over $X_{\gamma_+,h}$ with $\tilde{H}_0=f_+$. Now $f_-\cup H_s \cup \tilde{H}_s$ gives a based $S^1$-homotopy from $\tilde{f}$ to a map $\tilde{f}_1$ whose restriction to the central cylinder is a constant homotopy. Hence $\tilde{f}_1$ factors through $q$ and induces a based $S^1$-map $\hat{f}:X_{\widehat{\gamma},h}\to \operatorname{SWF}(Y)$ such that $\tilde{f}_1=\hat{f}\circ q$. Again, $\hat{f}$ induces an $S^1$-stable equivalence by Theorem \ref{thm:lattice_SWF_equivalence}. It follows that $\tilde{f}$ also induces an $S^1$-stable equivalence, completing the proof. 
\end{proof}

\subsection{Equivariant algebraic topology}
\label{subsec:equivariant_algebraic_topology}

Before giving the proof of Theorem \ref{thm:main}, we recall
some facts from equivariant algebraic topology. First, the following result from Smith theory is well-known:
\begin{theorem}[{\cite[Theorem III.4.23 and Remark III.4.30]{tom_Dieck-1987}}]\label{thm:Smith_fixed_points}
Let $P$ be a finite $2$-group and let $X$ be a based $P$-space having the homotopy type of a finite based $P$-CW complex. Suppose that $\widetilde{H}^*(X;\mathbb{F}_2)\cong \widetilde{H}^*(S^n;\mathbb{F}_2)$ for some integer $n\geq0$. Then, for every subgroup $K\leq P$, there exists an integer
$n(K)$ with $0\leq n(K)\leq n$ such that $\widetilde{H}^*(X^K;\mathbb{F}_2) \cong \widetilde{H}^*(S^{n(K)};\mathbb{F}_2)$. 
\end{theorem}

Next, we recall Borel's dimension formula. 

\begin{theorem}[{\cite[Theorem III.4.40]{tom_Dieck-1987}}]\label{thm:Borel_dimension_formula}
Let $X$ be a based $(\mathbb{Z}_2\times\mathbb{Z}_2)$-space having the homotopy type of a finite based $(\mathbb{Z}_2\times\mathbb{Z}_2)$-CW complex. Suppose that $\widetilde{H}^*(X;\mathbb{F}_2)\cong \widetilde{H}^*(S^n;\mathbb{F}_2)$ for some integer $n\geq0$.
For each subgroup $K$ of $\mathbb{Z}_2\times\mathbb{Z}_2$, let $n(K)$ be the integer given by Theorem \textup{\ref{thm:Smith_fixed_points}}. Then
\[
n(\{1\})-n(\mathbb{Z}_2\times\mathbb{Z}_2) =\sum_{|K|=2} \left(n(K)-n(\mathbb{Z}_2\times\mathbb{Z}_2)\right).
\]
\end{theorem}

Note that $\mathbb{Z}_2\times \mathbb{Z}_2$ has exactly three distinct subgroups of order $2$. 

Finally, we have the following result on Borel cohomology. The argument is essentially the same as that in the proof of \cite[Lemma 3.6]{Kang-Park-Taniguchi-2026}.

\begin{proposition}\label{prop:Borel_cohomology_single_degree}
Let $X$ be a based $\mathbb{Z}_2$-space having the homotopy type of a finite based $\mathbb{Z}_2$-CW complex. Suppose that $\widetilde{H}^*(X;\mathbb{F}_2)\cong \widetilde{H}^*(S^n;\mathbb{F}_2)$ for some integer $n\geq0$. Then there is an isomorphism of graded $\mathbb{F}_2[W]$-modules \begin{equation}\label{eq:Borel_cohomology}
    \widetilde{H}_{\mathbb{Z}_2}^*(X;\mathbb{F}_2) \cong (\mathbb{F}_2[W])_{[n]} \quad (|W|=1),
\end{equation}
where the subscript $[n]$ means grading shift by $n$, so that the element $1$ has degree $n$. Consequently, if $X$ is of type $\mathbb{Z}_2$-$\operatorname{SWF}$, then $d(X)=n$.
\end{proposition}
Here, $d(X)$ is the integer defined in Section \ref{subsec:type_SWF}.
\begin{proof}
For the universal $\mathbb{Z}_2$-bundle $E\mathbb{Z}_2\to B\mathbb{Z}_2$, consider the following bundle of pairs: \[
(X,\{*\}) \longrightarrow (E\mathbb{Z}_2\times_{\mathbb{Z}_2}X,\, E\mathbb{Z}_2\times_{\mathbb{Z}_2}\{*\}) \longrightarrow B\mathbb{Z}_2.
\]
Applying the relative Serre spectral sequence gives \[
E_2^{p,q}=H^p(B\mathbb{Z}_2;H^q(X,\{*\};\mathbb{F}_2)) \Longrightarrow \widetilde H_{\mathbb{Z}_2}^{p+q}(X;\mathbb{F}_2).
\]
By hypothesis, we have \[
E_2^{p,q} = \begin{cases}
H^p(B\mathbb{Z}_2;\mathbb{F}_2) &\text{if $q=n$},\\
0 & \text{if $q\neq n$}.
\end{cases}
\]
The spectral sequence therefore degenerates at the $E_2$ page, proving the isomorphism \eqref{eq:Borel_cohomology}; see \cite[Proposition III.1.18]{tom_Dieck-1987}. The last statement about $d(X)$ follows directly from its definition.
\end{proof} 

\subsection{Proof of Theorem \texorpdfstring{\ref{thm:main}}{1.1}}\label{subsec:main_proof}

We now complete the proof of Theorem \ref{thm:main}. Let \[
\gamma=(k_{-n},\dots,k_{-1},k_1,\dots,k_n)
\] be an almost $I$-invariant path carrying the lattice homology of $(\Gamma,[k])$. Fix a cutoff $h\in \mathcal{A}_{[k]}$ as in Section \ref{subsec:O2_path_lattice_spectrum}. By \eqref{eq:O(2)-stable_equivalence} and Lemma \ref{lem:O2_SWF_fixed_points}, we have \[
\operatorname{SWF}_{\operatorname{O}(2)}(Y,\mathfrak{s},\tau_Y) \simeq_{\operatorname{O}(2),\mathrm{st}} (F,a,b,c),
\]
where: \begin{itemize}[leftmargin=2em]
    \item $F$ has the homotopy type of a finite based $\operatorname{O}(2)$-CW complex, $a,b\in \mathbb{Z}_{\geq 0}$, and $c\in \mathbb{Q}$.
    \item The $S^1$-action is free on $F-F^{S^1}$, and there is a based $\langle I\rangle$-homotopy equivalence $F^{S^1}\simeq_{\langle I\rangle} S^{a\mathbb{R}\oplus b\widetilde{\mathbb{R}}}$. 
    \item $F^I$ is of type $\mathbb{Z}_2$-SWF at level $a$.
\end{itemize}
Moreover, $c+w(k)/2$ is an integer since \[
c\equiv n(Y,\mathfrak{s},g) \equiv - \frac{k^2-\sigma(W)}8 = -\frac{w(k)}{2} \pmod{\mathbb{Z}}
\]
by taking $X=W$ and $c_1(\mathfrak{t})=k$ in \eqref{eq:correction_term}. Since $h\in \mathcal{A}_{[k]}=w(k)+2\mathbb{Z}$, it follows that $c+h/2$ is also an integer. 

Observe that, for any integers $p,q,r\geq 0$, the suspension $F':=\Sigma^{p\mathbb{R}\oplus q\widetilde{\mathbb{R}}\oplus r\mathbb{C}}F$ satisfies: \begin{itemize}[leftmargin=2em]
    \item $\operatorname{SWF}_{\operatorname{O}(2)}(Y,\mathfrak{s},\tau_Y) \simeq_{\operatorname{O}(2),\mathrm{st}} (F',a+p,b+q,c+r)$.
    \item The $S^1$-action is free on $F'-(F')^{S^1}$ and $(F')^{S^1}\simeq_{\langle I\rangle} S^{(a+p)\mathbb{R}\oplus (b+q)\widetilde{\mathbb{R}}}$. 
    \item $(F')^I$ is of type $\mathbb{Z}_2$-SWF at level $a+p$.
\end{itemize}
Therefore, by taking suspensions if necessary, we may assume that $a,b$, and $c+h/2$ are all nonnegative integers. Furthermore, we may assume that the morphisms in Proposition \ref{prop:almost_I_comparison} are represented by maps into $F$. More precisely, letting $V:=a\mathbb{R}\oplus b\widetilde{\mathbb{R}}\oplus (c+\frac{h}{2})\mathbb{C}$, the proof of Proposition \ref{prop:almost_I_comparison} gives a based $\operatorname{O}(2)$-map $f:S^V\wedge A_{\gamma,h}\to F$ and a based $S^1$-map $\tilde{f}:S^V\wedge X_{\gamma,h}\to F$ such that $\tilde{f}|_{S^V\wedge A_{\gamma,h}}=f$, and $\tilde{f}$ induces an $S^1$-stable equivalence.

By Lemma \ref{lem:central_quotient}, we have \[
\frac{S^V\wedge X_{\gamma,h}}{S^V\wedge A_{\gamma,h}} \cong_{\operatorname{O}(2)} S^V\wedge \frac{X_{\gamma,h}}{A_{\gamma,h}} \cong_{\operatorname{O}(2)} S^V\wedge S^{\frac{w(k_{\pm 1})-h}{2}\mathbb{C}\oplus \widetilde{\mathbb{R}}} \cong_{\operatorname{O}(2)} S^{a\mathbb{R}\oplus (b+1)\widetilde{\mathbb{R}}\oplus (c+\frac{w(k_{\pm1})}{2})\mathbb{C}},
\]
and $(S^V\wedge A_{\gamma,h})^I \cong (S^V)^I\wedge (A_{\gamma,h})^I = \{*\}$. Now consider the reduced mapping cone  \[
\Theta:=\operatorname{Cone}(f:S^V\wedge A_{\gamma,h}\longrightarrow F),
\]
equipped with its induced $\operatorname{O}(2)$-action. Then we get the following homotopy-commutative diagram of $S^1$-equivariant cofiber sequences: \[
\begin{tikzcd}
S^V\wedge A_{\gamma,h} \ar[r, hook] \ar[d, equal] & S^V\wedge X_{\gamma,h} \ar[d, "\tilde{f}"] \ar[r] & S^{a\mathbb{R}\oplus (b+1)\widetilde{\mathbb{R}}\oplus (c+\frac{w(k_{\pm1})}{2})\mathbb{C}} \ar[d] \ar[r] & \cdots \\
S^V\wedge A_{\gamma,h} \ar[r, "f"] & F \ar[r] & \Theta \ar[r] &\cdots .
\end{tikzcd}
\]
Since $\tilde{f}$ induces an $S^1$-stable equivalence, it follows that \begin{equation}\label{eq:stable_equivalence_of_Theta}
    \Theta \simeq_{S^1,\mathrm{st}} S^{a\mathbb{R}\oplus (b+1)\widetilde{\mathbb{R}}\oplus (c+\frac{w(k_{\pm1})}{2})\mathbb{C}}. 
\end{equation}
Let $H$ denote the subgroup $\langle -1,I\rangle$ of $\operatorname{O}(2)=S^1\rtimes \langle I\rangle$, which is isomorphic to $\mathbb{Z}_2\times \mathbb{Z}_2$. 

\begin{lemma}\label{lem:fixed_point_sets} The following statements hold.
\begin{enumerate}[label=\textup{(\alph*)}, leftmargin=2em] 
    \item $\widetilde{H}^*(\Theta;\mathbb{F}_2) \cong \widetilde{H}^*(S^{a+b+2c+w(k_{\pm 1})+1};\mathbb{F}_2)$. 
    \item $\widetilde{H}^*(\Theta^{S^1};\mathbb{F}_2) \cong \widetilde{H}^*(S^{a+b+1};\mathbb{F}_2)$.
    \item $\Theta^{\langle -1\rangle}=\Theta^{S^1}$; hence $\widetilde{H}^*(\Theta^{\langle -1\rangle};\mathbb{F}_2) \cong \widetilde{H}^*(S^{a+b+1};\mathbb{F}_2)$.
    \item $\Theta^I=F^I$.
    \item $\widetilde{H}^*(\Theta^H;\mathbb{F}_2) \cong \widetilde{H}^*(S^{a};\mathbb{F}_2)$.
    \item $\Theta^I\cong \Theta^{-I}$. 
\end{enumerate}
\end{lemma}
In part (f), $-I$ denotes the product $-1\cdot I$ in $\operatorname{O}(2)$, which is different from $I^{-1}=I$. 
\begin{proof}
    The $S^1$-stable equivalence \eqref{eq:stable_equivalence_of_Theta} implies that there exists an $S^1$-representation $W= \mathbb{R}^p\oplus \mathbb{C}^q$ for some integers $p,q\geq 0$ such that there is a based $S^1$-homotopy equivalence \[
    \Sigma^W\Theta \simeq_{S^1} \Sigma^W S^{a\mathbb{R}\oplus (b+1)\widetilde{\mathbb{R}}\oplus (c+\frac{w(k_{\pm1})}{2})\mathbb{C}}\cong_{S^1} S^{(a+b+p+1)\mathbb{R}\oplus (c+q+\frac{w(k_{\pm1})}{2})\mathbb{C}}.
    \]
    Taking $\widetilde{H}^*(-;\mathbb{F}_2)$ proves part (a). Taking $S^1$-fixed-point sets gives a based homotopy equivalence \[
    \Sigma^{p\mathbb{R}} \Theta^{S^1} \simeq  S^{(a+b+p+1)\mathbb{R}},
    \]
    which proves part (b). 

    Next, recall that both $F$ and $S^V\wedge A_{\gamma,h}$ have free $S^1$-actions outside their $S^1$-fixed-point sets. Thus, the same property holds for the reduced mapping cone $\Theta$. This proves part (c). Taking the $I$-fixed-point set of $\Theta$ gives part (d); recall that $(S^V\wedge A_{\gamma,h})^I=\{*\}$. 

    For part (e), recall that $F^{S^1}\simeq_{\langle I\rangle} S^{a\mathbb{R}\oplus b\widetilde{\mathbb{R}}}$. Therefore, taking the $H$-fixed-point set of $\Theta$, and using $F^{\langle -1\rangle}=F^{S^1}$ give \[
    \Theta^H=F^H=(F^{\langle -1\rangle})^{I} = (F^{S^1})^{I} \simeq (S^{a\mathbb{R}\oplus b\widetilde{\mathbb{R}}})^I =S^{a\mathbb{R}},
    \]
    which proves part (e). 

    Finally, for part (f), note that the element $i\in S^1$ satisfies $iIi^{-1}=-I$. Thus the action of $i$ gives a based homeomorphism $\Theta^I\cong \Theta^{-I}$. 
\end{proof}

\begin{corollary}\label{cor:cohomology_of_F_I}
     $\widetilde{H}^*(F^I;\mathbb{F}_2)\cong \widetilde{H}^*(S^{a+c+\frac{w(k_{\pm 1})}{2}};\mathbb{F}_2)$.
\end{corollary}
\begin{proof}
    For each subgroup $K$ of $H=\langle -1,I\rangle \cong \mathbb{Z}_2\times \mathbb{Z}_2$, let $n(K)\geq 0$ denote the integer determined by $\widetilde{H}^*(\Theta^K;\mathbb{F}_2) \cong \widetilde{H}^*(S^{n(K)};\mathbb{F}_2)$; see Theorem \ref{thm:Smith_fixed_points}. Parts (a), (c), (e), and (f) of Lemma \ref{lem:fixed_point_sets} give \begin{align*}
        n(\{1\})&=a+b+2c+w(k_{\pm1})+1,\\
     n(\langle -1\rangle)&=a+b+1,\\
      n(\langle I\rangle)&=n(\langle -I\rangle),\\
       n(H)&=a.
    \end{align*}
    On the other hand, Theorem \ref{thm:Borel_dimension_formula} gives \[
n(\{1\})-n(H)=n(\langle -1\rangle)+n(\langle I\rangle)+n(\langle -I\rangle)-3n(H).
\]
Substituting the above values, we obtain \[
n(\langle I\rangle)=a+c+\frac{w(k_{\pm1})}{2}.
\]
Combining with Lemma \ref{lem:fixed_point_sets}(d) completes the proof. 
\end{proof}

Recall now that the space $F^I$ is of type $\mathbb{Z}_2$-SWF at level $a$. Therefore, Proposition \ref{prop:Borel_cohomology_single_degree} and Corollary \ref{cor:cohomology_of_F_I} give \[
d(F^I)=a+c+\frac{w(k_{\pm 1})}{2}.
\]
On the other hand, since $\operatorname{SWF}_{\operatorname{O}(2)}(Y,\mathfrak{s},\tau_Y) \simeq_{\operatorname{O}(2),\mathrm{st}} (F,a,b,c)$, applying the functor $\Phi^I$ gives \[
\operatorname{SWF}_R(Y,\mathfrak{s},\tau_Y)
\simeq_{\mathbb Z_2,\mathrm{st}}
(F^I,a,c).
\]
It follows that \begin{align*}
    \delta_R(Y,\mathfrak{s},\tau_Y) & = \frac{1}{2}d(F^I,a,c)=\frac{1}{2}(d(F^I)-a-c) \\
    & = \frac{w(k_{\pm1})}{4}=\frac{(k_{\pm1})^2-\sigma(W_\Gamma)}{16}.
\end{align*}
This completes the proof of Theorem \ref{thm:main}.

\section{Proof of Theorem \texorpdfstring{\ref{thm:main_existence}}{1.2}}\label{sec:proof_main_existence}

In this section, we prove Theorem \ref{thm:main_existence}. The arguments are purely combinatorial and essentially based on N\'emethi's work \cite{Nemethi-2005}, adapted to our equivariant setting.

\subsection{Increasing paths}\label{subsec:increasing_paths}

Let $\Gamma$ and $v_0$ be as in Theorem \ref{thm:main_existence}. We use the notation $L$, $L^*$, and $q$ from Section \ref{subsec:lattice_homology}, and write $Q:=Q_\Gamma$. For $x=\sum_v x_vv\in L$, we denote by $x_0:=x_{v_0}$ the coefficient of the vertex $v_0$. We define a partial order $\leq$ on $L$ by declaring $x\leq y$ if $x_v\leq y_v$ for every vertex $v$. Write $x\wedge y$ and $x\vee y$ for the coordinatewise minimum and maximum, respectively.

Fix a characteristic covector $\kappa\in\operatorname{Char}(\Gamma)$, and define $\chi_{\kappa}:L\to \mathbb{Z}$ by \[
\chi_\kappa(x):= -\frac{Q(x,x)+\kappa(x)}2.
\]
The map $L\to [\kappa]=\kappa+2q(L)$, $x\mapsto \kappa+2q(x)$, is a bijection. Moreover, the weight function $w$ satisfies \begin{equation}\label{eq:weight_function_and_chi}
    w(\kappa+2q(x))=w(\kappa)-2\chi_\kappa(x) \quad (x\in L).
\end{equation}
Thus, minimizing $\chi_\kappa$ on $L$ corresponds to maximizing $w$ on $[\kappa]$.

Let $K$ denote the canonical characteristic covector. We claim that \begin{equation}\label{eq:canonical_ineq}
\chi_K(z)\geq 1 \quad\text{for every nonzero effective }z\in L \text{ with }z_0\in\{0,1\}.
\end{equation}
To see this, let $\Gamma'$ be a rational graph obtained from $\Gamma$ by decreasing the weight of $v_0$ by an integer $d\geq 0$. Let $Q'$ and $K'$ be its intersection form and canonical characteristic covector, respectively. Then
\[
\chi_{K'}^{\Gamma'}(z)-\chi_K(z) = \frac{d}{2}z_0(z_0-1).
\]
The right-hand side vanishes when $z_0=0$ or $1$, and \eqref{eq:canonical_ineq} follows from the rationality of $\Gamma'$.

We will also use the inequality \begin{equation}\label{eq:chi_ineq}
\chi_\kappa(x\wedge y)+\chi_\kappa(x\vee y) \leq \chi_\kappa(x)+\chi_\kappa(y) \quad(x,y\in L).
\end{equation}
To see this, first observe that
\[
\chi_\kappa(x)+\chi_\kappa(y) -\chi_\kappa(x\wedge y)-\chi_\kappa(x\vee y) = Q(x-x\wedge y,y-x\wedge y).
\]
Both $x-x\wedge y$ and $y-x\wedge y$ are effective and have disjoint supports, so the right-hand side is nonnegative. (Here, the \emph{support} of an element $x=\sum_v x_v v\in L$ is the set of all vertices $v$ such that $x_v\neq 0$.) 

We say that a finite sequence $(z_0,\ldots,z_N)$ in $L$ is an \emph{increasing path} if each successive difference $z_{j+1}-z_j$ is a vertex of $\Gamma$. Under the above bijection $L\cong [\kappa]$, an increasing path in $L$ gives a lattice path in $[\kappa]$ in the sense of Section \ref{subsec:carrying_paths}. We say that a finite sequence $(z_0,\ldots,z_N)$ in $L$ is an \emph{edge path} if each successive difference $z_{j+1}-z_j$ is $\pm v$ for some vertex $v$. We say that a subset $\mathcal{S}\subset L$ is \emph{connected} if, for every $u,v\in \mathcal{S}$, there exists an edge path contained in $\mathcal{S}$ from $u$ to $v$. We have the following elementary observation:

\begin{lemma}\label{lem:increasing_path}
Let $\mathcal{S}\subset L$ be a connected finite set which is closed under $\wedge$ and $\vee$. If $u,v\in\mathcal{S}$ with $u\leq v$, then there is an increasing path from $u$ to $v$ contained in $\mathcal{S}$.
\end{lemma}
\begin{proof}
Choose an edge path $(u=z_0,z_1,\ldots,z_M=v)$ in $\mathcal{S}$, and put \[
\widetilde{z}_j:=(z_j\vee u)\wedge v, \quad y_j:=\bigvee_{t=0}^j\widetilde{z}_t \quad (j=0,\dots,M).
\]
All these points belong to $\mathcal{S}$. For each $j=0,\dots,M$, we have $u\leq \widetilde{z}_j\leq v$. Also, each difference $\widetilde{z}_{j+1}-\widetilde{z}_j$ is either zero or $\pm w$ for some vertex $w$. Consequently, every successive difference $y_{j+1}-y_j$ is either zero or a vertex of $\Gamma$. Since $y_0=u$ and $y_M=v$, removing repetitions gives the required increasing path.
\end{proof}

For each $i\in\mathbb{Z}$, define \[
L_i:=\{x\in L:x_0=i\}.
\]
Then define \[
\mu_i:=\min_{x\in L_i}\chi_\kappa(x), \quad \mathcal{M}_i:= \{x\in L_i:\chi_\kappa(x)=\mu_i\}.
\]
Since $Q$ is negative definite, $\mu_i$ is finite, and $\mathcal{M}_i$ is a nonempty finite set.

\begin{lemma}\label{lem:minimum} The following statements hold.
\begin{enumerate}[label=\textup{(\alph*)}, leftmargin=2em]
\item For every $c\in\mathbb{Z}$, $\{x\in L_i:\chi_\kappa(x)\leq c\}$ is a connected subset of $L$.
\item The set $\mathcal{M}_i$ is closed under $\wedge$ and $\vee$. In particular, it has a (unique) maximal element $m_i$ with respect to $\leq$. Every $x\in\mathcal{M}_i$ can be connected to $m_i$ by an increasing path contained in $\mathcal{M}_i$.
\item We have $m_i+v_0\leq m_{i+1}$.
\end{enumerate} \end{lemma}
\begin{proof}
For part (a), let $\Gamma_1,\ldots,\Gamma_s$ be the components of $\Gamma-\{v_0\}$. The inequality \eqref{eq:canonical_ineq} implies that each $\Gamma_j$ is a rational graph. Let $L(\Gamma_j)$ denote the sublattice of $L$ generated by the vertices of $\Gamma_j$. For each $j=1,\dots,s$, let $\kappa_{j}\in \operatorname{Char}(\Gamma_j)$ be the characteristic covector such that \[
\kappa_{j}(v):=\kappa(v)+2iQ(v_0,v) \quad \text{for every vertex $v$ of $\Gamma_j$},
\]
and fix $z_j\in L(\Gamma_j)$ minimizing $\chi_{\kappa_{j}}^{\Gamma_j}$. Now let $x\in L_i$ with $\chi_\kappa(x)\leq c$. Writing $x=iv_0+y_1+\cdots+y_s$ with $y_j\in L(\Gamma_j)$ gives \[
\chi_\kappa(x) = \chi_\kappa(iv_0) + \sum_{j=1}^s \chi_{\kappa_{j}}^{\Gamma_j}(y_j)
\]
and $\chi_{\kappa_{j}}^{\Gamma_j}(z_j)\leq \chi_{\kappa_{j}}^{\Gamma_j}(y_j)$. By \cite[Corollary 4.1.5]{Nemethi-2008}, the set \[
A_j:=\left\{ y\in L(\Gamma_j):\chi_{\kappa_j}^{\Gamma_j}(y)\leq \chi_{\kappa_j}^{\Gamma_j}(y_j)\right\}
\]
is a connected subset of $L(\Gamma_j)$. Thus, there is an edge path from $y_j$ to $z_j$ contained in $A_j$. These paths give an edge path from $x=iv_0+y_1+\cdots+y_s$ to $iv_0+z_1+\cdots+z_s$ contained in $\{\chi_\kappa \leq c\}\cap L_i$. Therefore, we have shown that every point in $\{\chi_\kappa \leq c\}\cap L_i$ can be joined to $iv_0+z_1+\cdots+z_s$ by an edge path contained in $\{\chi_\kappa \leq c\}\cap L_i$, proving part (a). 

For part (b), let $x,y\in\mathcal{M}_i$. Both $x\wedge y$ and $x\vee y$ belong to $L_i$, so
\eqref{eq:chi_ineq} implies $x\wedge y,x\vee y\in\mathcal{M}_i$. Since $\mathcal{M}_i$ is finite and closed under $\vee$, it has a maximal element $m_i$. Taking $c=\mu_i$ shows that $\mathcal{M}_i$ is connected. Therefore, the last statement of part (b) follows from Lemma \ref{lem:increasing_path}.

For part (c), note that $m_i\wedge m_{i+1}\in L_i$ and $m_i\vee m_{i+1}\in L_{i+1}$. Hence \eqref{eq:chi_ineq} gives \begin{align*}
    \mu_i+\mu_{i+1} &\leq \chi_\kappa(m_i\wedge m_{i+1}) + \chi_\kappa(m_i\vee m_{i+1})\leq \chi_\kappa(m_i)+\chi_\kappa(m_{i+1}) = \mu_i+\mu_{i+1}.
\end{align*}
In particular, $m_i\vee m_{i+1}\in\mathcal{M}_{i+1}$. The maximality of $m_{i+1}$ then implies $m_i\vee m_{i+1}=m_{i+1}$, and therefore $m_i\leq m_{i+1}$. Their $v_0$-coefficients are $i$ and $i+1$, respectively, so we get $m_i+v_0\leq m_{i+1}$.
\end{proof}

Next, the following lemma constructs increasing paths from $m_i\in \mathcal{M}_i$ to $m_{i+1}\in \mathcal{M}_{i+1}$ with precise control of $\chi_\kappa$. It is an analog of the results in \cite[Lemmas 7.7 and 9.1]{Nemethi-2005}.

\begin{lemma}\label{lem:paths_minimum}
For every $i\in\mathbb{Z}$, we have $\chi_\kappa(m_i+v_0)=\mu_{i+1}$. Furthermore, there exists an increasing path \[
(m_i=z_0,z_1,\dots,z_N=m_{i+1})
\] such that $z_1=m_i+v_0$ and $\chi_\kappa(z_j)=\mu_{i+1}$ for $j=1,\dots,N$. 
\end{lemma}
\begin{proof}
Put $\ell_i:=\kappa+2q(m_i)\in \operatorname{Char}(\Gamma)$. For every vertex $v\neq v_0$, we have $m_i+v\in L_i$, so $\chi_\kappa(m_i+v)-\chi_\kappa(m_i)\geq0$. Equality would imply that $m_i+v\in\mathcal{M}_i$, contradicting the maximality of $m_i$. Therefore, $\chi_\kappa(m_i+v)-\chi_\kappa(m_i)\geq 1$. Since \[
\chi_\kappa(m_i+v)-\chi_\kappa(m_i)=-\frac{\ell_i(v)+Q(v,v)}2,
\]
we obtain $\ell_i(v)\leq -Q(v,v)-2$, i.e., \begin{equation}\label{eq:canonical_bound}
\ell_i(v)\leq K(v) \quad(v\neq v_0).
\end{equation}

By Lemma \ref{lem:minimum}(c), we can write $m_{i+1}=m_i+v_0+y$ with $y\geq 0$ and $y_0=0$. It follows that \begin{align*}
    \chi_\kappa(m_{i+1})-\chi_\kappa(m_i+v_0) &= -\frac{Q(y,y)+\ell_i(y)+2Q(v_0,y)}2\\
 & \geq -\frac{Q(y,y)+K(y)+2Q(v_0,y)}2= \chi_K(v_0+y)-1\geq0.
\end{align*}
The first inequality follows from \eqref{eq:canonical_bound}, since $y$ is effective and its support does not contain $v_0$. The last inequality follows from \eqref{eq:canonical_ineq}. On the other hand, since $m_i+v_0\in L_{i+1}$, we have $\chi_{\kappa}(m_i+v_0)\geq \chi_{\kappa}(m_{i+1})$. Thus, we must have $\chi_{\kappa}(m_i+v_0)=\chi_{\kappa}(m_{i+1})=\mu_{i+1}$. In particular, $m_i+v_0\in\mathcal{M}_{i+1}$. Now Lemma \ref{lem:minimum}(b) gives an increasing path from $m_i+v_0$ to $m_{i+1}$ contained in $\mathcal{M}_{i+1}$. Adjoining $m_i$ gives the required path.
\end{proof}

Recall the bijection $L\to [\kappa]$, $x\mapsto\kappa+2q(x)$. Lemma \ref{lem:paths_minimum} gives a lattice path from $\kappa+2q(m_i)$ to $\kappa+2q(m_{i+1})$ whose minimum weight is \[
\min \{ w(\kappa)-2\chi_\kappa(m_i), \, w(\kappa)-2\chi_\kappa(m_{i+1})\}=w(\kappa)-2\max \{\mu_i,\mu_{i+1}\}.
\]

\subsection{Construction of the central pair}\label{subsec:construction_of_central_pair}

Fix a Real $\operatorname{spin}^c$ structure $\mathfrak{s}$ for $(Y_\Gamma,\tau_Y)$, and let $[k]\subset \operatorname{Char}(\Gamma)$ be the corresponding orbit. Then $I([k])=[k]$, where $I$ is the involution on $\operatorname{Char}(\Gamma)$ given by $I(\ell)=-\tau^*\ell$; see Section \ref{subsec:almost_I_invariant_paths}. In this subsection, we construct an almost $I$-invariant path $(\mathfrak{t}_{-1},\mathfrak{t}_1)$ with some additional properties. In the next subsection, we will extend this pair to a carrying almost $I$-invariant path.

Since $\tau$ permutes the vertices of $\Gamma$, $\tau_*:L\to L$ permutes the vertex basis. Also, it preserves $Q$ and satisfies $\tau^*q(x)=q(\tau_*x)$ for $x\in L$. We first choose a representative of $[k]$ adapted to the involution $\tau$.

\begin{lemma}\label{lem:normalization}
There exist $\kappa\in[k]$ and a (possibly empty) set $B$ of pairwise nonadjacent $\tau$-fixed vertices such that \[
I(\kappa)=\kappa+2q\left(\sum_{v\in B}v \right).
\]
\end{lemma}
\begin{proof} Fix any $\kappa'\in [k]$. Since $I([k])=[k]$, there is a unique $a\in L$ such that $I(\kappa')=\kappa'+2q(a)$. Then \[
\kappa'=I^2(\kappa')=I(\kappa'+2q(a))=I(\kappa')-2\tau^* q(a)=(\kappa'+2q(a))-2q(\tau_*a),
\]
which gives $\tau_*a=a$. Write $a=\sum_v a_vv$. For each $\tau$-fixed vertex $v$, put $t_v:=\lfloor a_v/2\rfloor$. For each pair of vertices exchanged by $\tau_*$, pick one $v$ and put $t_v=a_v$ and $t_{\tau_*v}=0$. Now define \[
b:=a-\sum_{v} t_v v - \tau_*  \sum_v t_vv =\sum_v (a_v-t_v-t_{\tau_*v})v.
\]
Then \begin{equation}\label{eq:b_v}
    b_v:=a_v-t_v-t_{\tau_*v} =\begin{cases}
    \text{0 or 1} & \text{if $\tau_*v=v$}, \\
    0 & \text{otherwise}.
\end{cases}
\end{equation}
Set $\kappa:=\kappa'+2q(\sum_v t_vv)$ and let $B:=\{v:b_v=1\}$. Then we have $I(\kappa)=\kappa+2q(b)$. It remains to show that the vertices in $B$ are pairwise nonadjacent. Note that $\kappa+\tau^*\kappa=-2q(b)$. For every $v\in B$, we have $\tau_*v=v$, and hence \[
\kappa(v)=\frac{1}{2}(\kappa(v)+(\tau^*\kappa)(v))=-Q(b,v) =-Q(v,v)-\deg_B(v),
\]
where $\deg_B(v)$ is the number of vertices in $B$ adjacent to $v$. Since $\kappa$ is characteristic, $\deg_B(v)$ must be even. This means that every vertex in the subgraph $\Gamma_B$ of $\Gamma$ induced by $B$ has even degree. Since $\Gamma$ is a tree, we conclude that $\Gamma_B$ has no edges, completing the proof. 
\end{proof}

Fix $\kappa$ and $B$ as in Lemma \ref{lem:normalization}, and put $b:=\sum_{v\in B} v$. We use the notation $\chi_\kappa$, $L_i$, $\mu_i$, and $\mathcal{M}_i$ from Section \ref{subsec:increasing_paths} for this choice of $\kappa$. Under the bijection $L\to [\kappa]$,  $x\mapsto\kappa+2q(x)$, the involution $I$ corresponds to the involution $\iota:L\to L$, $\iota(x)=b-\tau_*x$. Since $I$ preserves the weight function $w$, \eqref{eq:weight_function_and_chi} gives \begin{equation}\label{eq:chi_I_invariance}
\chi_\kappa(\iota(x))=\chi_\kappa(x) \quad (x\in L).
\end{equation}
Moreover, $\iota$ reverses the partial order $\leq$, i.e., $x\leq y$ implies $\iota(y)\leq \iota(x)$. It follows that $\iota(x\wedge y)=\iota(x)\vee \iota(y)$ and $\iota(x\vee y)=\iota(x)\wedge \iota(y)$.

Now, let $b_0\in \{0,1\}$ be the $v_0$-coefficient of $b$. Since $v_0$ is a $\tau$-fixed vertex, we have \begin{equation}\label{eq:i_symmetry}
\iota(L_i)=L_{b_0-i},\quad \iota(\mathcal{M}_i)=\mathcal{M}_{b_0-i}, \quad \mu_i=\mu_{b_0-i} \quad (i\in\mathbb{Z}).
\end{equation}
\begin{lemma}\label{lem:choosing_x}
There exists $x\in\mathcal{M}_0$ such that $\iota(x)=x+b\in\mathcal{M}_{b_0}$. 
\end{lemma}
\begin{proof} By \eqref{eq:i_symmetry}, $\mathcal{M}:=\mathcal{M}_0\cup\mathcal{M}_{b_0}$ is a finite $\iota$-invariant set, and we have $\chi_\kappa=\mu_0$ on $\mathcal{M}$. Moreover, the inequality \eqref{eq:chi_ineq} implies that $\mathcal{M}$ is closed under $\wedge$ and $\vee$.

We claim that $\mathcal{M}$ is a connected subset of $L$. By Lemma \ref{lem:minimum}(a), both $\mathcal{M}_0$ and $\mathcal{M}_{b_0}$ are connected. Thus, $\mathcal{M}$ is connected if $b_0=0$. If $b_0=1$, then Lemma \ref{lem:paths_minimum} gives $m_0+v_0\in\mathcal{M}_1$, so the edge path $(m_0,m_0+v_0)$ joins $\mathcal{M}_0$ and $\mathcal{M}_1$. Therefore, $\mathcal{M}=\mathcal{M}_0\cup \mathcal{M}_1$ is connected, as claimed.

For any $y\in \mathcal{M}$, we have $y\wedge \iota(y)\in \mathcal{M}$, and $y\wedge\iota(y) \leq y\vee\iota(y)=\iota(y\wedge\iota(y))$. Thus, we can choose $x\in\mathcal{M}$ with $x\leq\iota(x)$ minimizing
\[
D(x):=\sum_v(\iota(x)-x)_v.
\]
Put $d:=\iota(x)-x=b-\tau_*x-x\geq 0$. Then $\tau_*d=d$, and $d_v\equiv b_v \pmod{2}$ for every $\tau$-fixed vertex $v$.

By Lemma \ref{lem:increasing_path}, there is an increasing path $\gamma$ from $x$ to $\iota(x)$ contained in $\mathcal{M}$. For every point $y$ on this path, we have $x\leq y\leq\iota(x)$, and hence $x\leq\iota(y)\leq\iota(x)$. Therefore, $z:=y\wedge\iota(y)$ belongs to $\mathcal{M}$ and satisfies $x\leq z\leq\iota(z)$. Since $\tau_*$ permutes the vertices, $z\neq x$ implies \[
D(z)=\sum_v b_v-2\sum_v z_v < \sum_v b_v-2\sum_v x_v=D(x),
\]
which contradicts the minimality of $D(x)$. We conclude that $y\wedge\iota(y)=x$ for every point $y$ on the path.

We now claim that \[
d_v=\begin{cases}
       \text{0 or 1} & \text{if $\tau_*v=v$}, \\
    0 & \text{otherwise}.
\end{cases}
\]
Since $d_v\equiv b_v \pmod{2}$ for every $\tau$-fixed vertex $v$, this claim, together with \eqref{eq:b_v}, will imply that $d=b$, which gives $\iota(x)=x+b$, or equivalently, $\tau_*x=-x$. Since $v_0$ is $\tau$-fixed, it follows that $x_0=0$, and hence $x\in \mathcal{M}_0$ and $\iota(x)=x+b\in \mathcal{M}_{b_0}$, completing the proof.

It remains to prove the claim. Write the increasing path $\gamma$ as \[
\gamma=(x=y_0,y_1,\dots,y_N=\iota(x)).
\]
First, let $v$ be a $\tau$-fixed vertex and suppose $d_v\geq 2$. Let $i\geq 1$ be the minimal index such that $(y_i-x)_v=1$. Then \[
\iota(y_i)_v=(b-\tau_*(y_i))_v=b_v-(y_i)_v=(d_v+2x_v)-(x_v+1)\geq x_v+1,
\]
which contradicts that $y_i\wedge \iota(y_i)=x$. Since $d\geq 0$, we conclude that $d_v\in \{0,1\}$.

Next, let $v$ be a vertex not fixed by $\tau_*$, and suppose $d_v=d_{\tau_*v}>0$. Let $i\geq 1$ be the minimal index such that either $(y_i-x)_v=1$ or $(y_i-x)_{\tau_*v}=1$. Interchanging $v$ and $\tau_*v$ if necessary, we may assume that $(y_i)_v=x_v+1$ and $(y_i)_{\tau_*v}=x_{\tau_*v}$. Since $b_v=0$, it follows that \[
\iota(y_i)_v=-(y_i)_{\tau_*v}=-x_{\tau_*v}=d_v+x_v\geq x_v+1,
\]
again contradicting $y_i\wedge \iota(y_i)=x$. This proves the claim. \end{proof}

Choose $x$ as in Lemma \ref{lem:choosing_x}, and put \[
k_{-1}:=\kappa+2q(x),\quad k_1:=\kappa+2q(x+b).
\]
Then we have \[
I(k_{-1})=k_1,\quad k_1-k_{-1}=2\sum_{v\in B}v^*.
\]
Moreover, \eqref{eq:weight_function_and_chi} gives \begin{equation}\label{eq:central_weights}
w(k_{-1})=w(k_1)=w(\kappa)-2\mu_0=w(\kappa)-2\mu_{b_0}.
\end{equation}
Let $\mathfrak{t}_{-1}$ and $\mathfrak{t}_1$ denote the $\operatorname{spin}^c$ structures on $W_\Gamma$ corresponding to $k_{-1}$ and $k_1$, respectively. They restrict to $\mathfrak{s}$ and satisfy \begin{equation}\label{eq:t_1-t_-1}
    \mathfrak{t}_{-1}=\tau^*\overline{\mathfrak{t}}_1, \quad \mathfrak{t}_1-\mathfrak{t}_{-1} =\sum_{v\in B}v^*.
\end{equation}
The spheres corresponding to the vertices in $B$ are pairwise disjoint and are preserved setwise and orientation-preservingly by $\tau$. Thus $(\mathfrak{t}_{-1},\mathfrak{t}_1)$ is an almost $I$-invariant path in the sense of Section \ref{subsec:almost_I_invariant_paths}. To prove Theorem \ref{thm:main_existence}, it remains to extend this path to a carrying almost $I$-invariant path.

We end this subsection by noting the following properties of $k_{\pm 1}$: \begin{equation}\label{eq:properties_of_k_pm1}
    w(k_{-1}) =\max_{y\in L_0}\,w(\kappa+2q(y)), \quad w(k_1) =\max_{y\in L_{b_0}}\,w(\kappa+2q(y)),
\end{equation}
which follow from \eqref{eq:weight_function_and_chi} and \eqref{eq:central_weights}.

\subsection{Extension to a carrying path}\label{subsec:extension}

In this subsection, we complete the proof of Theorem \ref{thm:main_existence} by extending the pair $(\mathfrak{t}_{-1},\mathfrak{t}_1)$ constructed above to a carrying almost $I$-invariant path. We use the notation in Section \ref{subsec:construction_of_central_pair}.

For $c\in\mathbb{Z}$, let $G_c$ denote the graph whose vertices are the elements of \[
\{z\in L:\chi_\kappa(z)\leq c\},
\]
with an edge between two vertices whenever their difference is $\pm v$ for some vertex $v$ of $\Gamma$. We prove the following lemma using the argument in the proof of \cite[Proposition 4.3]{Nemethi-2005}.

\begin{lemma}\label{lem:finite_sublevel}
The set \[
\mathcal{C}:=\left\{ z\in L: \chi_\kappa(z)\leq\chi_\kappa(z\pm v) \text{ for every vertex $v$ of $\Gamma$} \right\}
\]
is finite and meets every connected component of $G_c$ for every $c\in\mathbb{Z}$.
\end{lemma}
\begin{proof} For any $z\in L$, put $\ell_z:=\kappa+2q(z)$. Note that $z\in \mathcal{C}$ if and only if \[
Q(v,v)\leq\ell_z(v)\leq-Q(v,v)
\]
for every vertex $v$ of $\Gamma$. Thus, $\ell_z(v)$ has only finitely many choices, which implies that $\mathcal{C}$ is finite. 

Next, let $C$ be a connected component of $G_c$. Since $Q$ is negative definite, $G_c$ has only finitely many vertices. Choose a vertex $z\in C$ minimizing $\chi_\kappa$ on $C$. If $\chi_\kappa(z\pm v)<\chi_\kappa(z)$ for some vertex $v$, then $z\pm v$ belongs to $G_c$ and is joined to $z$ by an edge. It therefore belongs to $C$, contradicting the choice of $z$. We conclude that $z\in\mathcal{C}$.
\end{proof}

Since $\mathcal{C}$ is finite, we can choose an integer $N\geq1$ such that \begin{equation}\label{eq:finite_slice_range}
b_0-N\leq z_0\leq N \quad (z\in \mathcal{C}).
\end{equation}
By Lemma \ref{lem:choosing_x}, $x+b\in \mathcal{M}_{b_0}$. Lemma \ref{lem:minimum}(b) gives an increasing path from $x+b$ to $m_{b_0}$ contained in $\mathcal{M}_{b_0}$. For each $i=b_0,\dots,N-1$, choose an increasing path from $m_i$ to $m_{i+1}$ as in Lemma \ref{lem:paths_minimum}. Concatenating all these paths gives an increasing path \[
P_+=(u_1=x+b,u_2,\dots,u_n=m_N).
\]
By construction, for each $j=1,\dots,n$, we have $u_j\in \mathcal{M}_{(u_j)_0}$, so $\chi_\kappa(u_j)=\mu_{(u_j)_0}$. 

Next, put $u_{-j}:=\iota(u_j)$ for $j=1,\dots,n$. In particular, $u_{-1}=x$ by Lemma \ref{lem:choosing_x}. Since $\iota(y+v)=\iota(y)-\tau_*v$ for any $y\in L$ and vertex $v$, it follows that \[
P_-:=(u_{-n},\dots,u_{-1}=x)
\]
is also an increasing path. 

For $j\in\{\pm1,\dots,\pm n\}$, let $\mathfrak{t}_j$ be the $\operatorname{spin}^c$ structure on $W_\Gamma$ such that $c_1(\mathfrak{t}_j)=\kappa+2q(u_j)$. For $j=\pm 1$, this agrees with $\mathfrak{t}_{\pm 1}$ constructed in Section \ref{subsec:construction_of_central_pair}.
\begin{lemma}
    $\gamma:=(\mathfrak{t}_{-n},\dots,\mathfrak{t}_{-1},\mathfrak{t}_1,\dots,\mathfrak{t}_n)$ is an almost $I$-invariant path for the Real spin$^c$ structure $\mathfrak{s}$.
\end{lemma}
\begin{proof}
    For every $j\in \{\pm 1,\dots,\pm n\}$, we have $\mathfrak{t}_j|_{Y_\Gamma}=\mathfrak{s}$ since $c_1(\mathfrak{t}_j)\in [\kappa]$, and $[\kappa]$ corresponds to $\mathfrak{s}$. For $j=1,\dots,n$, we have $\mathfrak{t}_{-j}=\tau^*\overline{\mathfrak{t}}_j$ since \[
    c_1(\tau^*\overline{\mathfrak{t}}_j)=-\tau^*c_1(\mathfrak{t}_j)=-\tau^*(\kappa+2q(u_j))=I(\kappa+2q(u_j))=\kappa+2q(\iota(u_j))=c_1(\mathfrak{t}_{-j}).
    \]
    By construction, $(c_1(\mathfrak{t}_{-n}),\dots,c_1(\mathfrak{t}_{-1}))$ and $(c_1(\mathfrak{t}_1),\dots,c_1(\mathfrak{t}_n))$ are lattice paths in the sense of Section \ref{subsec:carrying_paths}. Finally, the condition on $\mathfrak{t}_1-\mathfrak{t}_{-1}$ was established in \eqref{eq:t_1-t_-1}. 
\end{proof}

Choose a completion $\widehat{\gamma}$ of $\gamma$ (Section \ref{subsec:almost_I_invariant_paths}), and let $P$ be the corresponding increasing path in $L$ under the bijection $L\to [\kappa]$, $y\mapsto\kappa+2q(y)$. Let us also regard $P$ and $P_\pm$ as the graphs consisting of their vertices and the edges between successive vertices. Note that $P_\pm$ are subgraphs of $P$. 

\begin{lemma}\label{lem:graph_P} The following statements hold. \begin{enumerate}[label=\textup{(\alph*)}, leftmargin=2em]
    \item $\chi_\kappa(y)=\mu_{y_0}$ for every vertex $y$ of $P$.
    \item $P$ meets $L_i$ for every $b_0-N\leq i\leq N$.
\end{enumerate}
\end{lemma}
\begin{proof}
    We have already seen that $\chi_\kappa(y)=\mu_{y_0}$ for every vertex $y$ of $P_+$. The same holds for $P_-$ by \eqref{eq:chi_I_invariance} and \eqref{eq:i_symmetry}, since every vertex of $P_-$ is of the form $\iota(y)$ for some vertex $y$ of $P_+$. Write \[
    P=(\underbrace{u_{-n},\dots,u_{-1}=x}_{P_-},\dots,\underbrace{x+b=u_1,\dots,u_n}_{P_+}).
    \]
    If $y$ is a vertex between $x$ and $x+b$, Remark \ref{rmk:central_cube} implies \[
    w(\kappa+2q(y))=w(\kappa+2q(x))=w(\kappa+2q(x+b)).
    \]
    Thus, $\chi_\kappa(y)=\mu_0=\mu_{b_0}$ by \eqref{eq:weight_function_and_chi}. Since the $v_0$-coefficient $y_0$ of $y$ is either $0$ or $b_0$, we conclude that $\chi_\kappa(y)=\mu_{y_0}$, proving part (a). 

    For part (b), observe that $P$ starts in $L_{b_0-N}$ and ends in $L_N$. Since $P$ is an increasing path, it must meet every $L_i$ with $b_0-N\leq i\leq N$.
\end{proof}

\begin{lemma}\label{lem:P_c_to_G_c}
   For each $c\in \mathbb{Z}$, $P_c:=P\cap G_c\hookrightarrow G_c$ induces a bijection on connected components.
\end{lemma}
\begin{proof}
    Let $C$ be a connected component of $G_c$. By Lemma \ref{lem:finite_sublevel}, we can choose $z\in C\cap\mathcal{C}$. Then $b_0-N\leq z_0\leq N$ by \eqref{eq:finite_slice_range}, so $P$ contains a vertex $y\in L_{z_0}$ by Lemma \ref{lem:graph_P}(b). Then by Lemma \ref{lem:graph_P}(a), we have $\chi_\kappa(y)=\mu_{y_0}=\mu_{z_0}\leq\chi_\kappa(z)\leq c$. Now Lemma \ref{lem:minimum}(a) gives an edge path from $y$ to $z$ contained in $\{\chi_\kappa \leq c\}\cap L_{z_0}$. Thus $y\in C$, proving surjectivity.
    
Next, suppose two vertices $p,q$ of $P_c$ are contained in the same connected component of $G_c$. We may assume $p$ occurs before $q$ in the path $P$, so that $p_0\leq q_0$. An edge path from $p$ to $q$ in $G_c$ must meet every $L_i$ with $p_0\leq i\leq q_0$. Thus, $\mu_i\leq c$ for $p_0\leq i\leq q_0$. Now note that $p_0\leq y_0\leq q_0$ for every vertex $y$ of the subpath of $P$ from $p$ to $q$. Therefore, $\chi_\kappa(y)=\mu_{y_0}\leq c$ by Lemma \ref{lem:graph_P}(a). It follows that this subpath lies in $P_c$, proving injectivity.
\end{proof}

Now we give proofs of Theorem \ref{thm:main_existence} and Corollary \ref{cor:Froyshov_formula}.

\begin{proof}[Proof of Theorem \textup{\ref{thm:main_existence}}]
    Under the bijection $L\to [\kappa]$, $y\mapsto\kappa+2q(y)$, \eqref{eq:weight_function_and_chi} shows that the condition $\chi_\kappa(y)\leq c$ corresponds to $w(\kappa+2q(y))\geq w(\kappa)-2c$. Now recall from Section \ref{subsec:lattice_homology} the superlevel set description of lattice homology. The bijection $L\to [\kappa]$ restricts to a bijection between the vertices of $G_c$ and the vertices of the superlevel set $S_{w(\kappa)-2c}$. Hence $G_c$ is identified with the $1$-skeleton $S_{w(\kappa)-2c}^{(1)}$ of $S_{w(\kappa)-2c}$. Similarly, the path $P$ in $L$ is identified with the lattice path $\widehat{\gamma}$, and therefore $P_c$ is identified with the intersection $S_{w(\kappa)-2c}^{(1)}\cap \widehat{\gamma}$. In summary, we have a commutative diagram \[
\begin{tikzcd}
    P_c \ar[d, "\cong"'] \ar[r, hook] & G_c \ar[d, "\cong"] \\ 
    S_{w(\kappa)-2c}^{(1)}\cap \widehat{\gamma} \ar[r, hook] & S_{w(\kappa)-2c}^{(1)},
\end{tikzcd}
\]
where the vertical maps are identifications induced by the bijection $L\to [\kappa]$. As $c$ ranges over $\mathbb{Z}$, $w(\kappa)-2c$ ranges over $\mathcal{A}_{[\kappa]}$. Note that $S_{w(\kappa)-2c}^{(1)}\cap \widehat{\gamma}=S_{w(\kappa)-2c}\cap \widehat{\gamma}$, since $\widehat{\gamma}$ is $1$-dimensional. Therefore, Lemma \ref{lem:P_c_to_G_c} implies that the inclusion $S_a\cap \widehat{\gamma}\hookrightarrow S_a$ induces a bijection on $\pi_0$ for every $a\in \mathcal{A}_{[\kappa]}$. Lemma \ref{lem:carrying_criterion} now implies that the lattice path $\widehat{\gamma}$ carries the lattice homology of $(\Gamma,[\kappa])$. This means that $\gamma$ is a carrying almost $I$-invariant path; see Section \ref{subsec:almost_I_invariant_paths}. This completes the proof.
\end{proof}

\begin{proof}[Proof of Corollary \textup{\ref{cor:Froyshov_formula}}]
Choose $\kappa$ and $B$ as in Lemma \ref{lem:normalization}. Put $b:=\sum_{v\in B}v$. Then $\kappa+\tau^*\kappa=-2q(b)$, so $I(\kappa)+\tau^*I(\kappa)=2q(b)$, proving existence. For uniqueness, suppose $B'$ is another such set. Put $b':=\sum_{v\in B'}v$, and choose $\ell\in[k]$ such that $\ell+\tau^*\ell =2q(b')$. Writing $\ell=\kappa+2q(y)$ gives $b+b'=y+\tau_*y$. Therefore, the coefficient $b_v+b'_v$ is even for every $\tau$-fixed vertex $v$. Since $b_v,b'_v\in\{0,1\}$, we must have $b_v=b'_v$. Since $b$ and $b'$ are supported on $\tau$-fixed vertices, we conclude that $b=b'$; hence $B=B'$. This proves that $B_{\mathfrak{s}}=B$ is uniquely determined.

Next, let $\ell\in\mathcal{K}_{[k]}$. Then $\ell+\tau^*\ell =2q(b)$. In particular, this implies that $b$ is a Wu class of $\Gamma_{\mathrm{fix}}$ as in Remark \ref{rmk:Wu_class}. Writing $\ell=\kappa+2q(y)$, we get $y+\tau_*y=2b$. Since the vertex $v_0$ is $\tau$-fixed, this implies $y_0=b_0$. Therefore, \eqref{eq:properties_of_k_pm1} gives \[
w(\ell) \leq \max_{z\in L_{b_0}}\,w(\kappa+2q(z)) = w(k_1).
\]
This upper bound is attained since $k_1+\tau^*k_1=k_1-k_{-1}=2q(b)$. Theorem \ref{thm:main} now gives \[
\delta_R(Y_\Gamma,\mathfrak{s},\tau_Y) =\frac{k_1^2-\sigma(W_\Gamma)}{16} = \frac{w(k_1)}4 = \max_{\ell\in\mathcal{K}_{[k]}} \frac{\ell^2+|\Gamma|}{16}.
\]

It remains to show that $\det Q_{\Gamma_{\mathrm{fix}}}$ is odd if $Y_\Gamma$ is a $\mathbb{Z}_2$-homology sphere. If $Y_\Gamma$ is a $\mathbb{Z}_2$-homology sphere, then the mod $2$ intersection form \[
Q_2:H_2(W_\Gamma;\mathbb{Z}_2)\times H_2(W_\Gamma;\mathbb{Z}_2)\longrightarrow \mathbb{Z}_2
\]
is nondegenerate. Since the homomorphism $\operatorname{id}+\tau_*$ on $H_2(W_\Gamma;\mathbb{Z}_2)$ is self-adjoint with respect to $Q_2$, it follows that \[
\left( \ker (\operatorname{id}+\tau_*)\right)^\perp = \operatorname{im}(\operatorname{id}+\tau_*) \subset H_2(W_\Gamma;\mathbb{Z}_2).
\]
Note that $(\operatorname{id}+\tau_*)^2=0$, so $\operatorname{im}(\operatorname{id}+\tau_*)\subset \ker (\operatorname{id}+\tau_*)$. Therefore, $Q_2$ induces a nondegenerate form \[
\overline{Q}_2: \frac{\ker (\operatorname{id}+\tau_*)}{\operatorname{im}(\operatorname{id}+\tau_*)}\times \frac{\ker (\operatorname{id}+\tau_*)}{\operatorname{im}(\operatorname{id}+\tau_*)}\longrightarrow \mathbb{Z}_2.
\]
On the other hand, $\ker (\operatorname{id}+\tau_*)$ is spanned by the $\tau$-fixed vertices and the sums of the form $v+\tau_*v$ for non-fixed vertices $v$, while $\operatorname{im}(\operatorname{id}+\tau_*)$ is spanned by the sums $v+\tau_*v$. Therefore, the $\tau$-fixed vertices form a $\mathbb{Z}_2$-basis of $\ker(\operatorname{id}+\tau_*)/\operatorname{im}(\operatorname{id}+\tau_*)$, and $\overline{Q}_2$ is precisely the reduction of $Q_{\Gamma_{\mathrm{fix}}}$ modulo $2$. We conclude that $\det Q_{\Gamma_{\mathrm{fix}}}$ is odd.
\end{proof}

\section{Even torus knots}\label{sec:even_torus_knots}

We apply Corollary \ref{cor:Froyshov_formula} to compute the Real Fr{\o}yshov invariants of double branched covers of even torus knots. We use the notation in Section \ref{sec:lattice_spectra}. Throughout this section, let $p\geq 1$ and $q>1$ be integers with $q$ odd and $\gcd (p,q)=1$. Consider the Brieskorn manifold \[
Y_{p,q}:=\Sigma(2,2p,q)=\left\{(x,y,z)\in \mathbb{C}^3:x^2+y^{2p}+z^q=0\right\}\cap S^5,
\]
oriented as the link of the corresponding complex surface singularity. It is well-known that $Y_{p,q}$ is the double cover of $S^3$ branched along the torus knot $T_{2p,q}$, where the deck transformation $\tau_{p,q}$ is given by $(x,y,z)\mapsto (-x,y,z)$; see \cite[Lemma 1.1]{Milnor-1975}. 

\subsection{Plumbing graphs}\label{subsec:plumbing_graphs}

We first express $Y_{p,q}$ as the boundary of a negative definite plumbing; see \cite{Neumann-Raymond-1978} and \cite[Section 4.6]{Kang-Park-Taniguchi-2026}. If $p>1$, let $\alpha$ and $\beta$ be the unique integers satisfying \[
0<\alpha<p,\quad 0<\beta<q,\quad q\alpha\equiv-1\pmod{p},\quad 2p\beta\equiv-1\pmod{q},
\]
and set \[
c:=\frac{1+q\alpha+2p\beta}{pq}.
\]
Then $c\in\{1,2\}$ and \begin{equation}\label{eq:torus_seifert}
-c+\frac{\alpha}{p}+\frac{2\beta}{q}=-\frac1{pq}.
\end{equation}
Suppose that the negative continued fractions are given by \begin{equation}\label{eq:continued_fraction}
    \frac{p}{\alpha}=[a_1,\ldots,a_m]^-:=a_1-\frac{1}{a_2-\dfrac{1}{\cdots -\dfrac{1}{a_m}}}, \quad \frac{q}{\beta}=[b_1,\ldots,b_n]^-:=b_1-\frac{1}{b_2-\dfrac{1}{\cdots -\dfrac{1}{b_n}}} \quad (a_i,b_j\geq 2).
\end{equation}
When $p=1$, we set $\alpha=0$, $c=1$, $\beta=(q-1)/2$, and $m=0$. ($T_{2,q}$ is a $2$-bridge knot, so $Y_{1,q}=\Sigma(2,2,q)$ is a lens space.) 

\begin{figure}[ht]
\centering
\begin{tikzpicture}[scale=1.4]
\draw (-5,0.5) node[circle, fill, inner sep=1.2pt, black]{};
\draw (-3,0.5) node[circle, fill, inner sep=1.2pt, black]{};
\draw (-2,0.5) node[circle, fill, inner sep=1.2pt, black]{};
\draw (-1,0) node[circle, fill, inner sep=1.2pt, black]{};
\draw (1,0) node[circle, fill, inner sep=1.2pt, black]{};
\draw (-1,1) node[circle, fill, inner sep=1.2pt, black]{};
\draw (1,1) node[circle, fill, inner sep=1.2pt, black]{};
\draw (-4,0.5) node{$\cdots$};
\draw (0,1) node{$\cdots$};
\draw (0,0) node{$\cdots$};
\draw (-1,0) node[above]{$-b_1$};
\draw (1,0) node[above]{$-b_n$};
\draw (-1,1) node[above]{$-b_1$};
\draw (1,1) node[above]{$-b_n$};
\draw (-2.1,0.5) node[above]{$-c$};
\draw (-3,0.5) node[above]{$-a_1$};
\draw (-5,0.5) node[above]{$-a_m$};
\draw (-3,0.5) node[below]{\textcolor{Emerald}{$u_1$}};
\draw (-5,0.5) node[below]{\textcolor{Emerald}{$u_m$}};
\draw (-2.1,0.5) node[below]{\textcolor{Emerald}{$v_0$}};
\draw (-0.9,1) node[below]{\textcolor{Emerald}{$v_1^+$}};
\draw (1,1) node[below]{\textcolor{Emerald}{$v_n^+$}};
\draw (-1,0) node[below]{\textcolor{Emerald}{$v_1^-$}};
\draw (1,0) node[below]{\textcolor{Emerald}{$v_n^-$}};
\draw (-5,0.5)--(-4.3,0.5) (-1,0)--(-0.3,0) (-3.7,0.5)--(-2,0.5) (-2,0.5)--(-1,0) (-2,0.5)--(-1,1) (-1,1)--(-0.3,1) (0.3,1)--(1,1) (0.3,0)--(1,0) ;
\draw[<->, Emerald] (2,1) to[bend left] (2,0);
\draw (2.15,0.5) node[right] {\textcolor{Emerald}{$\tau$}};
\end{tikzpicture}
\caption{The graph $\Gamma_{p,q}$ and the involution $\tau$.}\label{fig:plumbing}
\end{figure}
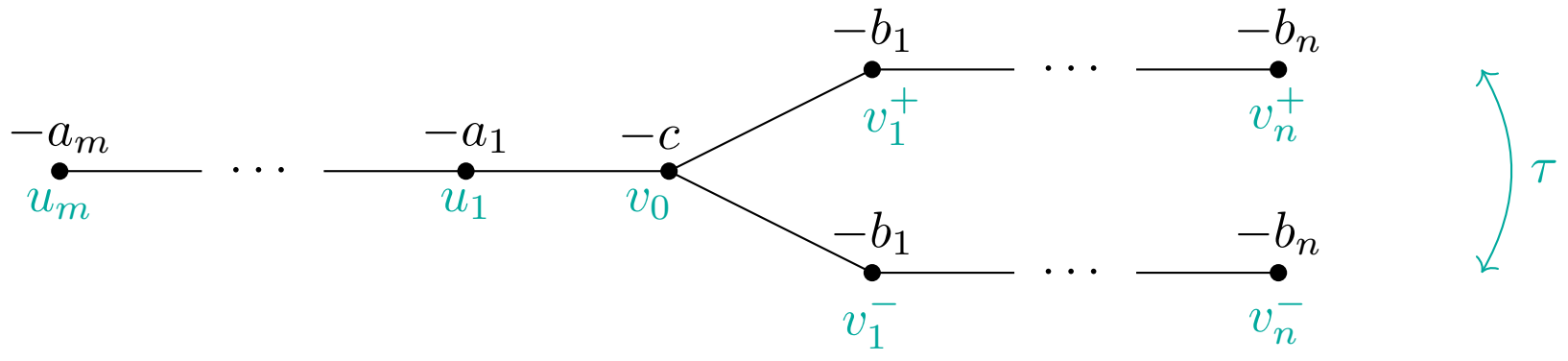

Let $\Gamma_{p,q}$ be the star-shaped graph described in Figure \ref{fig:plumbing}. It has central vertex $v_0$ of weight $-c$, one leg $u_1,\dots,u_m$ of weights $-a_1,\dots,-a_m$, and two legs $v_1^+,\dots,v_n^+$ and $v_1^-,\dots,v_n^-$, each of weights $-b_1,\dots,-b_n$. (When $p=1$, there is no $u$-leg.) This graph is negative definite, and its oriented boundary is $Y_{p,q}$. We denote the associated $4$-manifold by $W_{p,q}:=W_{\Gamma_{p,q}}$. 

The construction in \cite[Section 4.6]{Kang-Park-Taniguchi-2026} gives a smooth orientation-preserving involution $\tau$ on $W_{p,q}$ that permutes the vertex spheres with \[
\tau_*v_0=v_0,\quad \tau_*u_i=u_i,\quad \tau_*v_j^{\pm}=v_j^{\mp},
\]
such that the restriction $\tau|_{Y_{p,q}}$ is identified with the involution $\tau_{p,q}$ \cite[Lemma 4.4]{Kang-Park-Taniguchi-2026}. We can decrease the central weight $-c$ sufficiently to make the graph $\Gamma_{p,q}$ rational \cite[Example 8.2(3)]{Nemethi-2005}, so the pair $(\Gamma_{p,q},\tau)$ satisfies the assumptions of Theorem \ref{thm:main_existence}.

Let $\{\eta_v\}$ denote the Hom-dual basis of $L^*=H^2(W_{p,q};\mathbb{Z})$, so that $\eta_v(w)=\delta_{vw}$ for any two vertices $v,w$ of $\Gamma_{p,q}$. (Note that $\eta_v$ is different from $v^*=q(v)=\operatorname{PD}(v)$.) Consider the following positive definite matrices: \begin{equation}\label{eq:matrices}
    U=\begin{pmatrix}
a_1 & -1 & 0 & \cdots & 0 & 0   \\
-1&a_2&  -1& \cdots  & 0 & 0 \\
0 &  -1 & a_3 & \cdots & 0  &0 \\
\vdots &\vdots & \vdots &\ddots& \vdots & \vdots \\
0 & 0 & 0 & \cdots &a_{m-1}& -1 \\
0 & 0 & 0 & \cdots & -1 & a_m
\end{pmatrix}, \quad 
V=\begin{pmatrix}
b_1 & -1 & 0 & \cdots & 0 & 0   \\
-1&b_2&  -1& \cdots  & 0 & 0 \\
0 &  -1 & b_3 & \cdots & 0  &0 \\
\vdots &\vdots & \vdots &\ddots& \vdots & \vdots \\
0 & 0 & 0 & \cdots &b_{n-1}& -1 \\
0 & 0 & 0 & \cdots & -1 & b_n
\end{pmatrix}.
\end{equation}
Using the continued fraction \eqref{eq:continued_fraction}, define $r_0,\dots,r_{n+1}$ recursively by: \begin{equation}\label{eq:r_j}
    r_n=1,\quad r_{n+1}=0,\quad r_{j-1}=b_jr_j-r_{j+1}\quad(j=1,\dots,n).
\end{equation}
For $j=0,\dots,n-1$, if $V_j$ denotes the $(n-j)\times (n-j)$ submatrix of $V$ obtained from the last $n-j$ rows and columns, then we have $r_j=\det V_j$. In particular, $r_0=q$ and $r_1=\beta$. 

\begin{lemma}\label{lem:torus_labels}
The group $H^2(Y_{p,q};\mathbb{Z})$ is cyclic of order $q$, generated by $h:=\eta_{v_n^+}|_{Y_{p,q}}$. Also, we have \begin{equation}\label{eq:restriction_formula}
    \eta_{v_j^+}|_{Y_{p,q}}=r_jh,\quad \eta_{v_j^-}|_{Y_{p,q}}=-r_jh,\quad \eta_{v_0}|_{Y_{p,q}}=\eta_{u_i}|_{Y_{p,q}}=0.
\end{equation}
Furthermore, the canonical characteristic covector $K$ represents the unique spin structure $\mathfrak{s}_0$.  \end{lemma}

See \cite[Lemma 5.3]{BKST-2025} for an alternative proof of $H^2(Y_{p,q};\mathbb{Z})\cong \mathbb{Z}_q$.

\begin{proof} Let $\phi:L^*\to\mathbb{Z}_q$ be the homomorphism determined by \[
\phi(\eta_{v_j^\pm})=\pm r_j \quad (j=1,\dots,n), \quad \phi(\eta_{v_0})=\phi(\eta_{u_i})=0 \quad (i=1,\dots,m).
\]
We claim that $\phi$ vanishes on $q(L)$. It suffices to show that $\phi(q(v))=0$ for every vertex $v$. Let $Q:=Q_{\Gamma_{p,q}}$ denote the intersection form. Then \[
q(v)=Q(v,v)\eta_v +\sum_{\text{vertices $w$ with $Q(v,w)=1$}} \eta_w.
\]
Thus, we have: \begin{itemize}[leftmargin=2em]
    \item If $2\leq j\leq n-1$, then \[
    \phi(q(v_j^\pm))=\phi(-b_j \eta_{v_j^\pm}+\eta_{v_{j-1}^\pm}+\eta_{v_{j+1}^\pm})=\pm(-b_j r_j+r_{j-1}+r_{j+1})=0.
    \]
    \item If $n\geq 2$, then \[
    \phi(q(v_n^\pm))=\phi(-b_n \eta_{v_n^\pm}+\eta_{v_{n-1}^\pm})=\pm(-b_nr_n+r_{n-1})=0
    \]
    and \[
    \phi(q(v_1^\pm))=\phi( -b_1\eta_{v_1^\pm}+\eta_{v_0}+\eta_{v_2^\pm})=\pm(-b_1r_1+ r_2)=\mp r_0=\mp q=0. 
    \]
    \item If $n=1$, then $\beta=r_1=1$ and $b_1=q$, so \[
    \phi(q(v_1^\pm))=\phi(-b_1\eta_{v_1^\pm}+\eta_{v_0})=\mp b_1r_1=\mp q=0.
    \]
    \item If $p>1$, then \[
    \phi(q(v_0))=\phi(-c\eta_{v_0}+\eta_{u_1}+\eta_{v_1^+}+\eta_{v_1^-})=r_1-r_1=0.
    \]
    If $p=1$, then similarly $\phi(q(v_0))=0$. 
    \item $\phi(q(u_i))=0$ since $q(u_i)$ is a linear combination of $\eta_{v_0},\eta_{u_1},\dots,\eta_{u_m}$.
\end{itemize}
This proves the claim. Therefore, $\phi$ induces a homomorphism $\bar{\phi}:L^*/q(L)\to \mathbb{Z}_q$, which is surjective since $\phi(\eta_{v_n^+})=r_n=1$.

From now on, we assume $p>1$; the case $p=1$ is similar but easier. Order the vertices by \[
\{v_0,u_1,\dots,u_m,v_1^+,\dots,v_n^+,v_1^-,\dots,v_n^-\}. \]
Then the intersection matrix is given by \begin{equation}\label{eq:intersection_matrix}
    Q=-\begin{pmatrix}
    c & -\mathbf{e}_m^t & -\mathbf{e}_n^t & -\mathbf{e}_n^t \\
    -\mathbf{e}_m & U & 0 & 0 \\
    -\mathbf{e}_n & 0 & V & 0 \\
    -\mathbf{e}_n & 0 & 0 & V
\end{pmatrix},
\end{equation}
where $\mathbf{e}_m=(1,0,\dots,0)^t\in\mathbb{R}^m$, $\mathbf{e}_n=(1,0,\dots,0)^t\in \mathbb{R}^n$. Put $\mathbf{z}:=(\mathbf{e}_m^t,\mathbf{e}_n^t,\mathbf{e}_n^t)^t\in \mathbb{R}^{m+2n}$ and $D:=\operatorname{diag}(U,V,V)$ so that \[
Q=-\begin{pmatrix}
    c & -\mathbf{z}^t \\
    -\mathbf{z} & D
\end{pmatrix}.
\]
It follows that \begin{align*}
    \det (-Q)& =(c-\mathbf{z}^tD^{-1}\mathbf{z})\cdot \det D= \left(c-(U^{-1})_{11}-2(V^{-1})_{11} \right) \cdot pq^2 = \left(c-\frac{\alpha}{p}-\frac{2\beta}{q}\right)pq^2=q,
\end{align*}
where the last equality uses \eqref{eq:torus_seifert}. Hence $H^2(Y_{p,q};\mathbb{Z})\cong L^*/q(L)$ has order $|\det Q|=q$, which implies that $\bar{\phi}$ must be an isomorphism. This proves the first statement, and \eqref{eq:restriction_formula} follows from the definition of $\phi$. Finally, the canonical characteristic covector is \[
K=(c-2)\eta_{v_0} +\sum_{i=1}^m(a_i-2)\eta_{u_i} +\sum_{j=1}^n(b_j-2)  \left(\eta_{v_j^+}+\eta_{v_j^-}\right),
\]
so $K|_{Y_{p,q}}=0$ by \eqref{eq:restriction_formula}. This proves that the orbit $[K]$ corresponds to the unique spin structure $\mathfrak{s}_0$. ($Y_{p,q}$ has a unique spin structure since $|H^2(Y_{p,q};\mathbb{Z})|=q$ is odd.)
\end{proof}

Define \begin{equation}\label{eq:torus_labels}
\mathfrak{s}_i:=[K+2i\eta_{v_n^+}]=\mathfrak{s}_0+ih \in \operatorname{Spin}^c(Y_{p,q}) \quad (i\in\mathbb{Z}_q).
\end{equation}
In particular, $c_1(\mathfrak{s}_i)=2ih\in H^2(Y_{p,q};\mathbb{Z})$. By Lemma \ref{lem:torus_labels}, we have $\operatorname{Spin}^c(Y_{p,q})=\{\mathfrak{s}_i:i\in \mathbb{Z}_q\}$. Moreover, note that \begin{equation}\label{eq:action_of_tau}
    \tau^*\eta_{v_0}=\eta_{v_0}, \quad \tau^*\eta_{u_i}=\eta_{u_i},\quad \tau^*\eta_{v_j^\pm}=\eta_{v_j^\mp}.
\end{equation}
It follows that \[
\tau_{p,q}^*\mathfrak{s}_i=\overline{\mathfrak{s}}_i=\mathfrak s_{-i}.
\]
Hence every $\mathfrak{s}_i$ is a Real $\operatorname{spin}^c$ structure for $(Y_{p,q},\tau_{p,q})$. (More generally, for any knot $K\subset S^3$, every $\operatorname{spin}^c$ structure on the double branched cover is Real with respect to the deck transformation; see Section \ref{subsec:branched_covers}.) 

\subsection{Real \texorpdfstring{Fr{\o}yshov}{Froyshov} invariants}
\label{subsec:torus_formula}

In this subsection, we derive a formula for the Real Fr{\o}yshov invariants $\delta_R(Y_{p,q},\mathfrak{s}_i,\tau_{p,q})$. We use the notation in the previous subsection.

Let $\Gamma_{\mathrm{fix}}$ denote the full subgraph of $\Gamma_{p,q}$ spanned by the $\tau$-fixed vertices $v_0,u_1,\ldots,u_m$. Put \[
A:=-Q_{\Gamma_{\mathrm{fix}}} = \begin{pmatrix} 
c & -\mathbf{e}_m^t\\
-\mathbf{e}_m & U
\end{pmatrix}.
\]
(When $p=1$, we have $A=(1)$.) The matrix $A$ is positive definite, and \[
\det A = p\left(c-\frac{\alpha}{p}\right) = pc-\alpha =\frac{1+2p\beta}{q},
\]
which is an odd integer since both $1+2p\beta$ and $q$ are odd. (In fact, since $Y_{p,q}$ is a $\mathbb{Z}_2$-homology sphere by Lemma \ref{lem:torus_labels}, the proof of Corollary \ref{cor:Froyshov_formula} shows that $\det A=\pm \det Q_{\Gamma_{\mathrm{fix}}}$ is odd.) Thus, $\Gamma_{\mathrm{fix}}$ has a unique Wu class \[
\omega=\omega_0v_0+\sum_{i=1}^m \omega_i u_i\quad (\omega_0,\dots,\omega_m\in \{0,1\}),
\]
where the coefficients are determined by \[
A \begin{pmatrix}
     \omega_0 \\ \omega_1 \\ \vdots \\ \omega_m
\end{pmatrix} \equiv \begin{pmatrix}
    c \\ a_1 \\ \vdots \\ a_m
\end{pmatrix} \pmod{2}.
\]
By Corollary \ref{cor:Froyshov_formula}, we have $\omega=b_{\mathfrak{s}_i}$ for every $i\in\mathbb{Z}_q$.

Now recall \eqref{eq:continued_fraction}, \eqref{eq:matrices}, and \eqref{eq:r_j}. Put $\mathbf{r}=(r_1,\dots,r_n)^t$. For each $i\in\mathbb{Z}_q$, define \begin{equation}\label{eq:torus_finite_set}
\mathcal{D}_i := \left\{ \mathbf{z}=(z_1,\ldots,z_n)^t\in\mathbb{Z}^n:
\begin{array}{l}
|z_j|\leq b_j\quad (j=1,\dots,n),\\
z_1\equiv b_1+\omega_0 \pmod{2}, \\
z_j\equiv b_j \pmod{2} \quad (j=2,\dots,n),\\
\mathbf{r}^t\mathbf{z}\equiv i\pmod{q}
\end{array} \right\}.
\end{equation}
The conditions $|z_j|\leq b_j$ imply that this set is finite. Put $\boldsymbol{\omega}:=(\omega_0,\dots,\omega_m)^t$.

\begin{proposition}\label{prop:torus_Froyshov_formula}
For every $i\in\mathbb{Z}_q$, we have \[
\delta_R(Y_{p,q},\mathfrak{s}_i,\tau_{p,q}) =\frac{(1+m+2n)-\boldsymbol{\omega}^tA\boldsymbol{\omega}}{16}-\frac{1}{8}\min_{\mathbf{z}\in\mathcal{D}_i}\mathbf{z}^tV^{-1}\mathbf{z}.
\]
\end{proposition}
\begin{proof}
Corollary \ref{cor:Froyshov_formula} shows that \begin{equation}\label{eq:formula_from_corollary}
    \delta_R(Y_{p,q},\mathfrak{s}_i,\tau_{p,q}) = \max_\ell  \frac{\ell^2+|\Gamma_{p,q}|}{16}=\max_\ell \frac{\ell^2+(1+m+2n)}{16},
\end{equation}
where the maximum is taken over all $\ell\in \operatorname{Char}(\Gamma_{p,q})$ such that $[\ell]$ represents $\mathfrak{s}_i$ and $\ell+\tau^*\ell=2q(\omega)$. 

We first characterize (not necessarily characteristic) $\ell\in L^*$ such that $\ell+\tau^*\ell =2q(\omega)$. Since $\tau_*\omega=\omega$, this condition is equivalent to \[
\tau^*(\ell-q(\omega))=-(\ell-q(\omega)).
\]
Recall the basis $\{\eta_v\}$ for $L^*$. By \eqref{eq:action_of_tau}, this equality holds if and only if \[
\ell-q(\omega)=\sum_{j=1}^n z_j\left(\eta_{v_j^+}-\eta_{v_j^-}\right)
\]
for some $\mathbf{z}=(z_1,\dots,z_n)^t\in \mathbb{Z}^n$. Thus, letting $\ell_{\mathbf{z}}:=q(\omega)+\sum_{j=1}^n z_j(\eta_{v_j^+}-\eta_{v_j^-})$, we have \[
\left\{\ell\in L^*: \ell+\tau^*\ell =2q(\omega) \right\}=\left\{\ell_{\mathbf{z}}:\mathbf{z}\in \mathbb{Z}^n \right\}.
\]

Next, we show that $\ell_{\mathbf{z}}$ is characteristic if and only if \begin{equation}\label{eq:ell_z_characteristic}
    z_1\equiv b_1+\omega_0\pmod{2}, \quad z_j\equiv b_j\pmod{2} \quad (j=2,\dots,n).
\end{equation}
Since $\omega$ is a Wu class for $\Gamma_{\mathrm{fix}}$, we have \[
\ell_{\mathbf{z}}(v)=q(\omega)(v) \equiv v^2 \pmod{2}
\]
for every $\tau$-fixed vertex $v\in \{v_0,u_1,\dots,u_m\}$. For the other vertices, we have \[
\ell_{\mathbf{z}}(v_j^\pm) = q(\omega)(v_j^{\pm})\pm z_j= \omega_0\delta_{1j}\pm z_j.
\]
This proves \eqref{eq:ell_z_characteristic}, since $(v_j^\pm)^2=-b_j$. 

Now, note that $q(\omega)|_{Y_{p,q}}=0$ since $q(\omega)\in q(L)$. Therefore, \eqref{eq:restriction_formula} gives \[
\ell_{\mathbf{z}}|_{Y_{p,q}} = 2\left(\sum_{j=1}^n r_jz_j\right)h = 2(\mathbf{r}^t\mathbf{z})h.
\]
Since $H^2(Y_{p,q};\mathbb{Z})$ has odd order $q$, the first Chern class determines a $\operatorname{spin}^c$ structure on $Y_{p,q}$. Recalling that $c_1(\mathfrak{s}_i)=2ih$, we conclude that $\ell_{\mathbf{z}}$ represents $\mathfrak{s}_i$ if and only if \[
\mathbf{r}^t\mathbf{z}\equiv i\pmod{q}.
\]

Define \[
\widetilde{\mathcal{D}}_i := \left\{ \mathbf{z}=(z_1,\ldots,z_n)^t\in\mathbb{Z}^n:
\begin{array}{l}
z_1\equiv b_1+\omega_0 \pmod{2}, \\
z_j\equiv b_j \pmod{2} \quad (j=2,\dots,n),\\
\mathbf{r}^t\mathbf{z}\equiv i\pmod{q}
\end{array} \right\}.
\]
We have seen that there is a bijection \begin{align*}
    \widetilde{\mathcal{D}}_i &\longrightarrow \left\{\ell\in \operatorname{Char}(\Gamma_{p,q}):\text{$[\ell]$ represents $\mathfrak{s}_i$ and $\ell+\tau^*\ell=2q(\omega)$} \right\}, \\
    \mathbf{z} & \longmapsto \ell_{\mathbf{z}}.
\end{align*}
Therefore, \eqref{eq:formula_from_corollary} gives \[
\delta_R(Y_{p,q},\mathfrak{s}_i,\tau_{p,q}) = \max_{\mathbf{z}\in \widetilde{\mathcal{D}}_i} \frac{\ell_{\mathbf{z}}^2+(1+m+2n)}{16}.
\]
Observe that \[
\ell_{\mathbf{z}}^2 = \omega^2 + \bigg( \sum_{j=1}^n z_j(\eta_{v_j^+}-\eta_{v_j^-})\bigg)^2 = -\boldsymbol{\omega}^tA\boldsymbol{\omega} -2\mathbf{z}^tV^{-1}\mathbf{z}.
\]
Since the matrix $V^{-1}$ is positive definite, the function $\mathbf{z}\mapsto\mathbf{z}^tV^{-1}\mathbf{z}$ attains its minimum on $\widetilde{\mathcal{D}}_i$. It follows that \[
\delta_R(Y_{p,q},\mathfrak{s}_i,\tau_{p,q}) =\frac{(1+m+2n)-\boldsymbol{\omega}^tA\boldsymbol{\omega}}{16}-\frac{1}{8}\min_{\mathbf{z}\in \widetilde{\mathcal{D}}_i}\mathbf{z}^tV^{-1}\mathbf{z}.
\]

It remains to show that the minimum can be taken over $\mathcal{D}_i$. For $\mathbf{z}\in \widetilde{\mathcal{D}}_i$, suppose $|z_j|>b_j$ for some $j$. We show that $\mathbf{z}$ cannot minimize $\mathbf{z}^tV^{-1}\mathbf{z}$. Let $\{\mathbf{f}_1,\ldots,\mathbf{f}_n\}$ denote the standard basis of $\mathbb{Z}^n$. The relations \eqref{eq:r_j} imply $\mathbf{r}^tV=q\mathbf{f}_1^t$, so it follows that $\mathbf{z}\pm 2V\mathbf{f}_j \in \widetilde{\mathcal{D}}_i$. Moreover, we have \[
(\mathbf{z}\pm2V\mathbf{f}_j)^tV^{-1} (\mathbf{z}\pm2V\mathbf{f}_j) -\mathbf{z}^tV^{-1}\mathbf{z} = 4(\pm z_j+b_j).
\]
Since $|z_j|>b_j$, either $z_j+b_j$ or $-z_j+b_j$ is negative. This completes the proof. 
\end{proof}

\begin{remark} (1) Proposition \ref{prop:torus_Froyshov_formula} gives a finite algorithm for computing $\delta_R(Y_{p,q},\mathfrak{s}_i,\tau_{p,q})$ from $p$, $q$, and $i$. Indeed, the continued fractions determine the matrices, the Wu class is obtained by solving a linear system, and these determine the finite set $\mathcal{D}_i$.

(2) Since $\mathcal{D}_{-i}=-\mathcal{D}_i$, Proposition \ref{prop:torus_Froyshov_formula} gives
\begin{equation}\label{eq:conjugation_invariance}
    \delta_R(Y_{p,q},\mathfrak{s}_i,\tau_{p,q}) = \delta_R(Y_{p,q},\mathfrak{s}_{-i},\tau_{p,q}).
\end{equation}
Thus, in explicit calculations, it suffices to consider $0\leq i\leq (q-1)/2$. In fact, it is a general property that $\delta_R$ is invariant under conjugation: $\delta_R(Y,\bar{\mathfrak{s}},\tau) = \delta_R(Y,\mathfrak{s},\tau)$.
\end{remark}

As a sanity check, let us show that the spin case of Proposition \ref{prop:torus_Froyshov_formula} recovers the formula of \cite[Theorem 4.7]{Kang-Park-Taniguchi-2026}; the following proposition will not be used elsewhere. Recall that $\mathfrak{s}_0$ is the unique spin structure on $Y_{p,q}$.

\begin{proposition}\label{prop:spin_case}
We have \[
\delta_R(Y_{p,q},\mathfrak{s}_0,\tau_{p,q}) = -\frac{1}{2}\bar{\mu}(Y_{p,q}),
\]
where $\bar{\mu}$ denotes the Neumann--Siebenmann invariant \cite{Neumann-1980, Siebenmann-1980, Saveliev-2002}.
\end{proposition}
\begin{proof}
Let $\Omega$ denote the unique Wu class of $\Gamma_{p,q}$. Since $\tau_*\Omega$ is also a Wu class, uniqueness implies $\tau_*\Omega=\Omega$. It is easy to see that the restriction of $\Omega$ to the $\tau$-fixed vertices is a Wu class for $\Gamma_{\mathrm{fix}}$, so we can write \[
\Omega = \omega+\sum_{j=1}^n\epsilon_j(v_j^++v_j^-), \quad \epsilon_j\in \{0,1\}. \]
Put $\boldsymbol{\epsilon}:=(\epsilon_1,\dots,\epsilon_n)^t$. Then we have \[
V\boldsymbol{\epsilon} \equiv (b_1+\omega_0,b_2,\dots,b_n)^t \pmod{2}.
\]
Since the support of $\Omega$ consists of pairwise nonadjacent vertices, we have $\omega_0\epsilon_1=0$. It follows that \[
\Omega^2 = \omega^2+\bigg(\sum_{j=1}^n\epsilon_j(v_j^++v_j^-) \bigg) ^2= -\boldsymbol{\omega}^tA\boldsymbol{\omega} - 2\boldsymbol{\epsilon}^tV\boldsymbol{\epsilon}.
\]

Next, recall the identity $\mathbf{r}^tV=q\mathbf{f}_1^t$ from the proof of Proposition \ref{prop:torus_Froyshov_formula}. This shows that the subgroup $V\mathbb{Z}^n\subset \mathbb{Z}^n$ is contained in the kernel of the map $\mathbf{r}^t:\mathbb{Z}^n\to \mathbb{Z}_q$, which is surjective since $r_n=1$. Since $\det V=q$, both $V\mathbb{Z}^n$ and $\ker \mathbf{r}^t$ have index $q$, and therefore $V\mathbb{Z}^n=\ker \mathbf{r}^t$. Since $V$ is invertible modulo $2$, we see that the set $\widetilde{\mathcal{D}}_0$ defined in the proof of Proposition \ref{prop:torus_Froyshov_formula} is exactly \[
\widetilde{\mathcal{D}}_0=\left\{V\mathbf{y}: \mathbf{y}\in\mathbb{Z}^n,\,\mathbf{y}\equiv\boldsymbol{\epsilon}\pmod{2} \right\}.
\]

For any $\mathbf{y}\in \mathbb{Z}^n$ with $\mathbf{y}\equiv\boldsymbol{\epsilon}\pmod{2}$, we have
\[
\mathbf{y}^tV\mathbf{y} = y_1^2+y_n^2 + \sum_{j=1}^{n-1}(y_j-y_{j+1})^2 + \sum_{j=1}^n(b_j-2)y_j^2,
\]
which clearly shows that \[
\min_{\mathbf{y}\equiv\boldsymbol{\epsilon}\pmod{2}}
\mathbf{y}^tV\mathbf{y} = \boldsymbol{\epsilon}^tV\boldsymbol{\epsilon}.
\]
The proof of Proposition \ref{prop:torus_Froyshov_formula} shows that the minimum of $\mathbf{z}^tV^{-1}\mathbf{z}$ over $\mathcal{D}_0$ agrees with the minimum over $\widetilde{\mathcal{D}}_0$. Hence \begin{align*}
    \min_{\mathbf{z}\in \mathcal{D}_0} \mathbf{z}^tV^{-1}\mathbf{z}=  \min_{\mathbf{z}\in \widetilde{\mathcal{D}}_0} \mathbf{z}^tV^{-1}\mathbf{z}=\min_{\mathbf{y}\equiv\boldsymbol{\epsilon}\pmod{2}} (V\mathbf{y})^t V^{-1}(V\mathbf{y})=\boldsymbol{\epsilon}^tV\boldsymbol{\epsilon}.
\end{align*}
Now Proposition \ref{prop:torus_Froyshov_formula} gives \begin{align*}
    \delta_R(Y_{p,q},\mathfrak{s}_0,\tau_{p,q})&=\frac{(1+m+2n)-\boldsymbol{\omega}^tA\boldsymbol{\omega} }{16}-\frac{1}{8}\min_{\mathbf{z}\in \mathcal{D}_0} \mathbf{z}^tV^{-1}\mathbf{z} \\
    & =\frac{|\Gamma_{p,q}|-\boldsymbol{\omega}^tA\boldsymbol{\omega}-2\boldsymbol{\epsilon}^tV\boldsymbol{\epsilon}}{16} \\
    & = \frac{-\sigma(W_{p,q})+\Omega^2}{16}.
\end{align*}
The last term is precisely $-\frac{1}{2}\bar{\mu}(Y_{p,q})$. 
\end{proof}

\subsection{The case \texorpdfstring{$p=2$}{p=2}}
\label{subsec:torus_p2}

In this subsection, we derive an explicit formula for $\delta_R(Y_{p,q},\mathfrak{s}_i,\tau_{p,q})$ when $p=2$. Write \[
q=4k+l \geq 3, \quad l\in\{1,3\}.
\]
By conjugation invariance \eqref{eq:conjugation_invariance}, it suffices to consider $0\leq i\leq (q-1)/2$.

\begin{proposition}\label{prop:T4q_Froyshov_formula}
For every $0\leq i\leq (q-1)/2$, we have \begin{equation}\label{eq:T4q_formula}
    \delta_R(Y_{2,q},\mathfrak{s}_i,\tau_{2,q})= \frac{(q-2i)^2}{8q}- \frac{3-l}{16}- \max\left\{ 0,\left\lceil\frac{k-i}{2}\right\rceil \right\}.
\end{equation}
\end{proposition}
\begin{proof}
For $p=2$, direct computation gives \[
\alpha=1, \quad \beta=\frac{lq-1}{4}, \quad c=\frac{l+1}{2}, \quad U=(2), \quad m =1, \quad A= \begin{pmatrix}
c & -1\\ -1 & 2
\end{pmatrix}, \quad \boldsymbol{\omega}=\begin{pmatrix}
0\\  (3-l)/2
\end{pmatrix}.
\]
In particular, $\omega_0=0$ and $\boldsymbol{\omega}^tA\boldsymbol{\omega}=3-l$. The negative continued fraction for $q/\beta$ is given by \[
\frac{q}{\beta}= \begin{cases}
[5,\underbrace{2,\dots,2}_{k-1}]^- & \text{if $l=1$}, \\
[2,2,3,\underbrace{2,\dots,2}_{k-1}]^- & \text{if $l=3$ and $k\geq 1$}, \\
[2,2]^-  & \text{if $l=3$ and $k=0$},
\end{cases}
\]
and therefore $n=k+l-1$. 

Since $\omega_0=0$, the set  $\widetilde{\mathcal{D}}_i$ defined in the proof of Proposition \ref{prop:torus_Froyshov_formula} is: \[
\widetilde{\mathcal{D}}_i := \left\{ \mathbf{z}=(z_1,\ldots,z_n)^t\in\mathbb{Z}^n:
\begin{array}{l}
z_j\equiv b_j \pmod{2} \quad (j=1,\dots,n),\\
\mathbf{r}^t\mathbf{z}\equiv i\pmod{q}
\end{array} \right\}.
\]
Also, Proposition \ref{prop:torus_Froyshov_formula} gives \begin{equation}\label{eq:p2_delta_minimum}
\delta_R(Y_{2,q},\mathfrak{s}_i,\tau_{2,q}) = \frac{l+k}{8} -  \frac{3-l}{16}-\frac{1}{8}\min_{\mathbf{z}\in \widetilde{\mathcal{D}}_i}\mathbf{z}^tV^{-1}\mathbf{z}.
\end{equation}
Comparing with \eqref{eq:T4q_formula}, it remains to prove that \begin{equation}\label{eq:reduction}
    \min_{\mathbf{z}\in \widetilde{\mathcal{D}}_i}\mathbf{z}^tV^{-1}\mathbf{z}= (l+k) +8 \max\left\{ 0,\left\lceil\frac{k-i}{2}\right\rceil \right\} -\frac{(q-2i)^2}{q}.
\end{equation}

\vspace{1mm}
\textbf{Step 1:} To compute the minimum in an efficient way, we first realize $V$ as a Gram matrix of $n=k+l-1$ vectors in $\mathbb{Z}^{n+1}=\mathbb{Z}^{k+l}$. Let $\mathbf{f}_0,\dots,\mathbf{f}_{n}$ denote the standard basis of $\mathbb{Z}^{n+1}$, and put \[
\nu = (\underbrace{1,\dots,1}_{l},\underbrace{2,\dots,2}_{k})= \sum_{j=0}^{l-1}\mathbf{f}_j + 2\sum_{j=l}^{n}\mathbf{f}_j.
\]
Then $\|\nu\|^2=l+4k=q$. 

If $l=1$, define \[
\xi_1=-2\mathbf{f}_0+\mathbf{f}_1, \quad \xi_j=-\mathbf{f}_{j-1}+\mathbf{f}_j \quad(j=2,\dots,k).
\]
If $l=3$, define \[
\xi_1=\mathbf{f}_0-\mathbf{f}_1, \quad \xi_2=\mathbf{f}_1-\mathbf{f}_2, \quad \xi_3=-\mathbf{f}_0-\mathbf{f}_1+\mathbf{f}_3, \quad \xi_j=-\mathbf{f}_{j-1}+\mathbf{f}_j \quad(j=4,\dots,k+2).
\]
In either case, the vectors $\xi_1,\dots,\xi_{k+l-1}$ have Gram matrix $V$, and the orthogonal complement is generated by the vector $\nu$. 

\vspace{1mm}
\textbf{Step 2:} For every $\mathbf{z}\in\mathbb{Z}^n$ satisfying $z_j\equiv b_j\pmod{2}$, there exists $\mathbf{x}\in \mathbb{Z}^{n+1}$ such that each $x_i$ is odd and \begin{equation}\label{eq:Step_2}
    z_j=\langle\mathbf{x},\xi_j\rangle \quad(j=1,\dots,n).
\end{equation}
To see this, let $x_0$ be any odd integer. If $l=1$, define \[
x_1=z_1+2x_0, \quad x_j=z_j+x_{j-1} \quad (j=2,\dots,n).
\]
If $l=3$, define \[
x_1=x_0-z_1, \quad x_2=x_1-z_2, \quad x_3=z_3+x_0+x_1, \quad x_j=z_j+x_{j-1} \quad (j=4,\dots,n). 
\]
In either case, the vector $\mathbf{x}=(x_0,\dots,x_n)^t$ satisfies \eqref{eq:Step_2}. Observe that two such vectors $\mathbf{x}$ representing the same $\mathbf{z}$ differ by an element of $2\mathbb{Z}\nu$.

\vspace{1mm}
\textbf{Step 3:} For $\mathbf{x}$ as in \textbf{Step 2}, we have: \begin{equation}\label{eq:Step_3}
\mathbf{z}^tV^{-1}\mathbf{z}=\|\mathbf{x}\|^2-\frac{\langle\mathbf{x},\nu\rangle^2}{q}.
\end{equation}
To see this, view $\mathbf{x}$ as an element of $\mathbb{R}^{n+1}$, and consider the orthogonal decomposition $\mathbb{R}^{n+1}=\langle \xi_1,\dots,\xi_n\rangle \oplus \langle \nu\rangle$. Then $\mathbf{x}$ decomposes as $\mathbf{x}=(\mathbf{x}-\frac{\langle\mathbf{x},\nu\rangle}{q}\nu)+\frac{\langle\mathbf{x},\nu\rangle}{q}\nu$. Also, since $\langle \mathbf{x}-\frac{\langle\mathbf{x},\nu\rangle}{q}\nu, \xi_j\rangle = \langle \mathbf{x},\xi_j\rangle = z_j$, the coefficient vector of $\mathbf{x}-\frac{\langle\mathbf{x},\nu\rangle}{q}\nu$ with respect to the basis $\{\xi_1,\dots,\xi_n\}$ is precisely $V^{-1}\mathbf{z}$, and therefore \[
\left\|\mathbf{x}-\frac{\langle\mathbf{x},\nu\rangle}{q}\nu\right\|^2=  (V^{-1}\mathbf{z})^t V (V^{-1}\mathbf{z})=\mathbf{z}^tV^{-1}\mathbf{z},
\]
which gives \eqref{eq:Step_3}.

\vspace{1mm}
\textbf{Step 4:} We next consider the condition $\mathbf{r}^t\mathbf{z}\equiv i \pmod{q}$. A direct computation gives: \[
\langle\mathbf{x},\nu\rangle= \begin{cases}
qx_0+2\mathbf{r}^t\mathbf{z} &\text{if $l=1$},\\
2\mathbf{r}^t\mathbf{z}-q(x_0-x_1-x_2) & \text{if $l=3$}.
\end{cases}
\]
For example, if $l=1$, then $\mathbf{r}=(k,k-1,\dots,2,1)^t$, so \begin{align*}
    qx_0+2\mathbf{r}^t\mathbf{z} & =qx_0+2 \sum_{j=1}^k (k+1-j) z_j = qx_0+2\sum_{j=1}^k  (k+1-j)\langle \mathbf{x},\xi_j\rangle \\
    & = qx_0+2k(-2x_0+x_1)+2\sum_{j=2}^k (k+1-j)(-x_{j-1}+x_j) \\
    & = x_0+2(x_1+\cdots +x_k)= \langle \mathbf{x},\nu\rangle. 
\end{align*} 
The case $l=3$ is similar. Since $x_0$ and $x_0-x_1-x_2$ are odd, both cases give \[
 \langle \mathbf{x},\nu\rangle \equiv q+2\mathbf{r}^t\mathbf{z} \pmod{2q},
\]
which shows that $\mathbf{r}^t\mathbf{z} \equiv i\pmod{q}$ if and only if $\langle \mathbf{x},\nu\rangle \equiv q+2i \pmod{2q}$. 

Adding a multiple of $2\nu$ to $\mathbf{x}$ preserves $\mathbf{z}$ and changes $\langle \mathbf{x},\nu\rangle$ by a multiple of $2q$. Thus, for every $\mathbf{z}\in\widetilde{\mathcal{D}}_i$, there is a unique $\mathbf{x}$ as in \textbf{Step 2} such that $\langle \mathbf{x},\nu\rangle=2i-q$. Using \eqref{eq:Step_3}, with $\mathbf{x}$ replaced by $-\mathbf{x},$ we obtain \[
\min_{\mathbf{z}\in \widetilde{\mathcal{D}}_i}\mathbf{z}^tV^{-1}\mathbf{z}= \min_{\mathbf{x}} \|\mathbf{x}\|^2 - \frac{(q-2i)^2}{q},
\]
where the minimum is taken over all $\mathbf{x}\in \mathbb{Z}^{n+1}$ such that each $x_i$ is odd and $\langle \mathbf{x},\nu\rangle=q-2i$. 

\vspace{1mm}
\textbf{Step 5:} Comparing with \eqref{eq:reduction}, it remains to prove that \begin{equation}\label{eq:Step_5}
\min_{\mathbf{x}} \|\mathbf{x}\|^2 = (l+k) +8 \max\left\{ 0,\left\lceil\frac{k-i}{2}\right\rceil \right\} .
\end{equation}
Since every $x_i$ is odd, we have \begin{equation}\label{eq:Step_5-1}
    \|\mathbf{x}\|^2\geq n+1=l+k, \quad \|\mathbf{x}\|^2\equiv l+k\pmod{8}.
\end{equation}
Note that every odd integer $x$ satisfies \[
x^2\geq 2x-1, \quad x^2\geq 4x-3.
\]
Applying these inequalities to the first $l$ and the last $k$ coordinates of $\mathbf{x}$, respectively, we get \begin{equation}\label{eq:Step_5-2}
    \|\mathbf{x}\|^2 \geq \sum_{j=0}^{l-1}(2x_j-1)+\sum_{j=l}^{k+l-1} (4x_j-3) = 2\langle\mathbf{x},\nu\rangle-l-3k = (l+k)+4(k-i).
\end{equation}
Combining \eqref{eq:Step_5-1} and \eqref{eq:Step_5-2} gives \[
\|\mathbf{x}\|^2\geq  (l+k) +8 \max\left\{ 0,\left\lceil\frac{k-i}{2}\right\rceil \right\}.
\]

\vspace{1mm}
\textbf{Step 6:} Finally, let us show that this bound can be attained. If $i<k$, let $\mathbf{x}$ be the vector obtained from $(1,\dots,1)$ as follows: \begin{itemize}[leftmargin=2em]
    \item Replace $\lfloor (k-i)/2\rfloor$ of the last $k$ coordinates by $3$.
    \item If $k-i$ is odd, also replace one of the first $l$ coordinates by $3$. 
\end{itemize}
Then $\mathbf{x}$ satisfies \[
\langle\mathbf{x},\nu\rangle=q-2i, \quad \|\mathbf{x}\|^2 =
(l+k)+8\left\lceil\frac{k-i}{2}\right\rceil = (l+k) +8 \max\left\{ 0,\left\lceil\frac{k-i}{2}\right\rceil \right\}.
\]

If $i\geq k$, let $\mathbf{x}$ be the vector obtained from $(1,\dots,1)$ as follows: \begin{itemize}[leftmargin=2em]
    \item Replace $\lfloor (i-k)/2\rfloor$ of the last $k$ coordinates by $-1$. (This is possible since $0\leq i\leq (q-1)/2$.)
    \item If $i-k$ is odd, also replace one of the first $l$ coordinates by $-1$. 
\end{itemize}
Then $\mathbf{x}$ satisfies \[
\langle\mathbf{x},\nu\rangle=q-2i, \quad \|\mathbf{x}\|^2 =l+k = (l+k) +8 \max\left\{ 0,\left\lceil\frac{k-i}{2}\right\rceil \right\}.
\]
This completes the proof. 
\end{proof}

\section{The nonorientable 4-genus}\label{sec:nonorientable_genus}

We apply the computations in Section \ref{sec:even_torus_knots} to obstruct certain torus knots from bounding smooth M\"obius bands in the $4$-ball $B^4$. For a knot $K\subset S^3$, its (smooth) \emph{nonorientable $4$-genus} is defined by \[
\gamma_4(K):=\min_F b_1(F),
\]
where the minimum is taken over all smoothly embedded compact nonorientable surfaces $F\subset B^4$  with $\partial F=K$. Thus, $\gamma_4(K)=1$ if and only if $K$ bounds a smooth M\"obius band in $B^4$.

\subsection{Double branched covers}
\label{subsec:branched_covers}

Let $K\subset S^3$ be a knot, and denote by $Y_K:=\Sigma_2(S^3,K)$ its double branched cover. Since $|H^2(Y_K;\mathbb{Z})|=\det K$ is odd, $Y_K$ is a $\mathbb{Z}_2$-homology sphere. Let $\tau_K:Y_K\to Y_K$ denote the deck transformation.

Suppose $K$ bounds a smoothly embedded M\"obius band $F\subset B^4$, and denote by $W:=\Sigma_2(B^4,F)$ the double cover of $B^4$ branched along $F$. Let $\tau_W:W\to W$ denote the deck transformation. Then $\partial W=Y_K$ and $\tau_W|_{Y_K}=\tau_K$.

\begin{lemma}\label{lem:mobius_cover}
We have \[
H_i(W;\mathbb{Z}_2)\cong \begin{cases}
\mathbb {Z}_2 &\text{if $i=0,2$}, \\
0 &\text{otherwise}.
\end{cases}
\]
Moreover, $\tau_W^*=-\operatorname{id}$ on $H^2(W;\mathbb{Z})$ and every $\mathfrak{t}\in \operatorname{Spin}^c(W)$ satisfies $\tau_W^*\mathfrak{t}\cong \overline{\mathfrak{t}}$. 
\end{lemma}
\begin{proof}
For a double branched cover, we have the following exact sequence; see \cite[Theorem 1]{Lee-Weintraub-1995}: \[
\cdots \longrightarrow H_{i+1}(B^4;\mathbb{Z}_2) \longrightarrow H_i(B^4,F;\mathbb{Z}_2) \longrightarrow H_i(W;\mathbb{Z}_2)\longrightarrow H_i(B^4;\mathbb{Z}_2)\longrightarrow \cdots.
\]
This proves the first statement. In particular, $H^2(W;\mathbb{Z})$ has no $2$-torsion. 

Next, consider the ring $\mathbb{Z}[1/2]$. Since $2$ is invertible in $\mathbb{Z}[1/2]$, the transfer homomorphism gives \[
H^2(W;\mathbb{Z}[1/2])^{\tau_W^*} \cong H^2(W/\tau_W;\mathbb{Z}[1/2]) =H^2(B^4;\mathbb{Z}[1/2])=0.
\]
It follows that $\tau_W^*=-\operatorname{id}$ on $H^2(W;\mathbb{Z}[1/2])$. Note that the homomorphism $H^2(W;\mathbb{Z})\to H^2(W;\mathbb{Z}[1/2])$ is injective since $H^2(W;\mathbb{Z})$ has no $2$-torsion. Thus, $\tau_W^*=-\operatorname{id}$ on $H^2(W;\mathbb{Z})$. Therefore, for any $\mathfrak t\in\operatorname{Spin}^c(W)$, we have \[
c_1(\tau_W^*\mathfrak{t})=-c_1(\mathfrak{t}) =c_1(\overline{\mathfrak{t}}).
\]
Since $H^2(W;\mathbb{Z})$ has no $2$-torsion, the map $c_1:\operatorname{Spin}^c(W)\to H^2(W;\mathbb{Z})$ is injective. Hence $\tau_W^*\mathfrak{t}\cong \overline{\mathfrak{t}}$. 
\end{proof}

The same argument, using $Y_K/\tau_K=S^3$, shows that $\tau_K^*=-\operatorname{id}$ on $H^2(Y_K;\mathbb{Z})$. Consequently, every $\operatorname{spin}^c$ structure on $Y_K$ is Real. Write $\mathfrak{s}_0$ for the unique spin structure on $Y_K$ and put \[
\mathfrak{s}_x:=\mathfrak{s}_0+x\in \operatorname{Spin}^c(Y_K) \quad (x\in H^2(Y_K;\mathbb{Z})).
\]

Since $b_2(W)=1$ by Lemma \ref{lem:mobius_cover}, the intersection form $Q_W$ is represented by a $1\times 1$ matrix \[
Q_W \cong \langle \epsilon d\rangle \quad (\epsilon \in \{\pm 1\}, \, d\in \mathbb{Z}_{>0}).
\]
Since $\partial W=Y_K$ is a $\mathbb{Z}_2$-homology sphere, the integer $d$ must be odd. Consider the cohomology long exact sequence of the pair $(W,Y_K)$: \[
\begin{tikzcd}[matrix scale=0.6]
    0 \ar[r] & H^2(W,Y_K;\mathbb{Z}) \ar[d, equal] \ar[r, "j"] & H^2(W;\mathbb{Z}) \ar[d, equal] \ar[r, "r"] & H^2(Y_K;\mathbb{Z}) \ar[r] & H^3(W,Y_K;\mathbb{Z}) \ar[r] \ar[d, equal] & H^3(W;\mathbb{Z}) \ar[r] \ar[d, equal] & 0 \\
     & \mathbb{Z}\oplus T_2 & \mathbb{Z}\oplus T_1 & & T_1 & T_2,
\end{tikzcd}
\]
where $T_1=\operatorname{Tor} H^2(W;\mathbb{Z})$ and $T_2=\operatorname{Tor}H^2(W,Y_K;\mathbb{Z})$. Note that the restriction $j:T_2\to T_1$ is injective. In particular, $|T_2|$ divides $|T_1|$. Moreover, the proof of \cite[Lemma 2.1]{Owens-Strle-2006} shows that: \begin{itemize}[leftmargin=2em]
    \item The image of $r|_{T_1}:T_1\to H^2(Y_K;\mathbb{Z})$ has order $|T_1|/|T_2|$.    
    \item The image of $r:H^2(W;\mathbb{Z})\to H^2(Y_K;\mathbb{Z})$ has order $d|T_1|/|T_2|$.
    \item $|H^2(Y_K;\mathbb{Z})|=d(|T_1|/|T_2|)^2$. 
\end{itemize}

Next, fix $\xi\in H^2(W;\mathbb{Z})$ whose class generates $H^2(W;\mathbb{Z})/{\operatorname{Tor}}\cong \mathbb{Z}$. For every $u,v\in H^2(W;\mathbb{Z})$, we can uniquely write $u=a\xi+u'$ and $v=b\xi +v'$ where $a,b\in \mathbb{Z}$ and $u',v'\in H^2(W;\mathbb{Z})$ are torsion. Moreover, the linking form \[
\lambda_{Y_K}:H^2(Y_K;\mathbb{Z})\times H^2(Y_K;\mathbb{Z})\longrightarrow \mathbb{Q}/\mathbb{Z}
\]
satisfies \begin{equation}\label{eq:linking_form}
    \lambda_{Y_K}(r(u),r(v)) \equiv -\frac{\epsilon ab}{d} \pmod{\mathbb{Z}};
\end{equation}
see \cite[pp.199-200]{Owens-Strle-2006}. 

\begin{proposition}\label{prop:mobius_Real_Froyshov}
Suppose $\mathfrak{s}\in \operatorname{Spin}^c(Y_K)$ satisfies \[
c_1(\mathfrak{s})=r(a\xi+\eta)
\]
for some odd integer $a$ and $\eta\in \operatorname{Tor}H^2(W;\mathbb{Z})$. Then \begin{equation}\label{eq:mobius_Real_Froyshov}
\epsilon\,\delta_R(Y_K,\mathfrak s,\tau_K)
\leq\frac{a^2/d-1}{16}.
\end{equation}
\end{proposition}

\begin{proof}
Since $H^2(W;\mathbb{Z})$ has no $2$-torsion, the map $c_1:\operatorname{Spin}^c(W)\to H^2(W;\mathbb{Z})$ is a bijection onto the set of characteristic classes. Note that characteristic classes in $H^2(W;\mathbb{Z})$ are of the form $a\xi+\eta$ with $a$ odd and $\eta\in \operatorname{Tor}H^2(W;\mathbb{Z})$. On the other hand, since $|H^2(Y_K;\mathbb{Z})|$ is odd, the map $c_1:\operatorname{Spin}^c(Y_K)\to H^2(Y_K;\mathbb{Z})$ is injective. Therefore, if $\mathfrak{s}\in \operatorname{Spin}^c(Y_K)$ satisfies $c_1(\mathfrak{s})=r(a\xi+\eta)$ for some odd $a\in \mathbb{Z}$ and $\eta\in \operatorname{Tor}H^2(W;\mathbb{Z})$, then $\mathfrak{s}=\mathfrak{t}|_{Y_K}$ where $c_1(\mathfrak{t})=a\xi+\eta$. Thus, $(W,\mathfrak{t},\tau_W)$ is a Real $\operatorname{spin}^c$ filling of $(Y_K,\mathfrak{s},\tau_K)$.  

Now, if $\epsilon=-1$, then Corollary \ref{cor:real_froyshov_filling} gives \[
\frac{1-a^2/d}{16} \leq\delta_R(Y_K,\mathfrak{s},\tau_K).
\]
If $\epsilon=1$, then applying Corollary \ref{cor:real_froyshov_filling} to the filling $(-W,\mathfrak{t},\tau_W)$ of $(-Y_K,\mathfrak{s},\tau_K)$ gives \[
\frac{1-a^2/d}{16} \leq\delta_R(-Y_K,\mathfrak{s},\tau_K)=-\delta_R(Y_K,\mathfrak{s},\tau_K).
\]
This proves \eqref{eq:mobius_Real_Froyshov} in both cases.
\end{proof}

\subsection{Applications}\label{subsec:T4q_application}

We now consider the torus knots $T_{4,q}$, where $q>1$ is odd. Since it is already known that $\gamma_4(T_{4,q})=2$ if $q\equiv 5,7\pmod{8}$ \cite{BKST-2025}, we mostly focus on the case $q\equiv 1,3\pmod{8}$. Recall from Lemma \ref{lem:torus_labels} that $H^2(Y_{2,q};\mathbb{Z})$ is cyclic of order $q$, generated by $h=\eta_{v_n^+}|_{Y_{2,q}}$. Recall also the labeling $\mathfrak{s}_i=\mathfrak{s}_0+ih$ from \eqref{eq:torus_labels}. Consider the linking form $\lambda_{Y_{2,q}}$ of $Y_{2,q}$.

\begin{lemma}\label{lem:torus_linking_form}
    $\lambda_{Y_{2,q}}(h,h)\equiv -\dfrac{2}{q}\pmod{\mathbb{Z}}$.
\end{lemma}

It is proved in \cite[Proposition 5.4]{BKST-2025} that $\lambda_{Y_{2,q}}(x,x)\equiv -\frac{2}{q}\pmod{\mathbb{Z}}$ for \emph{some} generator $x\in H^2(Y_{2,q};\mathbb{Z})$. Our lemma proves this for the explicitly chosen generator $h$. (It can happen that $\lambda_{Y_{2,q}}(h',h')\not\equiv -\frac{2}{q}\pmod{\mathbb{Z}}$ if we choose another generator $h'$.) 
\begin{proof}
We use the notation in Section \ref{subsec:plumbing_graphs}. Since $p=2$, we have $\alpha=1$ and $m=1$. The intersection matrix \eqref{eq:intersection_matrix} for the plumbing $W_{2,q}$ with respect to the basis $\{v_0,u_1,v_1^+,\dots,v_n^+,v_1^-,\dots,v_n^-\}$ for $H_2(W_{2,q};\mathbb{Z})$ is given by \[
  Q=-\begin{pmatrix}
    c & -1 & -\mathbf{e}_n^t & -\mathbf{e}_n^t \\
    -1 & 2 & 0 & 0 \\
    -\mathbf{e}_n & 0 & V & 0 \\
    -\mathbf{e}_n & 0 & 0 & V
\end{pmatrix}.
\]
Since $h=\eta_{v_n^+}|_{Y_{2,q}}$, we have \[
\lambda_{Y_{2,q}}(h,h)\equiv -(Q^{-1})_{v_n^+,v_n^+} \pmod{\mathbb{Z}},
\]
so we have to compute $-(Q^{-1})_{v_n^+,v_n^+}=(-Q)^{-1}_{v_n^+,v_n^+}$. As in the proof of Lemma \ref{lem:torus_labels}, write $\mathbf{z}:=(1,\mathbf{e}_n^t,\mathbf{e}_n^t)^t$ and $D:=\operatorname{diag}(2,V,V)$, so that \[
-Q=\begin{pmatrix}
    c & -\mathbf{z}^t \\ -\mathbf{z} & D 
\end{pmatrix}.
\]
Thus, the usual formula for the inverse of a block matrix gives \[
(-Q)^{-1}=\begin{pmatrix}
     * & * \\
     * & D^{-1}+D^{-1}\mathbf{z}(c-\mathbf{z}^tD^{-1}\mathbf{z})^{-1}\mathbf{z}^tD^{-1}
\end{pmatrix}.
\]
Clearly, $D^{-1}=\operatorname{diag}(1/2, V^{-1},V^{-1})$. First, note that \[
(V^{-1})_{1n}=\frac{(-1)^{1+n}\det V_{n,1}}{\det V}=\frac{1}{q},
\]
where $V_{n,1}$ is obtained from $V$ by removing the $n$-th row and the first column, which is a lower triangular matrix with $-1$'s on the diagonal. Similarly, we have \[
(V^{-1})_{11}=\frac{(-1)^{1+1}\det V_{1,1}}{\det V}=\frac{\beta}{q}, \quad (V^{-1})_{nn}=\frac{(-1)^{n+n} \det V_{n,n}}{\det V}=\frac{\beta'}{q},
\]
where $\beta'$ is the unique integer satisfying $0<\beta'<q$ and $\beta\beta'\equiv 1\pmod{q}$; this follows from the well-known property of negative continued fractions that $q/\beta=[b_1,\dots,b_n]^{-}$ implies $q/\beta'=[b_n,\dots,b_1]^{-}$. It follows that \begin{align*}
    c-\mathbf{z}^tD^{-1}\mathbf{z}& = c-\left(\frac{1}{2}+(V^{-1})_{11}+(V^{-1})_{11} \right)=c-\frac{1}{2}-\frac{2\beta}{q}=\frac{1}{2q},
\end{align*}
where the last equality uses \eqref{eq:torus_seifert}. Next, note that \[
D^{-1}\mathbf{z}=\left(1/2, (V^{-1})_{11},\dots,(V^{-1})_{n1},(V^{-1})_{11},\dots,(V^{-1})_{n1}\right)^t.
\]
We conclude that the $(v_n^+,v_n^+)$-entry of $-Q^{-1}$ is given by \[
(V^{-1})_{nn}+2q\left((V^{-1})_{n1}\right)^2=\frac{\beta'}{q}+2q\cdot \frac{1}{q^2}=\frac{\beta'+2}{q}.
\]
Recall that $4\beta\equiv -1\pmod{q}$ by the definition of $\beta$. Since $\beta\beta'\equiv 1\pmod{q}$, it follows that $\beta'\equiv -4\pmod{q}$. This completes the proof. 
\end{proof}

Suppose $T_{4,q}$ bounds a smooth M\"obius band $F\subset B^4$, and the intersection form of its double branched cover $W:=\Sigma_2(B^4,F)$ is isomorphic to $\langle\epsilon d\rangle$, where $\epsilon\in\{\pm1\}$ and $d>0$. Consider the homomorphism \[
r:H^2(W;\mathbb{Z})\longrightarrow H^2(Y_{2,q};\mathbb{Z}),
\]
and put $s:=|r(\operatorname{Tor}H^2(W;\mathbb{Z}))|$. We have seen in Section \ref{subsec:branched_covers} that \[
q=|H^2(Y_{2,q};\mathbb{Z})|=ds^2, \quad |\operatorname{im}r|=ds.
\]
Since $H^2(Y_{2,q};\mathbb{Z})$ is a cyclic group generated by $h$, it follows that $\operatorname{im}r$ and $r(\operatorname{Tor}H^2(W;\mathbb{Z}))$ are generated by $sh$ and $dsh$, respectively. 

As in Section \ref{subsec:branched_covers}, choose $\xi\in H^2(W;\mathbb{Z})$ whose class generates $H^2(W;\mathbb{Z})/{\operatorname{Tor}}\cong \mathbb{Z}$. Then \eqref{eq:linking_form} gives \begin{equation}\label{eq:generator_xi}
    \lambda_{Y_{2,q}}(r(\xi),r(\xi))\equiv -\frac{\epsilon}{d}\pmod{\mathbb{Z}}.
\end{equation}
On the other hand, if we write $r(\xi)=ush$ for some integer $u\in \mathbb{Z}$, then Lemma \ref{lem:torus_linking_form} gives \[
\lambda_{Y_{2,q}}(r(\xi),r(\xi))=\lambda_{Y_{2,q}}(ush,ush)\equiv -\frac{2(us)^2}{q}=-\frac{2u^2}{d}\pmod{\mathbb{Z}}.
\]
It follows that $2u^2\equiv \epsilon\pmod{d}$. 

Next, let $a\in \mathbb{Z}$ be an odd integer and let $i\in \mathbb{Z}_q$. Suppose $2i\equiv aus\pmod{ds}$, so that $2i=aus+bds$ for some integer $b$. Then for the $\operatorname{spin}^c$ structure $\mathfrak{s}_i\in \operatorname{Spin}^c(Y_{2,q})$, we have \[
c_1(\mathfrak{s}_i)=2ih =(aus+bds)h=r(a\xi)+bdsh \in r(a\xi)+r(\operatorname{Tor}H^2(W;\mathbb{Z})),
\]
and therefore Proposition \ref{prop:mobius_Real_Froyshov} implies \[
\epsilon \delta_R(Y_{2,q},\mathfrak{s}_i,\tau_{2,q}) \leq \frac{a^2/d-1}{16}.
\]
In summary, we have proved the following lemma.

\begin{lemma}\label{lem:T4q_filling_constraints}
Suppose $T_{4,q}$ bounds a smooth M\"obius band $F\subset B^4$, and the intersection form of its double branched cover $W:=\Sigma_2(B^4,F)$ is isomorphic to $\langle\epsilon d\rangle$, where $\epsilon\in\{\pm1\}$ and $d>0$. Then \begin{equation}\label{eq:T4q_filling_congruence}
q=ds^2,\quad 2u^2\equiv\epsilon\pmod{d}
\end{equation}
for some odd integer $s>0$ and some integer $u\in \mathbb{Z}$. Moreover, \begin{equation}\label{eq:T4q_rank_one_inequality}
\epsilon \delta_R(Y_{2,q},\mathfrak{s}_i,\tau_{2,q}) \leq \frac{a^2/d-1}{16}
\end{equation}
whenever $a\in \mathbb{Z}$ is odd and $2i\equiv aus \pmod{ds}$. 
\end{lemma}

The linking form condition \eqref{eq:T4q_filling_congruence} has the following interpretation.
\begin{corollary}\label{cor:torus_linking_form}
    Let $q\geq 1$ be an odd integer. \begin{enumerate}[label=\textup{(\alph*)}, leftmargin=2em]
        \item  If $q\equiv 1\pmod{8}$ and $\gamma_4(T_{4,q})=1$, then $q=X^2\pm 8Y^2$ for some $X,Y\in \mathbb{Z}$.
    \item If $q\equiv 3\pmod{8}$ and $\gamma_4(T_{4,q})=1$, then $q=X^2+2Y^2$ for some $X,Y\in \mathbb{Z}$.
    \end{enumerate}
\end{corollary}
\begin{proof} The case $q=1$ is immediate. Suppose $q>1$ and $\gamma_4(T_{4,q})=1$, and let $\epsilon,d,s,u$ be as in Lemma \ref{lem:T4q_filling_constraints}. Consider the quadratic ring $R_\epsilon:=\mathbb{Z}[\sqrt{2\epsilon}]$, and define \[
I_\epsilon:=\left\{x+y\sqrt{2\epsilon}\in R_\epsilon:y\equiv ux\pmod{d}\right\}.
\]
It is well-known that $R_\epsilon$ is a Euclidean domain with the norm function $x+y\sqrt{2\epsilon}\mapsto |x^2-2\epsilon y^2|$. Since $2u^2\equiv\epsilon\pmod{d}$, $I_\epsilon$ is an ideal of $R_\epsilon$. Since every Euclidean domain is a PID, the ideal $I_\epsilon$ is principal, say $I_\epsilon=(x+y\sqrt{2\epsilon})$. As an abelian group, $I_\epsilon$ is generated by $1+u\sqrt{2\epsilon}$ and $d\sqrt{2\epsilon}$, so it is of index $d$. Therefore, any generator of $I_\epsilon$ (as an ideal) must have norm $d$.

Thus, if $\epsilon=-1$, then $x^2+2y^2=d$. If $\epsilon=1$, then $|x^2-2y^2|=d$, and by multiplying by the unit $1+\sqrt{2}$ if necessary, we may assume that $x^2-2y^2=d$. Since $q=ds^2$, in either case we have \[
(sx)^2-2\epsilon(sy)^2= s^2(x^2-2\epsilon y^2)=q.
\]

Now, if $q\equiv 1 \pmod{8}$, then $d\equiv 1\pmod{8}$. Since $d=x^2-2\epsilon y^2$ and $d$ is odd, $x$ is also odd, so $x^2\equiv 1\pmod{8}$, and therefore $2y^2\equiv 0 \pmod{8}$, i.e., $y$ is even. Writing $y=2z$ gives \[
q=(sx)^2-2\epsilon(2sz)^2 = (sx)^2 \pm 8(sz)^2.
\]
If $q\equiv 3\pmod{8}$, then $d\equiv3\pmod{8}$. Since $x^2-2y^2\equiv 1$ or $7\pmod{8}$ for odd $x$, we must have $\epsilon =-1$, and \[
q=(sx)^2+2(sy)^2. \]
This completes the proof. 
\end{proof}

\begin{remark}\label{rmk:BKST}
The linking form condition \eqref{eq:T4q_filling_congruence} also follows from the proof of \cite[Lemma 5.10]{BKST-2025}, so Corollary \ref{cor:torus_linking_form} is essentially a consequence of their results. While they do not state this corollary explicitly, they point out that if $q\equiv 1, 3 \pmod{8}$ and $\gamma_4(T_{4,q})=1$, then any prime factor of $q$ congruent to $5$ modulo $8$ must have even exponent in the prime factorization of $q$ \cite[Remark 1.14 and Example 5.11]{BKST-2025}. Observe that this property can be easily recovered from Corollary \ref{cor:torus_linking_form} using elementary number theory.
\end{remark}

Now we strengthen the result of Corollary \ref{cor:torus_linking_form} in the case $q\equiv 1\pmod{8}$ using the Real Fr{\o}yshov inequality \eqref{eq:T4q_rank_one_inequality}. Put $k:=(q-1)/4$, which is even, and write \[
\rho(i)=\begin{cases}
    0 & \text{if $i$ is even},\\
    1 & \text{if $i$ is odd}.
\end{cases}
\]
By Proposition \ref{prop:T4q_Froyshov_formula}, we have \begin{equation}\label{eq:T4q_mod8_piecewise}
\delta_R(Y_{2,q},\mathfrak{s}_i,\tau_{2,q})=\begin{cases}
\dfrac{i^2}{2q}-\dfrac{\rho(i)}2 & \text{if $0\leq i\leq k$}, \\
\dfrac{(q-2i)^2-q}{8q} & \text{if $k\leq i\leq(q-1)/2$}.
\end{cases}
\end{equation}

\begin{proof}[Proof of Theorem \textup{\ref{thm:T4q_nonorientable_genus}}]
The case $q=1$ is immediate, so assume $q>1$. Suppose $\gamma_4(T_{4,q})=1$, and let $\epsilon,d,s,u$ be as in Lemma \ref{lem:T4q_filling_constraints}. If $d=1$, then $q=s^2$, so we are done. Thus, assume $d>1$. 

\textbf{Case 1:} $\epsilon=-1$. The proof of Corollary \ref{cor:torus_linking_form} shows that there exist integers $x,y\in \mathbb{Z}$ such that $d=x^2+2y^2$ and $y\equiv ux \pmod{d}$, and we may assume $y\geq 0$ by replacing $(x,y)$ with $(-x,-y)$ if necessary. Since $q\equiv 1\pmod{8}$ and $q=ds^2$, we have $d\equiv1\pmod{8}$, so $x$ is odd and $y$ is even. Writing $y=2z$ gives $d=x^2+8z^2$ and $2z\equiv ux\pmod{d}$. 

Set $a=x$ and $i=sz$. Then $2i\equiv aus \pmod{ds}$. Moreover, since \[
4k+1= q=ds^2=(x^2+8z^2)s^2=(sx)^2+8i^2 \geq 1+8i^2,
\]
we have $i\leq k$, and therefore \eqref{eq:T4q_mod8_piecewise} gives \[
\delta_R(Y_{2,q},\mathfrak{s}_i,\tau_{2,q})=\frac{i^2}{2q}-\frac{\rho(i)}{2} =\frac{z^2}{2d}-\frac{\rho(z)}2.
\]
On the other hand, \eqref{eq:T4q_rank_one_inequality} gives \[
\delta_R(Y_{2,q},\mathfrak{s}_i,\tau_{2,q}) \geq \frac{1-x^2/d}{16}=\frac{z^2}{2d}.
\]
Hence $\rho(z)=0$, i.e., $z$ is even. Writing $z=2w$ gives \[
q=(sx)^2+32(sw)^2.
\]

\textbf{Case 2:} $\epsilon=1$. We show that $d>1$ is impossible. Since $d\equiv 1\pmod{8}$, we have $d\geq 9$. Since $ds$ is odd, the congruence $2i\equiv us \pmod{ds}$ has a unique solution modulo $ds$, say $i_0\in \mathbb{Z}$. Thus, the congruence $2i\equiv us \pmod{ds}$ has exactly $s$ solutions modulo $q$, which are $i_0, i_0+ds,\dots, i_0+(s-1)ds$. Moreover, \eqref{eq:T4q_rank_one_inequality} with $a=1$ gives \[
\delta_R(Y_{2,q},\mathfrak{s}_{i_0+jds},\tau_{2,q}) \leq \frac{1/d-1}{16}<0 \quad (j=0,1,\dots,s-1).
\]

On the other hand, for any $i\in \mathbb{Z}_q$, \eqref{eq:T4q_mod8_piecewise} (and \eqref{eq:conjugation_invariance}) shows that if $\delta_R(Y_{2,q},\mathfrak{s}_i,\tau_{2,q})<0$, then $i$ is contained in the union of arcs \[
\left(-\sqrt{q},\sqrt{q}\right)\cup \left(\frac{q-\sqrt q}{2},\frac{q+\sqrt q}{2}\right) \subset \mathbb{R}/q\mathbb{Z},
\]
where we view $i\in \mathbb{Z}_q=\mathbb{Z}/q\mathbb{Z}\subset \mathbb{R}/q\mathbb{Z}$. These arcs have lengths $2\sqrt{q}$ and $\sqrt{q}$, respectively. Since $ds=\sqrt{dq}>2\sqrt{q}$, it follows that each arc can contain at most one among $i_0,i_0+ds,\dots,i_0+(s-1)ds$. Therefore, $s\leq 2$. Since $s$ is odd, we get $s=1$ and $q=d$.

Now let $i\in \mathbb{Z}$ represent the unique solution to $2i\equiv u\pmod{q}$. Replacing $\xi\in H^2(W;\mathbb{Z})$ in \eqref{eq:generator_xi} by $-\xi$, and hence $u$ by $-u$, if necessary, we may assume that $0\leq i\leq (q-1)/2$. Since $2u^2\equiv 1\pmod{q}$, we get $8i^2\equiv 1\pmod{q}$. Also, \eqref{eq:T4q_rank_one_inequality} gives \[
\delta_R(Y_{2,q},\mathfrak{s}_i,\tau_{2,q})\leq\frac{1-q}{16q}<0.
\]

If $i\leq k$, then \eqref{eq:T4q_mod8_piecewise} implies $i$ is odd, and \[
\frac{i^2}{2q}-\frac{1}{2}=\delta_R(Y_{2,q},\mathfrak{s}_i,\tau_{2,q})\leq\frac{1-q}{16q}
\]
gives $8i^2\leq 7q+1$. Thus, $(8i^2-1)/q$ is a positive integer with $\equiv 7\pmod{8}$ and $\leq 7$, so it must be $7$, i.e., $8i^2=7q+1$. Since $q>1$, we have $i>1$. Since $i$ is odd and $i^2\equiv 8i^2=7q+1\equiv 1\pmod{7}$, we get $i\geq13$, and therefore \[
k=\frac{q-1}{4}=\frac{(8i^2-1)/7-1}{4}=\frac{2i^2-2}{7} >3i.
\]
Now, applying \eqref{eq:T4q_rank_one_inequality} with $a=3$ gives \[
\frac{9i^2-q}{2q}=\delta_R(Y_{2,q},\mathfrak{s}_{3i},\tau_{2,q})\leq\frac{9-q}{16q},
\]
which gives $72i^2\leq7q+9$. This contradicts $8i^2=7q+1$.

If $i>k$, put $t=q-2i$, which is a positive odd integer. Since $t\equiv -2i\pmod{q}$, we have $2t^2\equiv 8i^2\equiv 1\pmod{q}$ as above. Also, \[
\frac{t^2-q}{8q}=\delta_R(Y_{2,q},\mathfrak{s}_i,\tau_{2,q})\leq\frac{1-q}{16q}
\]
gives $2t^2\leq q+1$, and therefore $q=2t^2-1$. Since $q>1$, we have $t\geq 3$. Now, $j:=3i-q$ satisfies \[
j=3i-q=3\cdot \frac{q-t}{2}-q =\frac{q-3t}{2}\leq \frac{q-1}{2}
\]
and \[
j-k=\frac{q-3t}{2}-\frac{q-1}{4}=\frac{q-6t+1}{4}=\frac{(2t^2-1)-6t+1}{4}=\frac{t(t-3)}{2}\geq 0,
\]
i.e., $k\leq j\leq (q-1)/2$. Thus, applying \eqref{eq:T4q_rank_one_inequality} with $a=3$ gives \[
\frac{9t^2-q}{8q}=\delta_R(Y_{2,q},\mathfrak{s}_{j},\tau_{2,q})\leq\frac{9-q}{16q},
\]
which gives $18t^2\leq q+9=2t^2+8$, a contradiction. This completes the proof. 
\end{proof}

\begin{corollary}
Let $q\geq 1$ be an odd integer. Then $\gamma_4(T_{4,q})=2$ in either of the following cases:
\begin{enumerate}[label=\textup{(\alph*)}, leftmargin=2em]
\item $q\equiv1\pmod{8}$ and $q$ is not represented by $x^2+32y^2$ with $x,y\in\mathbb{Z}$.
\item $q\equiv3\pmod{8}$ and $q$ is not represented by $x^2+2y^2$ with $x,y\in\mathbb{Z}$.
\end{enumerate}
\end{corollary}
\begin{proof}
    Since $\gamma_4(T_{4,q})\leq 2$ for every odd integer $q\geq 1$ \cite{BKST-2025}, the result follows from Theorem \ref{thm:T4q_nonorientable_genus} and Corollary \ref{cor:torus_linking_form}.
\end{proof}

\begin{corollary}\label{cor:T4q_infinite_family}
For every integer $r\geq0$, we have \[
\gamma_4(T_{4,\,17\cdot25^r})=2.
\]
\end{corollary}
\begin{proof}
Clearly $17\cdot 25^r\equiv 1\pmod{8}$. We prove by induction on $r$ that $17\cdot 25^r$ is not of the form $x^2+32y^2$. This is clear if $r=0$, so assume $r>0$. If $17\cdot 25^r=x^2+32y^2$, then reducing modulo $5$ gives $x^2+2y^2\equiv 0\pmod{5}$, and therefore $5$ divides both $x$ and $y$ since $-2$ is not a quadratic residue modulo $5$. Writing $x=5x_0$ and $y=5y_0$ gives $17\cdot 25^{r-1} = x_0^2+32y_0^2$, which is impossible by the induction hypothesis. 
\end{proof}
However, we have $17\cdot 25^r=(3\cdot 5^r)^2+8(5^r)^2$, so the condition of Corollary \ref{cor:torus_linking_form}(a) is satisfied. 

% \clearpage
\bibliography{references}{}
\bibliographystyle{alpha}
\end{document}